%% file: dp.tex
\documentclass[11pt]{amsart}
\usepackage{amsopn}
\usepackage{amssymb, amscd, verbatim}

\usepackage{graphics, epsfig}

\newcommand{\nc}{\newcommand}

\nc{\noi}{\noindent}
\nc{\ft}{\footnote}
\nc{\G}{\Gamma}
\nc{\g}{\gamma}
\nc{\Ld}{\Lambda}
\nc{\ld}{\lambda}
\nc{\al}{\alpha}
\nc{\be}{\beta}
\nc{\te}{\theta}
\nc{\ep}{\epsilon}

\nc{\la}{\langle}
\nc{\ra}{\rangle}
\nc{\ba}{\backslash}
\nc{\ke}{\hspace{-.3cm}}
\nc{\msl}{\{\!\{}
\nc{\msr}{\}\!\}}
\nc{\rlh}{\rightleftharpoons}
\nc{\fts}{\footnotesize}
\nc{\scs}{\scriptsize}

\nc{\gen}[2]{\langle #1,#2\rangle}
\nc{\im}{\mathtt{i}}

\nc{\mb}{\mathbf}
\nc{\m}{\mathcal}

\nc{\w}{{\mathbf w}}
\nc{\x}{{\mathbf x}}
\nc{\y}{{\mathbf y}}
\nc{\z}{{\mathbf z}}
\nc{\vv}{{\mathbf v}}
\nc{\bb}{{\mathbf b}}
\nc{\M}{{\mathbf M}}

\nc{\SO}{{\mathrm SO}}
\nc{\Spe}{{\mathrm Sp}}
\nc{\Sl}{{\mathrm SL}}
\nc{\SU}{{\mathrm SU}}
\nc{\Or}{{\mathrm O}}
\nc{\U}{{\mathrm U}}
\nc{\Gl}{{\mathrm GL}}
\nc{\PGl}{{\mathrm PGL}}
\nc{\Se}{{\mathrm S}}
\nc{\Cl}{{\mathrm Cl}}
\nc{\Spin}{{\mathrm Spin}}
\nc{\Pin}{{\mathrm Pin}}
\nc{\tr}{{\mathrm Tr}}

\nc{\R}{{\mathbb R}}
\nc{\HH}{{\mathbb H}}
\nc{\C}{{\mathbb C}}
\nc{\Z}{{\mathbb Z}}
\nc{\F}{{\mathbb F}}
\nc{\N}{{\mathbb N}}
\nc{\Q}{{\mathbb Q}}
\nc{\PP}{{\mathbb P}}

\nc{\rank}{\operatorname{rank}}

\newtheorem*{proposition}{Proposition}

\title{Describing the Platycosms}

\author[J. H. Conway \and J. P. Rossetti]{J. H. Conway$\,^\dag$ \and J. P. Rossetti$\,^\ddag$}

\address{Department of Mathematics, Princeton University, Fine Hall, Princeton, NJ 08544, USA.}
\email{conway@math.princeton.edu}
\address{FaMAF(Ciem), Universidad Nacional de C\'ordoba, Ciudad Universitaria, 5000-C\'ordoba, Argentina.
(Visiting Princeton University).}
\email{rossetti@mate.uncor.edu, rossetti@math.princeton.edu}

\thanks{$\dag$ NSF grant DMS-0072839.\qquad $\ddag$ Supported by a Guggenheim fellowship.}

\begin{document}

\maketitle

\begin{abstract} We study in detail the closed flat Riemannian 3-manifolds.
\end{abstract}

{\small
\setcounter{tocdepth}{1}
\tableofcontents
}

\section{Introduction}

That there are just 10 closed flat 3-manifolds --- we propose to call them {\it platycosms} --- has been known
since around 1933. Since then, they have been studied in many papers and several books, for example \cite{Th},
\cite{We}, \cite{Wo}.
They are of interest to speculative physicists as well as mathematicians, and indeed there are some
recent astronomical observations that suggest that the physical universe might actually be a platycosm!
(see for example \cite{E}, \cite{HS}, \cite{NYT}, \cite{Spe}, \cite{TOH}, \cite{URLW}).

We believe that these 10 manifolds should be more widely known. Accordingly, our principal aim in
this paper is to give a complete discussion of the ten compact platycosms.
In particular, we give them a uniform set of individual names, specify parameters in a systematic way, and
give presentations for their fundamental groups, which we use to derive their homology groups, automorphism
groups and list their double covers.

We then study their geometry, finding their Bravais types, diameters and injectivity radii, and listing
the compact embedded flat surfaces.

We have also included brief treatments of some related topics, namely the 2-dimensional analogues of the
platycosms, the~8 infinite platycosms and the 219 crystallographic groups.

Appendix I gives a brief proof that the list of 10 is complete, Appendix II describes the ``conorms'' we use
as parameters, while Appendix III is a dictionary between our names and others.

An accompanying paper \cite{hearing} proves that, up to scale, there is a unique isospectral pair of platycosms,
the `DDT-example' of \cite{DR}.

\noindent {\it Acknowledgements:}
We have benefitted from conversations with Jeff Weeks, Bill Thurston, and most particularly
with Peter Doyle, who has also helped us in many other ways.
J. P. Rossetti is grateful to the Mathematics Department of Princeton University for its hospitality during
the writing of this paper.

\section{Dimension 2}

Although our main aim in this paper is to describe the platycosms or 3-dimensional flat manifolds, particularly,
the compact ones, we briefly discuss their 2-dimensional analogues in this section.

\subsection*{The 2-dimensional flat manifolds}

\begin{figure}[!htb]
\input{torklein.pstex_t}
\caption{The torus $T_{A\,B\,C}$ and the Klein bottle $K^A_B$.}
\label{figtorklein}
\end{figure}
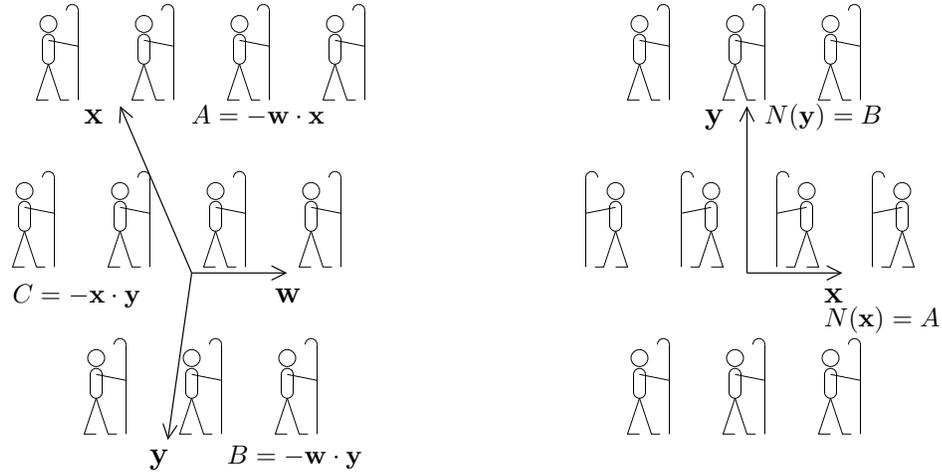

Everybody is familiar with the torus and Klein bottle. Each of these can be
given various locally Euclidean metrics --- that is to say, can be realized as a {\it flat} 2-manifold.
We can see this by ``rolling up'' the Euclidean plane in various ways:
one obtains the typical flat torus $T_{A\,B\,C}$  by dividing the plane by a 2-dimensional lattice
of translations (Figure~\ref{figtorklein}) --- this is mathematically more natural than identifying opposite sides of
a parallelogram, since each lattice can be defined by many different parallelograms.
Beware: the parameters $A,B,C$ for the torus are the negatives of the inner products of $\w,\x,\y$,
rather than their norms, which are $A+B,A+C,B+C$.

The most general flat Klein bottle $K^A_B$ (Figure~\ref{figtorklein}) can be obtained in a similar way
by dividing the plane by the group generated by the translation through a vector $\x$ with norm $N(\x)=\x\cdot\x=A$
together with a glide reflection based on an orthogonal vector $\y$ of norm $B$.

\subsection*{Voronoi and Bravais types}

The metric of a flat torus is determined by the shape of its translation lattice. Voronoi classified
lattices by the topological type of their Voronoi cell, which is either hexagonal or rectangular,
while Bravais classified them by their symmetries, which are controlled by the shape of the Delaunay cells.
The lattice parameters $A\,B\,C$ we recommend, called {\it conorms}
(see Figure \ref{figtorklein} and Appendix II), easily yield both classifications (Table~1),
since the Voronoi cell is rectangular just if some conorm vanishes, while the shape of the Delaunay cell is
determined by this and which conorms are equal.

\begin{table}[!htb]\label{tdim2}
\begin{tabular}{c|c|c|c|}
conorms & \begin{tabular}{c} topological type \\ of Voronoi cell \end{tabular} &
\begin{tabular}{c} shape of \\ Delaunay cell \end{tabular} & lattice shape \\
\hline
\begin{tabular}{c} $A\,B\,C$ \\ $A\,A\,B$ \\ $A\,A\,A$ \end{tabular} &
hexagonal &
\begin{tabular}{c} scalene triangle \\ isosceles triangle \\ equilateral triangle \end{tabular} &
\begin{tabular}{c} generic lattice \\ rhombic lattice  \\ hexagonal lattice \end{tabular} \\
\hline
\begin{tabular}{c} $A\,B\,0$ \\ $A\,A\,0$ \end{tabular} &
rectangular &
\begin{tabular}{c} rectangle \\ square \end{tabular} &
\begin{tabular}{c} rectangular lattice \\ square lattice \end{tabular} \\
\hline
\end{tabular}
\bigskip
\caption{(it is understood that $A,B,C$ are distinct and non-zero.)}
\end{table}

If one lattice can be continuously deformed into another without changing either its
Bravais or its Voronoi class, we say they are in the same `BraVo' class\ft{This is not the same
as saying that they are in the same Bravais class and also the same Voronoi class. The rhombic
lattices $\Ld_{A\,A\,B}$ split into two BraVo classes according as $A>B$ or $A<B$, between
which we cannot pass without encountering a hexagonal lattice $\Ld_{A\,A\,A}$.}

All these classifications can be applied to arbitrary flat manifolds and we discuss them
for the platycosms in Section \ref{sbravo}. All flat Klein bottles $K^A_B$ lie in a single
BraVo class since $A$ and $B$ can be continuously and independently varied through all positive
values.

\subsection*{Infinite flat 2-manifolds}

The torus and Klein bottle are finite 2-manifolds --- that is to say, they have finite volume, or
(equivalently for flat manifolds) are compact. There also exist just three types of
infinite (or non-compact) flat 2-manifolds without boundary, namely the

\begin{tabular}{rl}
{\it Euclidean Plane} & $\R^2$ ($\cong T_{\infty\,\infty}=K^\infty_\infty$), \\
{\it (infinite) Cylinder} & $C_A$ ($\cong T_{A\,\infty}\cong K^A_\infty$), \\
{\it M\"obius Cylinder} & $M_B$ ($\cong K^\infty_B$),
\end{tabular}

\begin{figure}[!htb]
\centerline{\mbox{\includegraphics*[scale=.5]{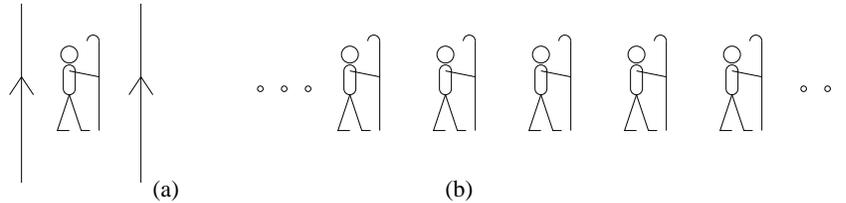}}}
\caption{The (infinite) cylinder (a) and its covering plane (b).}
\label{figinfcil}
\end{figure}

\begin{figure}[!htb]
\centerline{\mbox{\includegraphics*[scale=.5]{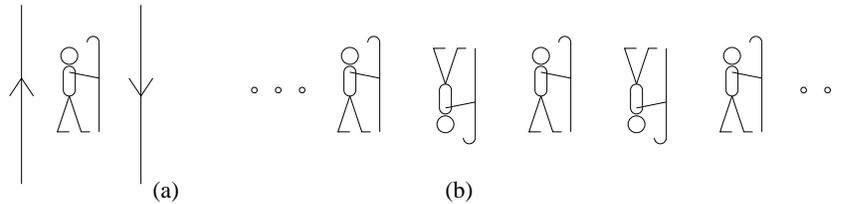}}}
\caption{The M\"obius cylinder (a) and its covering plane (b).}
\label{figinfklein}
\end{figure}

\

\noindent illustrated in Figures~\ref{figinfcil} and~\ref{figinfklein} .
Their 3-dimensional analogues are discussed in Section~\ref{sn-c}.

\subsection*{2-dimensional flat orbifolds}

Dividing a Euclidean space by a discrete subgroup of its symmetries yields a manifold only when
the group is {\it fixed-point-free} (i.e., only the identity element fixes any point of the space),
or equivalently, {\it torsion-free} (i.e., only the identity element has finite order).
More generally, one obtains a {\it flat orbifold}. There are just 17 types of compact flat 2-dimensional
orbifolds, corresponding to the 17 plane crystallographic groups, which are called
$$*632,\, 632,\, *442,\, 4{*}2,\, 442,\, *333,\, 3{*}3,\, 333$$
$$*2222,\,2{*}22,\, 22*,\, 22\times,\,2222,\,**,\,*\times,\,\times\times,\,\circ$$
in the orbifold notation \cite{orbnot}, the last two being the Klein bottle and the torus.
The analogous results in 3 dimensions are discussed in Section~\ref{spict}.

\section{Picturing the Platycosms - Space Groups}\label{spict}

\begin{figure}[!htb]
\centerline{\mbox{
\includegraphics*[scale=.5]{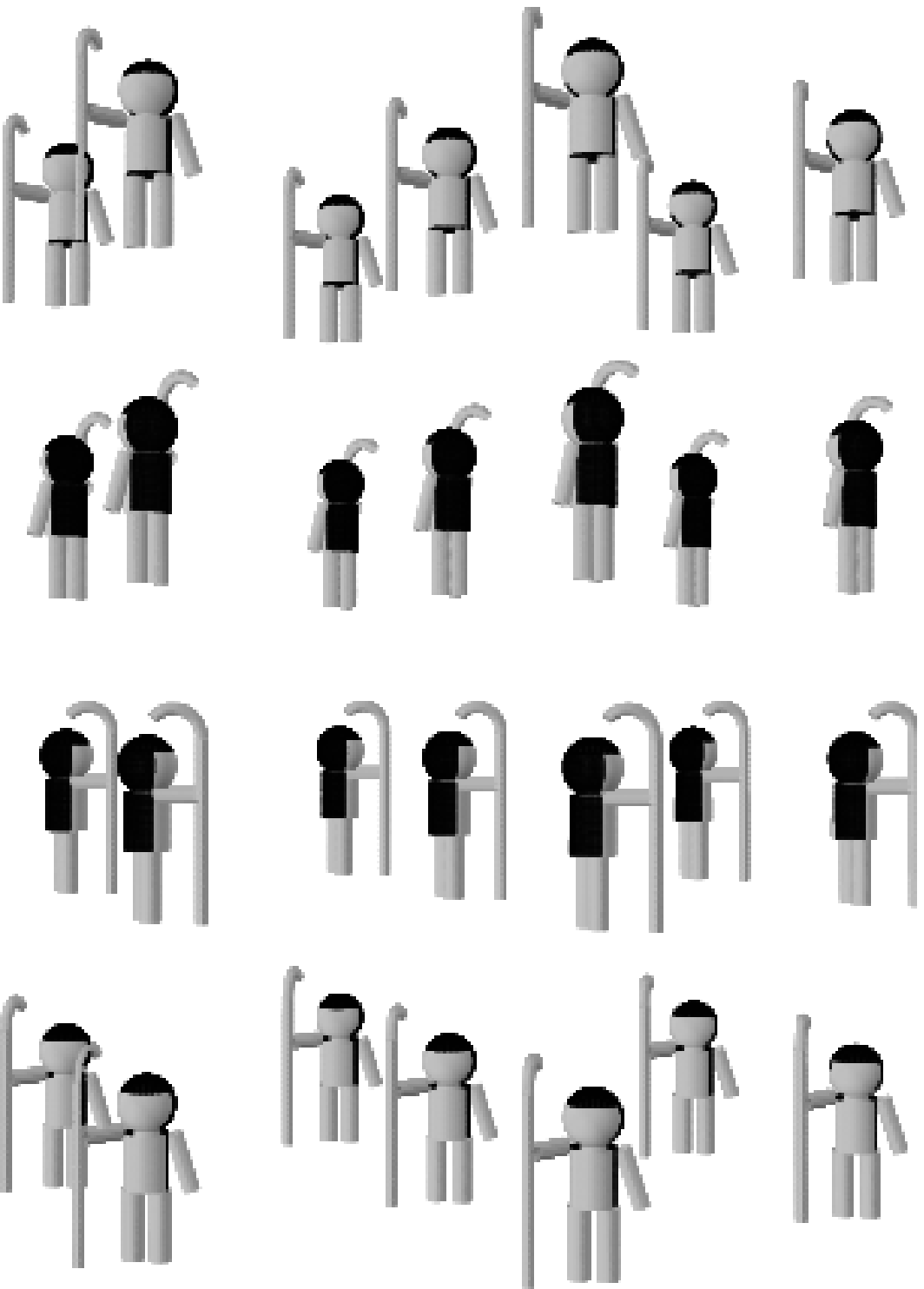} (i)  \hspace{2cm} (ii)
\includegraphics*[scale=.5]{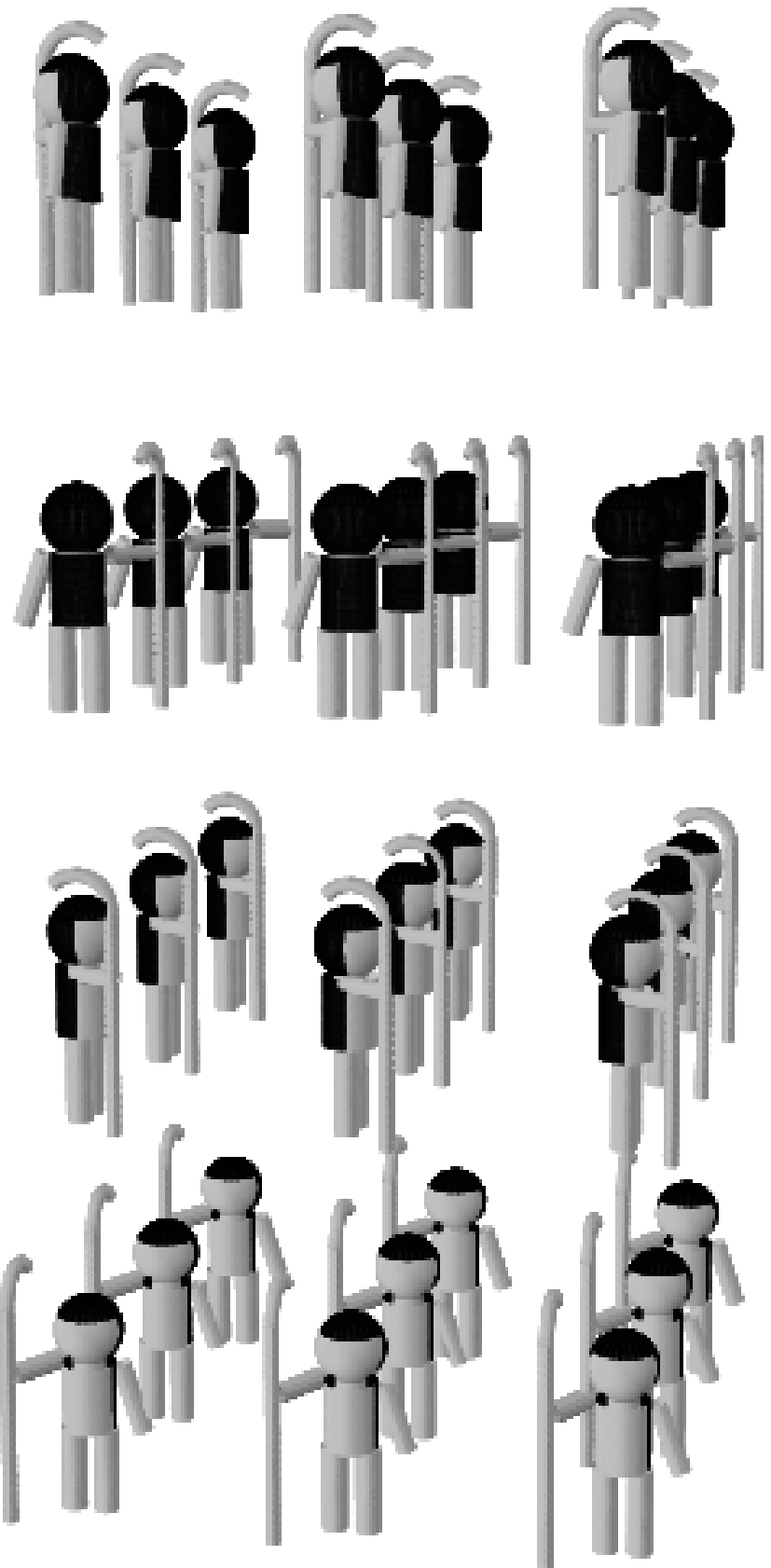} }}
\caption{(i) a tricosm, and (ii) a tetracosm.}
\label{figtritetra}
\end{figure}

\begin{figure}[!htb]
\centerline{\mbox{\includegraphics*[scale=.5]{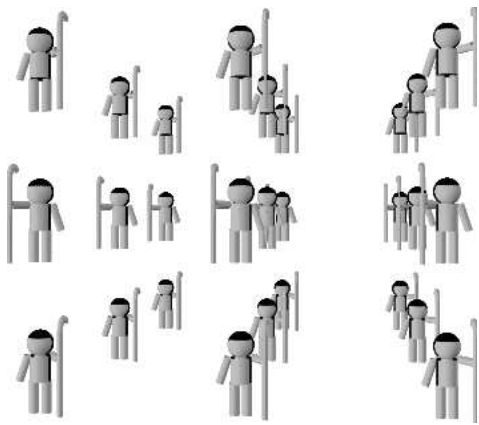}}}
\caption{A first amphicosm $+a1$.}
\label{fig+a1}
\end{figure}

\begin{figure}[!htb]
\centerline{\mbox{\includegraphics*[scale=.5]{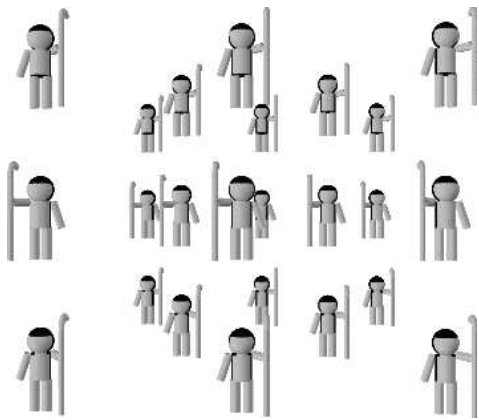}}}
\caption{A second amphicosm $-a1$.}
\label{fig-a1}
\end{figure}

We shall use the term {\it platycosm} (``flat universe'') for a compact locally Euclidean 3-manifold
without boundary, since these are the simplest alternative universes for us to think of living
in \cite{sciamer}. (The program {\it Curved Spaces} of Jeff Weeks --- accompanying \cite{We2} --- allows you
to `fly' through platycosms and other spaces.)
If you lived in a small enough platycosm, you would appear to be surrounded
by images of yourself which can be arranged in one of ten essentially different ways
(Figures \ref{figtoro}, \ref{figdicosm}, \ref{figtritetra}(i), \ref{figtritetra}(ii), \ref{figsinisdextral},
\ref{figdidi}, \ref{fig+a1}, \ref{fig-a1}, \ref{fig+a2-a2}(i), \ref{fig+a2-a2}(ii)).

The images one sees of oneself lie in the manifold's universal cover, which is
of course a Euclidean 3-space $\R^3$.  The symmetry operators that relate them
form a crystallographic space-group $\G$, and geometrically the manifold is the
quotient space $\R^3/\G$.

If you hold something in one hand, your images will either all hold things with the same hand
(if the manifold is orientable, Figures \ref{figtritetra}, \ref{figtoro},
\ref{figsinisdextral}, \ref{figdidi})
or half with their
left hand and half with their right (if not; Figures \ref{fig+a1}, \ref{fig-a1}, \ref{fig+a2-a2}).
We therefore call a platycosm {\it chiral} (``handed'') or {\it amphichiral}
(``either handed'') according as it is or is not orientable.
Selecting the images of a particular handedness in the picture for a non-orientable manifold $X$ yields
the picture for its orientable double cover, $Y$, so we regard $X$ as an `amphichiralized' version of $Y$,
and call it `an amphi-$Y$'.


The chiral platycosms are the {\it helicosms} $c1,c2,c3,c4,c6$ (individually called the
{\it torocosm, dicosm, tricosm, tetracosm, hexacosm}) and the {\it didicosm} $c22$,
also known as the {\it Hantzsche-Wendt manifold}.
In these notations, the letter ``$c$'' stands for ``chiral'' while the digits indicate the
point group, which is a cyclic group $C_N$ of order $N$ for $cN$, and $C_2\times C_2$ for $c22$.

\begin{figure}[!htb]
\centerline{\mbox{\includegraphics*[scale=.3]{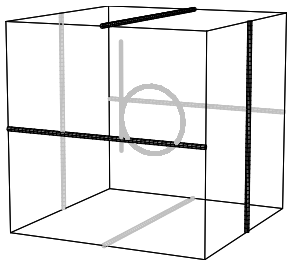}}}
\caption{The didicosm space group is generated by the half-turn screw motions corresponding to
the indicated lines.}
\label{figbox}
\end{figure}

A {\it torocosm} (usually called a 3-torus) is just the 3-dimensional analogue of a torus.
The space group for another helicosm $cN$ is generated by a suitable 2-dimensional lattice
of translations together with an orthogonal screw motion of period\ft{We say a screw motion has
period $N$ if the lowest power of it that is a translation is the~$N^{\mathrm{th}}$.} $N$, while the
space group for a didicosm is generated by orthogonal period 2 screw motions whose axes bisect
the faces of a `box' as in Figure~\ref{figbox}. The space is tessellated into such boxes.

\begin{figure}[!htb]
\centerline{\mbox{\includegraphics*[scale=.8]{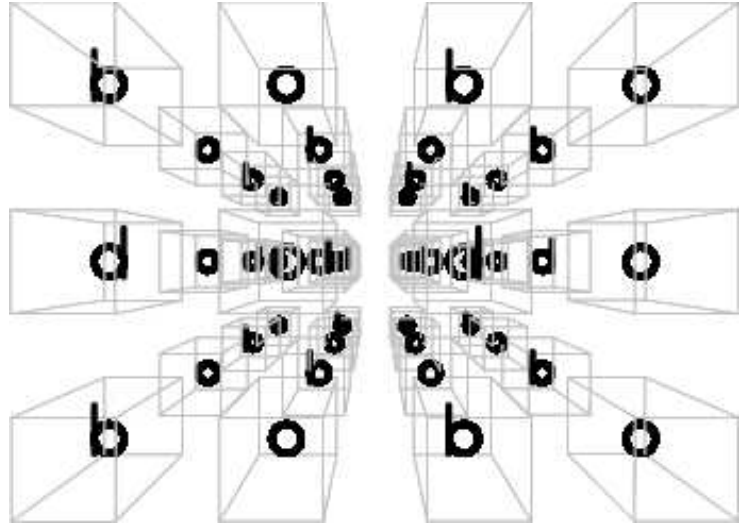}}}
\caption{Another picture of a second amphicosm.}
\label{figB}
\end{figure}

The amphichiral platycosms are the {\it first} and {\it second} (or {\it positive} and {\it negative})
{\it amphicosms}, $\pm a1$, and {\it amphidicosms}, $\pm a2$, whose orientable double covers are $c1$ and $c2$
respectively (Figures \ref{fig+a1}, \ref{fig-a1}, \ref{fig+a2-a2}).
The figures become easier to follow if we replace the human bodies (of Fig.~\ref{fig-a1}, say)
by boxes bearing the letters {\bf b},{\bf d},{\bf p},{\bf q} and insert some labelled {\bf o} to
indicate spacing (as in Fig.~\ref{figB}).

\begin{figure}[!htb]
\centerline{\mbox{
$+a1$ \includegraphics*[scale=.3]{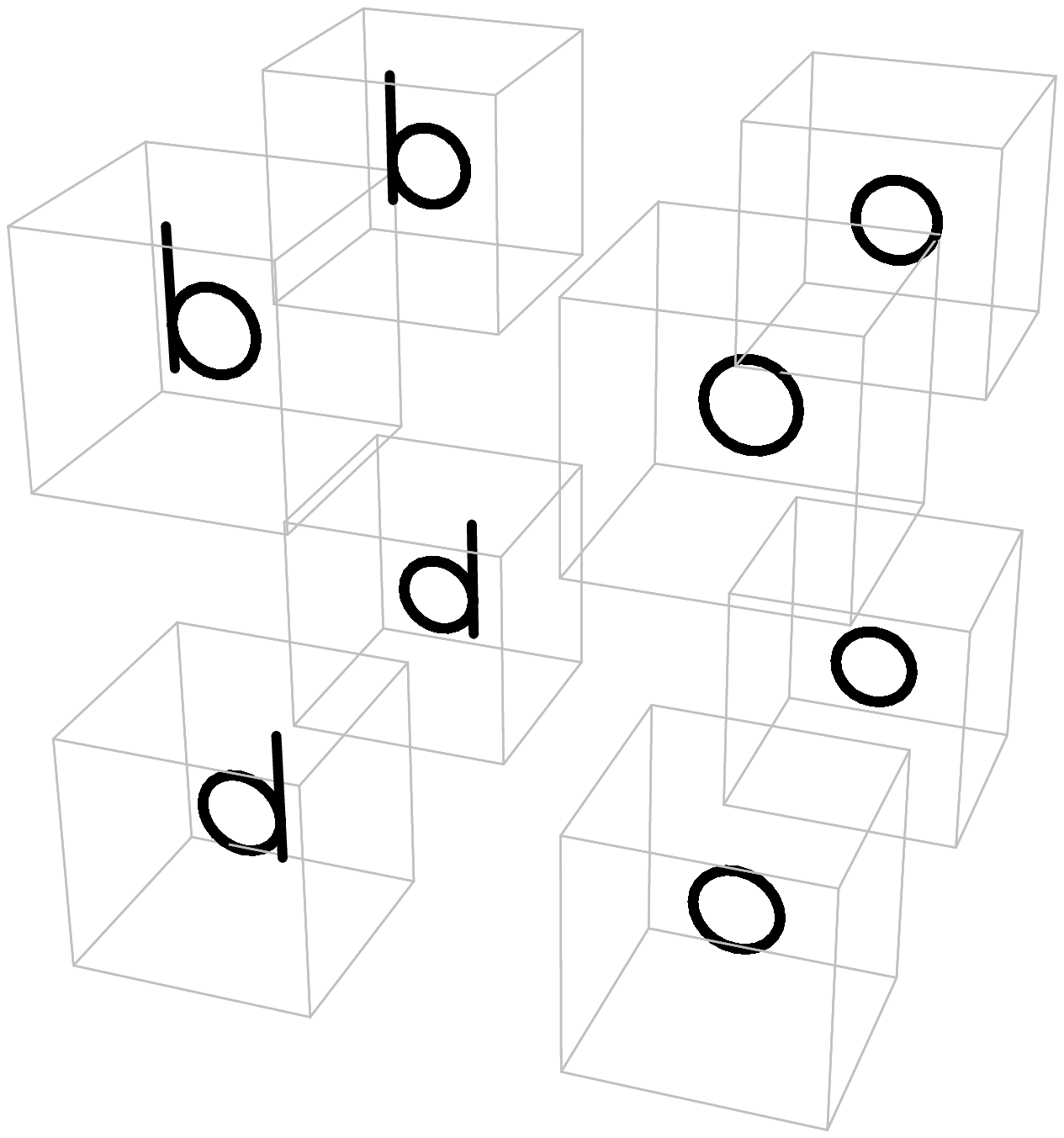} \qquad\qquad\qquad  $-a1$ \includegraphics*[scale=.3]{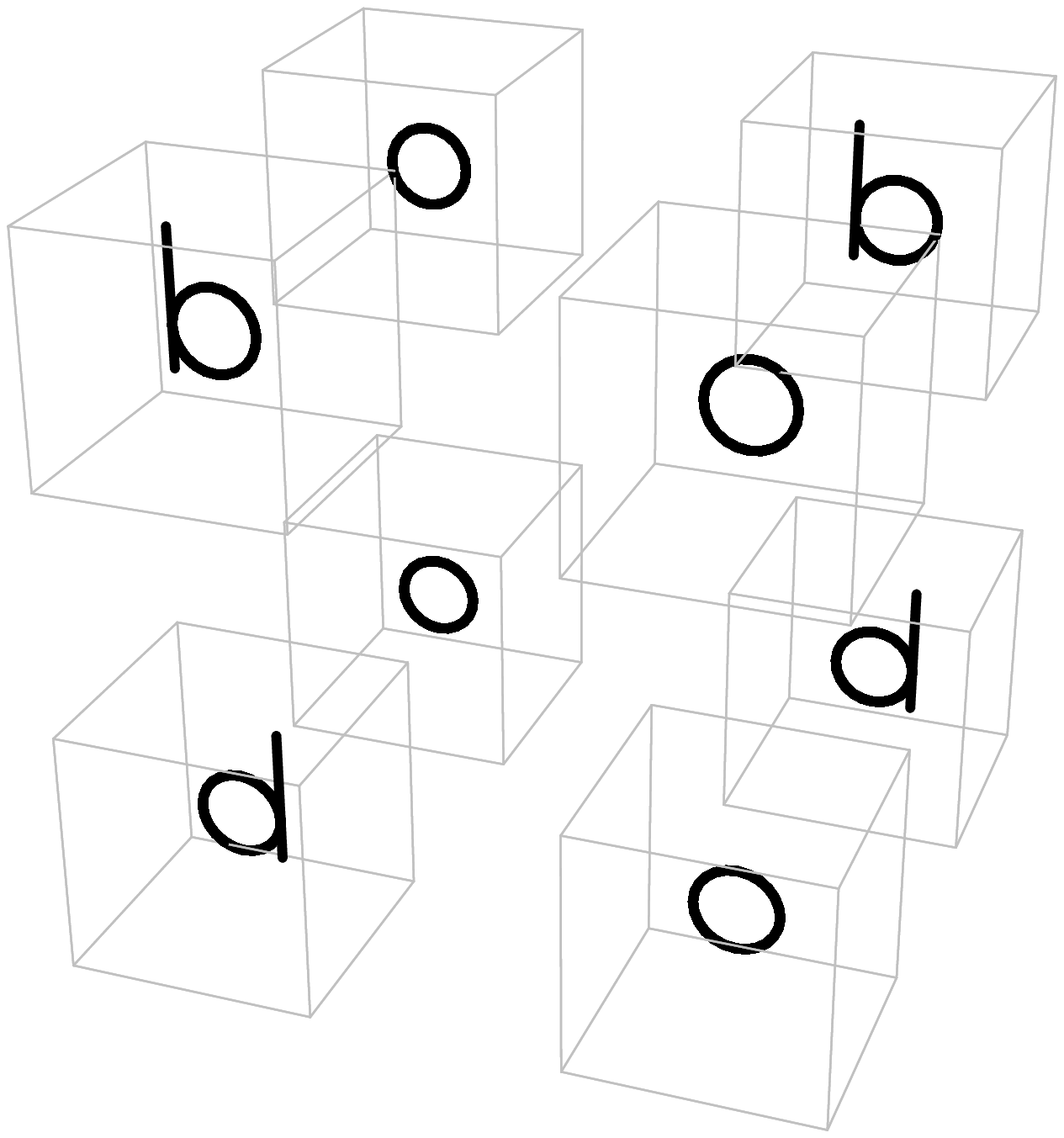}}}
\bigskip
\bigskip
\centerline{\mbox{$+a2$ \includegraphics*[scale=.3]{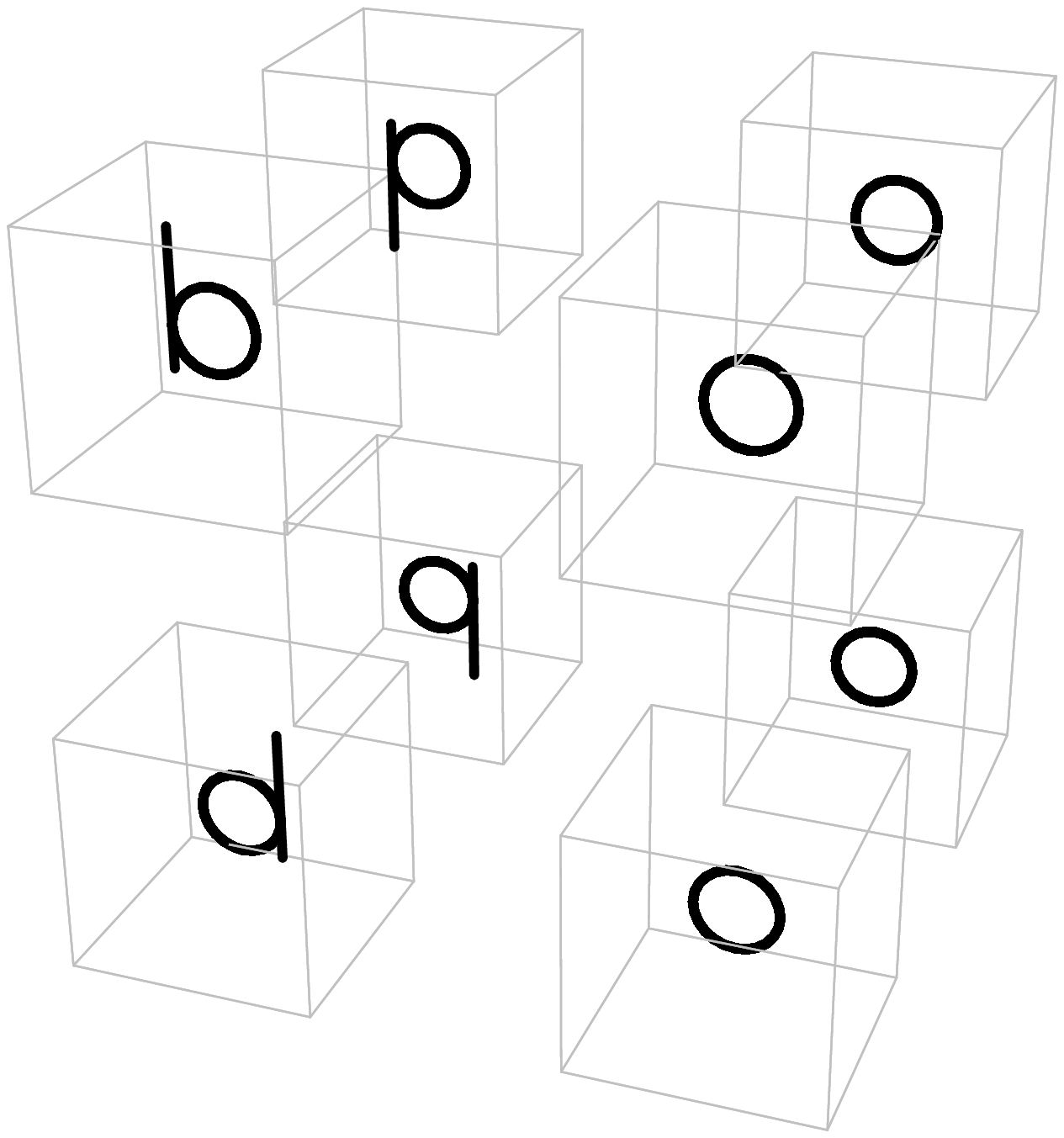}\qquad\qquad\qquad $-a2$ \includegraphics*[scale=.3]{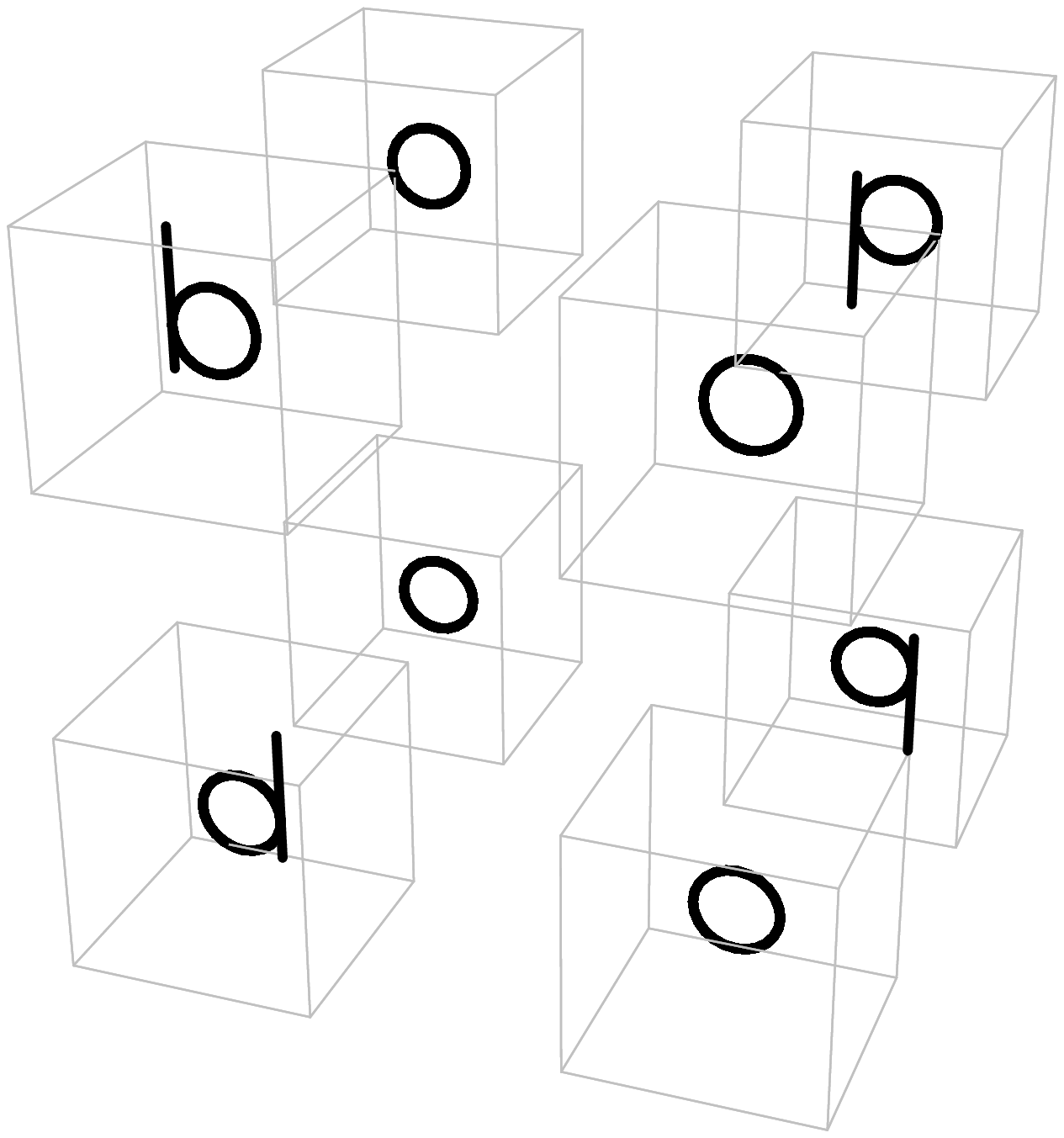}}}
\caption{The amphicosms $\pm a1$ and the amphidicosms $\pm a2$.}
\label{figfouramphis}
\end{figure}

Such shorter figures for the amphis (Fig.~\ref{figfouramphis}) abbreviate nicely to

\

\begin{equation}\label{bodo}
\begin{tabular}{ccccccc}
\begin{tabular}{|c c|} \hline \bf{b} & \bf{o} \\ \bf{d} & \bf{o} \\ \hline \end{tabular}
\raisebox{.5cm}[0pt][0pt]{\begin{tabular}{|c c|} \hline \bf{b} & \bf{o} \\ \bf{d} & \bf{o} \\ \hline \end{tabular}}
& &
\begin{tabular}{|c c|} \hline \bf{b} & \bf{o} \\ \bf{d} & \bf{o} \\ \hline \end{tabular}
\raisebox{.5cm}[0pt][0pt]{\begin{tabular}{|c c|} \hline \bf{o} & \bf{b} \\ \bf{o} & \bf{d} \\ \hline \end{tabular}}
& &
\begin{tabular}{|c c|} \hline \bf{b} & \bf{o} \\ \bf{d} & \bf{o} \\ \hline \end{tabular}
\raisebox{.5cm}[0pt][0pt]{\begin{tabular}{|c c|} \hline \bf{p} & \bf{o} \\ \bf{q} & \bf{o} \\ \hline \end{tabular}}
& &
\begin{tabular}{|c c|} \hline \bf{b} & \bf{o} \\ \bf{d} & \bf{o} \\ \hline \end{tabular}
\raisebox{.5cm}[0pt][0pt]{\begin{tabular}{|c c|} \hline \bf{o} & \bf{p} \\ \bf{o} & \bf{q} \\ \hline \end{tabular}}
\\
$+a1$ &  & $-a1$ &  & $+a2$ & & $-a2$
\end{tabular}
\end{equation}

\

\noindent in which it is understood that the letters continue with period two in all directions.
Their space groups consist precisely of the operations that take the leading box
labelled `{\bf b}' to the other boxes labelled {\bf b}, {\bf d}, {\bf p} or {\bf q}.

We chose these standard forms because when $x$ and $y$ are orthogonal the first amphicosm is the Cartesian
product of a Klein bottle and circle (in later notation: $+a1^D_{A:B}=K^D_A\times S^1_B$).
However, they can be continuously deformed into the {\it variant forms}:

\

\begin{equation}\label{bodo2}
\begin{tabular}{ccccccc}
\begin{tabular}{|c c|} \hline \bf{b} & \bf{o} \\ \bf{d} & \bf{o} \\ \hline \end{tabular}
\raisebox{.5cm}[0pt][0pt]{\begin{tabular}{|c c|} \hline \bf{d} & \bf{o} \\ \bf{b} & \bf{o} \\ \hline \end{tabular}}
& &
\begin{tabular}{|c c|} \hline \bf{b} & \bf{o} \\ \bf{d} & \bf{o} \\ \hline \end{tabular}
\raisebox{.5cm}[0pt][0pt]{\begin{tabular}{|c c|} \hline \bf{o} & \bf{d} \\ \bf{o} & \bf{b} \\ \hline \end{tabular}}
& &
\begin{tabular}{|c c|} \hline \bf{b} & \bf{o} \\ \bf{d} & \bf{o} \\ \hline \end{tabular}
\raisebox{.5cm}[0pt][0pt]{\begin{tabular}{|c c|} \hline \bf{q} & \bf{o} \\ \bf{p} & \bf{o} \\ \hline \end{tabular}}
& &
\begin{tabular}{|c c|} \hline \bf{b} & \bf{o} \\ \bf{d} & \bf{o} \\ \hline \end{tabular}
\raisebox{.5cm}[0pt][0pt]{\begin{tabular}{|c c|} \hline \bf{o} & \bf{q} \\ \bf{o} & \bf{p} \\ \hline \end{tabular}}
\\
$+a1$ &  & $-a1$ &  & $+a2$ & & $-a2$
\end{tabular}
\end{equation}
in necessarily two different ways. For the amphicosms, one can use the `shear motion' that lifts each layer `one-half a unit'
more than the one in front of it (or, in later notation, just change base from $\x,\y,\z$ to $\x,\pm\w,\z$).
The amphidicosms need only the `shift of motif' that leaves the metric undisturbed, but moves each letter half
way to the one above (in the $\mb{b},\mb{d}$ layers) or below it (in the $\mb{p},\mb{q}$ ones).
The variant forms are more useful in determining fundamental groups and automorphisms.

A similar figure for the didicosm is

\begin{figure}[!htb]
\centerline{\mbox{\includegraphics*[scale=.8]{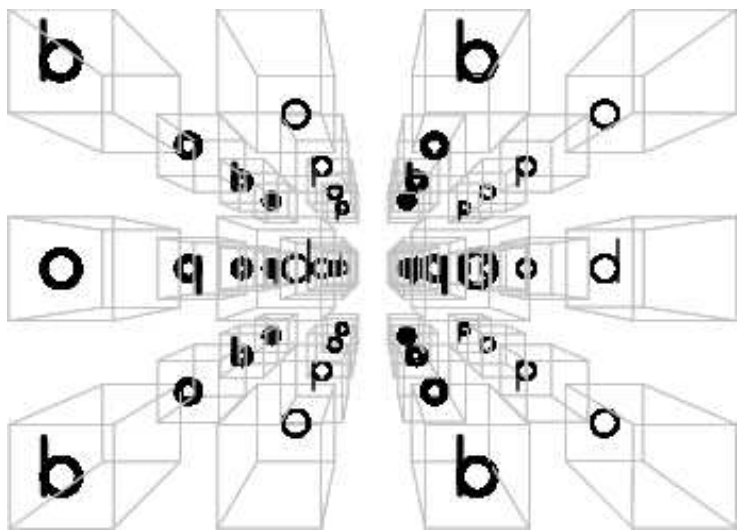}}}
\caption{A didicosm $c22$.}
\label{figD}
\end{figure}

\

\begin{equation}\label{bodoc22}
c22: \quad
\begin{tabular}{|c c|} \hline \bf{b} & {o} \\ \bf{o} & {d} \\ \hline \end{tabular} \,\,
\raisebox{.5cm}[0pt][0pt]{\begin{tabular}{|c c|} \hline \bf{o} & {p} \\ \bf{q} & {o} \\ \hline \end{tabular}}
\end{equation}
in which the non-bold letters are those seen from behind, through the boxes (Figure~\ref{figD}).

\subsection*{Platycosms and space groups}

Why are there only 10 platycosms?
The real reason is that just 10 of the crystallographic space groups are fixed-point-free,
analogously to the two plane crystallographic groups $\circ$ and $\times\times$ that
yield the torus and Klein bottle.
In fact, Nowacki \cite{No} found the 10 compact platycosms in 1934 by finding which of
the space groups of the International list were fixed-point-free.
A new and simple enumeration of the 219 space groups is given in \cite{CDHT}, from which
one can easily read off the numbers having various properties (and so, in particular,
pick out the platycosms). Of course one can obtain a much shorter enumeration by
restricting the argument to platycosms throughout, as was done by Hantzsche and Wendt \cite{HW} in 1934.
We give a simpler proof of this type in Appendix I.

\begin{table}[!htb]\label{tfibo}
\begin{tabular}{|c|ccccc|c|}
\hline
\# parameters & 1 & 2 & 3 & 4 & 6 & \# \\
\hline
\# fibrations & 0 & 1 & 3 & $\infty$ & $\infty$ & - \\
\hline
total & $35_0$ & $110_1$ & $12_1+26_2+21_3$ & $5_2+8_3$ & $2_1$ & 219 \\
\hline
chiral & $12_0$ & $29_1$ & $5_1+4_2$  & $2_2+1_3$ & $1_1$ & 54 \\
\hline
metachiral & 1 & 10 & 0 & 0 & 0 & 11 \\
\hline
\begin{tabular}{c} point groups \\ in this case\end{tabular} &
\begin{tabular}{c}$*432$, $432$ \\ $*332$, $3{*}2$ \\ 332 \end{tabular}
& \begin{tabular}{c} $*22N$, $22N$ \\ $*N\!N$, $N*$, $N\!N$ \\
$2{*}M$, $M\times$ \end{tabular} & $*222$, 222, $*22$ & $2*$, $22$, $*$  & $\times$, 1 & - \\
\hline
\# pt. groups & 5 & 19 & 3 & 3 & 2 & 32 \\
\hline
\multicolumn{3}{c}{} \\
\cline{3-7}
\multicolumn{2}{r|}{\# parameters} & 2 & 3 & 4 & 6 & \# \\
\multicolumn{2}{r|}{platycosms} & $c3$, $c4$, $c6$ & $c22$, $\pm a2$  & $c2$, $\pm a1$ & $c1$ & 10 \\
\multicolumn{2}{r|}{\# fibrations} & 1 & 3 & $\infty$ & $\infty$ & - \\
\multicolumn{2}{r|}{\text{total}} & $3_1$ & $1_1+2_3$ &  $1_2+2_3$ & $1_1$ & 10 \\
\multicolumn{2}{r|}{\text{chiral}} & $3_1$ & $1_1$ & $1_2$  & $1_1$  & 6  \\
\multicolumn{2}{r|}{metachiral} & 3 & 0 & 0 & 0 & 3 \\
\multicolumn{2}{r|}{point groups} & 33, 44, 66 & 222, $*22$ & 22, $*$ & 1 & \\
\multicolumn{2}{r|}{\# pt. groups} & 3 & 2 & 2 & 1 & 8 \\
\cline{3-7}
\end{tabular}
\bigskip
\caption{The upper table arranges space groups according to their point groups and
numbers of parameters and Seifert fibrations. The lower table restricts this information
to platycosms. The number $M$ is $2$ or $3$, and $N$ is $3$ or $4$ or $6$.}
\end{table}

We briefly describe the terms used here.
A space-group is {\it amphichiral} (``either handed'') or {\it achiral} (``not handed'')\ft{The
terms ``chiral'' and ``achiral'' were introduced by Lord Kelvin before
1896, when he used them in his Baltimore Lectures. Coxeter's enlargement of the latter to {\it amphichiral}
avoids using the negating prefix for a positive property.}
if it contains handedness-reversing operations; otherwise
it is {\it chiral} (``handed''). It is {\it metachiral} if the group itself is distinct from its mirror image.

The {\it point group} is the group obtained by identifying any two elements of the space group that
differ by a translation.
We specify these finite crystallographic groups in the orbifold notation \cite{orbnot}.
The space groups with a given point group constitute one of the 32 {\it crystal classes}.

The number of space-groups is often given as 230, because the metachiral ones are counted twice, once
for each of their two inequivalent orientations.  From the orbifold point of view this way of counting is
incorrect --- it would wrongly make us say there were 13 platycosms!  The numbers of types of oriented
orbifolds and platycosms are only  54+11=65  and 6+3 = 9  rather than 219 and 13.
The 9 oriented platycosm types arise as follows:
the cases $\pm a1,\pm a2$, are non-orientable, so yield no oriented types;
each of $c1,c2,c22$ has an orientation-reversing symmetry, so yields a single oriented type;
finally $c3,c4,c6$ give two oriented types each, since for them the defining screw motion may be either
{\it dextral} (like a corkscrew) or {\it sinistral} (like a reflected corkscrew).

\subsection*{Seifert fiber spaces}

The orbifolds of many space groups can be realized as Seifert fiber spaces, the one-dimensional fibers being the images
in the orbifold of a family of parallel lines that is fixed by the group.
The number of such `Seifert fibrations', if not 0 or 1, is at least 3 (in fact usually 3, otherwise $\infty$),
since if two families of parallel
lines are fixed so is the family of lines perpendicular to them\ft{and this third family will be rationally
related to the lattice if the first two are.}.

In Table~2, the numbers of cases with 3 or more fibrations appear with subscripts 1, 2 or 3 according
as the fibrations are of 1, 2 or 3 distinct types. Thus the number with just three fibrations is given
as $12_1+26_2+21_3$, meaning that there are
12 in which the three fibrations are all of 1 type, 26 in which the fibrations are of 2 distinct types,
and 21 in which the fibrations are of 3 distinct types.

It happens that every platycosm has at least one realization as a Seifert
fiber space. Each such fibration determines a plane crystallographic group
(by ``looking along the fibers'') and the ``types'' we now give (in Table~3) are the orbifold notations
for these groups:

\begin{table}[!htb]\label{tseifplat}
\begin{tabular}{|l|r|l|}
\hline
torocosm & $c1$      &  infinitely many fibrations, all of type $\circ$.  \\
dicosm & $c2$      &  one fibration of type $2222$, infinity of type $\times\times$.  \\
tricosm & $c3$      &  one fibration of type $333$.  \\
tetracosm & $c4$      &  one fibration of type $444$.  \\
hexacosm & $c6$      &  one fibration of type $632$.  \\
didicosm & $c22$      &  just three fibrations, all of type $22\times$.  \\
first amphicosm & $+a1$    &  one of type $\circ$, infinity of types $**$, $\times\times$.  \\
second amphicosm & $-a1$    &  one of type $\circ$, infinity of types $*\times$, $\times\times$.  \\
first amphidicosm & $+a2$  &  just three fibrations, of types $22*$, $**$, $\times\times$.  \\
second amphidicosm & $-a2$  &  just three fibrations, of types $22\times$, $*\times$, $\times\times$. \\
\hline
\end{tabular}
\bigskip
\caption{The platycosms and their Seifert fibrations.}
\end{table}

\section{Parameters for Platycosms}

In the interest of consistency, the parameters we use are always the conorms of what we call
the Naming lattice $\m{N}$, generated by all the translation vectors, screw vectors, glide vectors.
Since our parameters are often squared lengths, we introduce the convention that numbers inside a 
`square' should be squared, for instance $c22\framebox{$^{a\,b\,c}$}$ means $c22^{A\,B\,C}$, where $A=a^2$, $B=b^2$, $C=c^2$.

\subsection*{The Torocosm $c1=c1_{A\,B\,C}^{D\,E\,F}$}

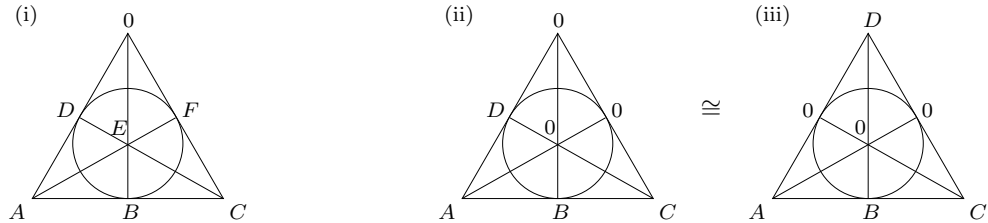
\begin{figure}[!ht]
\input{contor.pstex_t}
\caption{Conorms for $\Ld_{A\,B\,C}^{D\,E\,F}$ (i), and $\Ld_{A\,B\,C}^{D}$ (ii) and (iii).}
\label{figcontor}
\end{figure}

This is our proposed name for the 3-dimensional torus. Geometrically it is the quotient $\R^3/{\m T}$, where
${\m T}$ is the normal subgroup generated by translations --- ``the lattice of translations''.
Our parameters $A,B,\dots,F$ are the conorms (see Appendix II) of $\m{T}$  --- in other words $\m{T}$ has an obtuse
superbase $\mb{x},\mb{y},\mb{z},\mb{t}$ with $\mb{y}\cdot\mb{z}=-A$, $\mb{z}\cdot\mb{x}=-B$,
$\mb{x}\cdot\mb{y}=-C$, $\mb{x}\cdot\mb{t}=-D$, $\mb{y}\cdot\mb{t}=-E$, $\mb{z}\cdot\mb{t}=-F$.
In the corresponding conorm diagram ((i) in Figure \ref{figcontor})
the Voronoi vector corresponding to any line is a translation vector.

\begin{figure}[!htb]
\centerline{\mbox{\includegraphics*[scale=.5]{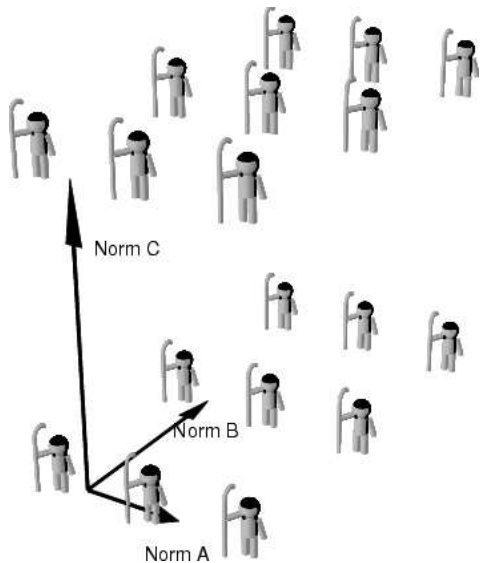}}}
\caption{The orthogonal torocosm $c1^{A\,B\,C}$.}
\label{figtoro}
\end{figure}

\begin{figure}[!htb]
\centerline{\mbox{\includegraphics*[scale=.5]{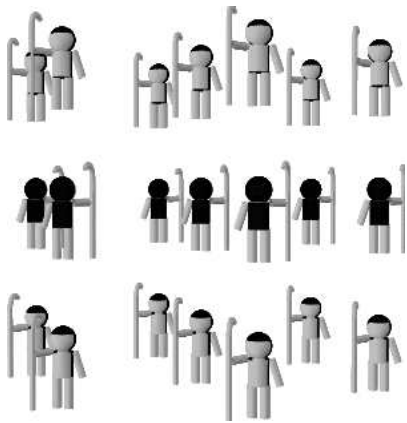}}}
\caption{The dicosm.}
\label{figdicosm}
\end{figure}

We omit conorms with value $0$: for instance the lattice for $c1^D_{A\,B\,C}$ has conorms (ii), or
equivalently (iii) in Figure~\ref{figcontor}. The lattice for $c1^{A\,B\,C}$ is spanned by three orthogonal
vectors of norms $A,B,C$ as in Figure~\ref{figtoro}.

The torocosm is the helicosm for $N=1$. The $N$-cosms for $N=2,3,4,6$ (the other helicosms) are generated
by the translations of a 2-dimensional lattice $\langle \mb{w},\mb{x},\mb{y}\rangle$
where $\mb{w}+\mb{x}+\mb{y}=\mb{o}$, whose conorms are the lower parameters together with a period $N$ screw
motion along a perpendicular vector $\mb{z}$ of norm $D$. The particular cases are:

\subsection*{The Dicosm $c2^D_{A\,B\,C}$}

The 2-dimensional lattice $\Ld_{A\,B\,C}$ has an obtuse superbase of three vectors $\mb{w},\mb{x},\mb{y}$ with
$\mb{w}+\mb{x}+\mb{y}=\mb{o}$ whose inner products $\mb{w}\cdot \mb{x}=-A$, $\mb{w}\cdot \mb{y}=-B$, $\mb{x}\cdot \mb{y}=-C$
are non-positive.
The half turn negating all these vectors is an order 2 rotational symmetry. The space group for the
dicosm $c2^D_{A\,B\,C}$ is generated by the translations of this lattice together with the period 2 screw motion
obtained by combining the above rotation with a translation through a vector $\mb{z}$ of norm $D$ perpendicular
to $\Ld_{A\,B\,C}$, so the conorm function for the naming
lattice $\Ld^D_{A\,B\,C}$ is (ii) or (iii) in Figure~\ref{figcontor}.

\subsection*{The Tricosm $c3^D_{A\,A\,A}$, Tetracosm $c4^D_{A\,A}$, Hexacosm $c6^D_{A\,A\,A}$}

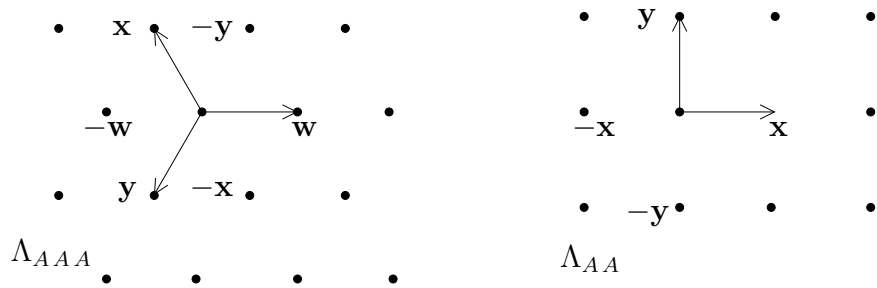
\begin{figure}[!ht]
\input{triang.pstex_t}
\caption{The hexagonal lattice $\Ld_{A\,A\,A}$ and square lattice $\Ld_{A\,A}$.}
\label{figtriang}
\end{figure}

Only certain special 2-dimensional lattices have higher order rotational symmetry. They are $\Ld_{A\,A\,A}$ which
has order 3 and 6 rotations taking $\w\to \x\to \y\to \w$ and $\w\to -\y\to \x\to -\w$
respectively (Figure~\ref{figtriang} left),
and $\Ld_{A\,A}$ (meaning $\Ld_{A\,A\,0}$), which has the order 4 rotation $\mb{x}\to \mb{y}\to-\mb{x}$
(Figure \ref{figtriang} right).

\begin{figure}[!htb]
\centerline{\mbox{\includegraphics*[scale=.6]{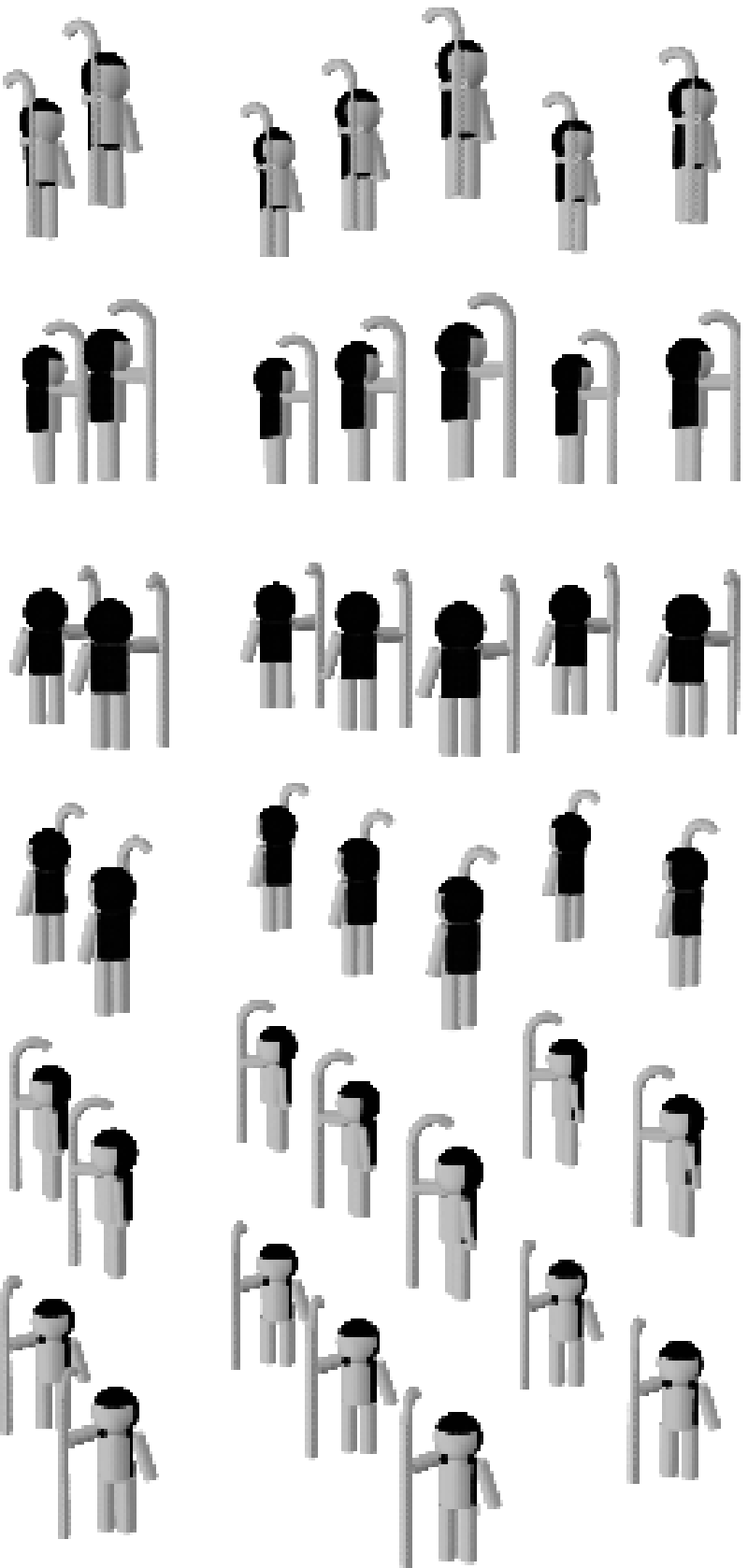} \hspace{2.5cm}
\includegraphics*[scale=.6]{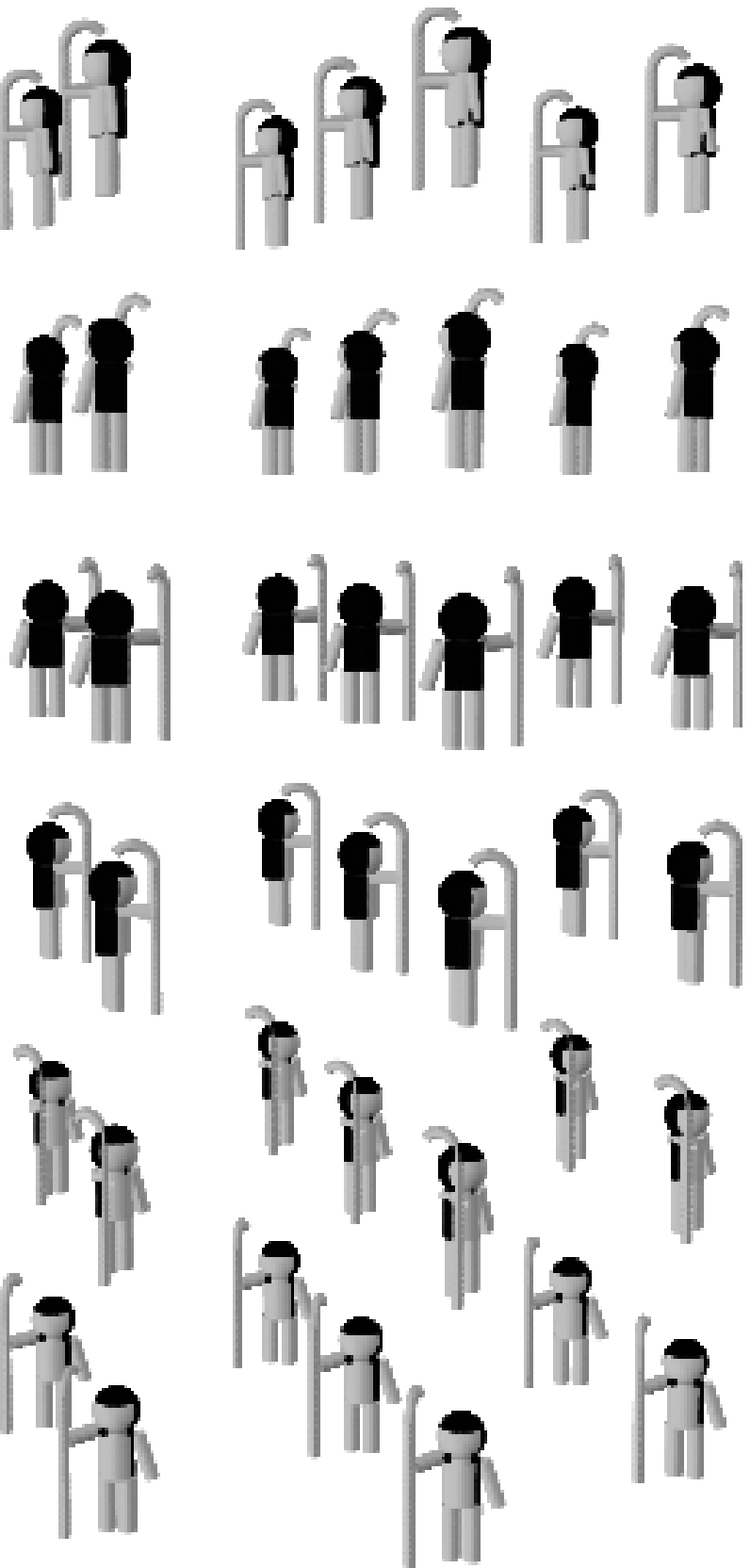}}}
\caption{A sinistral hexacosm (left) and a dextral hexacosm (right).}
\label{figsinisdextral}
\end{figure}

The space groups of the tricosm $c3^D_{A\,A\,A}$ (Figure~\ref{figtritetra}(i)), tetracosm $c4^D_{A\,A}$
(Figure~\ref{figtritetra}(ii)), hexacosm $c6^D_{A\,A\,A}$ (Figure~\ref{figsinisdextral})
are generated by the translations of the appropriate lattice together with the
period 3, 4 or 6 screw motion obtained by combining the corresponding rotation with a translation through a
perpendicular vector of norm $D$.

Since a period $N$ screw motion with $N>2$ is either {\it dextral} or {\it sinistral}
(see Figure~\ref{figsinisdextral}), the groups of these
three manifolds have two enantiomorphic forms each --- that is to say, they are metachiral.

\subsection*{The Didicosm $c22^{A\,B\,C}$}
\begin{figure}[!htb]
\centerline{\mbox{\includegraphics*[scale=.6]{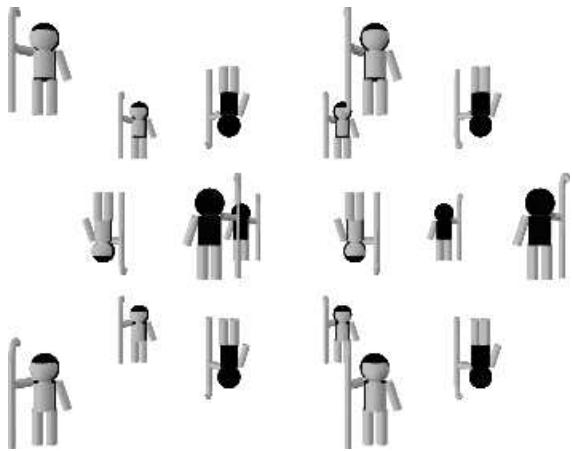}}}
\caption{The dizziness of the didicosm! (Figure~\ref{figD} may be more helpful.)}
\label{figdidi}
\end{figure}

For the didicosm (see Figures~\ref{figD} and~\ref{figdidi}), the naming lattice is generated by the three defining
screw vectors $\x,\y,\z$, and so is a rectangular lattice
$\Ld^{A\,B\,C}$, whose conorms $A,B,C$  are their squared lengths, say $a^2,b^2,c^2$.
The volume is $2abc$ (and so the {\it determinant}, or squared volume, is $4ABC$),
since a fundamental region consists of two $a\times b\times c$ `boxes', say the leading box
labelled {\bf b} in Figure~\ref{figD} and an adjacent one labelled {\bf o}.

\subsection*{The Amphicosms $\pm a1^{D}_{A:B\,C}$}

The naming lattice for the amphicosms $\pm a1^{D}_{A:B\,C}$ is the
lattice $\m{N}=\langle \w,\x,\y,\z\rangle$ (where $\w+\x+\y=0$) with Gram-matrix

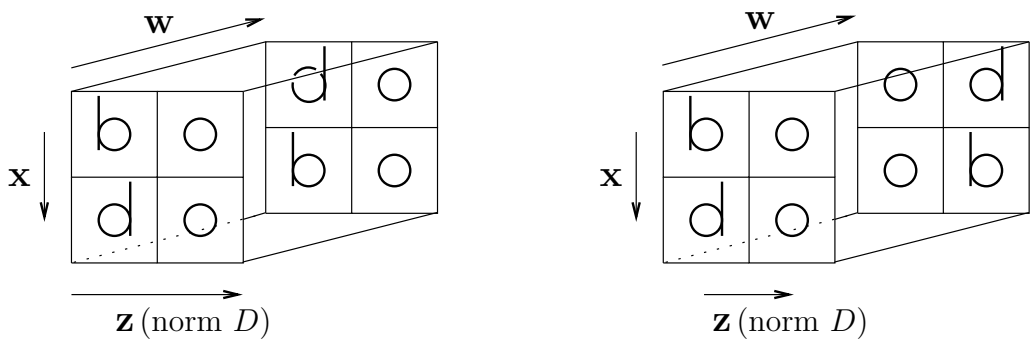
\begin{figure}[!htb]
\input{xfigamphic.pstex_t}
\caption{Defining vectors for the amphicosms.}
\label{xfigamphic}
\end{figure}

\begin{figure}[!htb]
\centerline{\mbox{\includegraphics*[scale=.45]{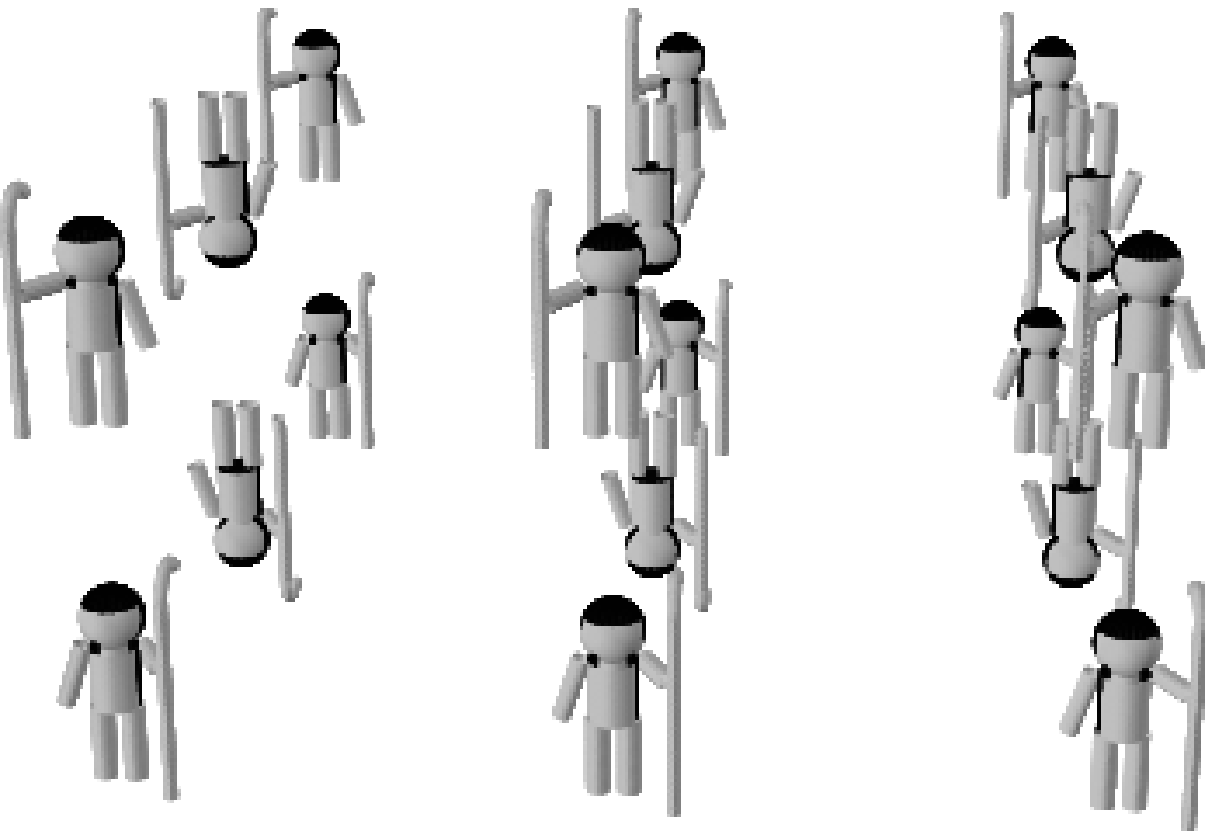} (i) \hspace{1.2cm} (ii)
\includegraphics*[scale=.5]{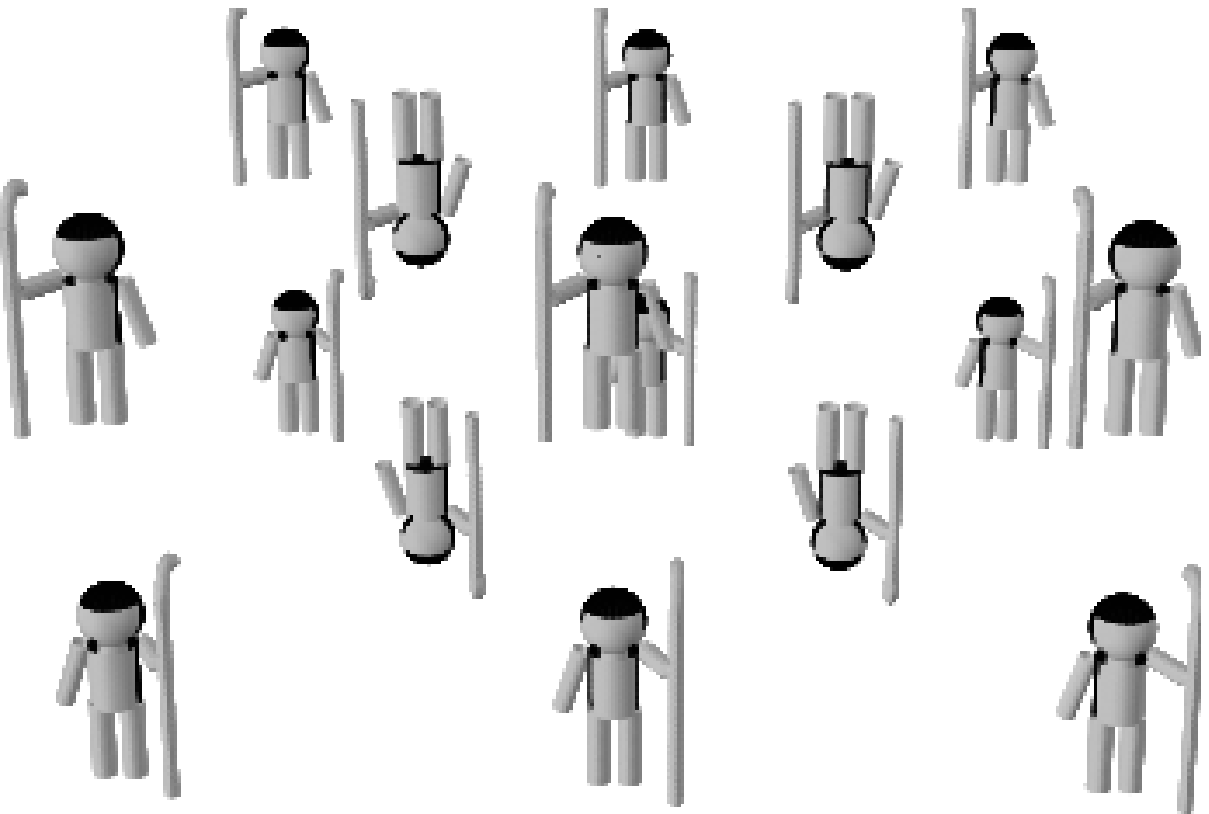}}}
\caption{(i) a first amphidicosm $+a2$, and (ii) a second amphidicosm $-a2$.}
\label{fig+a2-a2}
\end{figure}

{\small \begin{tabular}{c|cccc|}
& $\mb{w}$ & $\mb{x}$ & $\mb{y}$ & $\mb{z}$ \\
\hline
$\mb{w}$ & $A+B$ & $-A$ & $-B$ & $0$ \\
$\mb{x}$ & $-A$ & $A+C$ & $-C$ & $0$ \\
$\mb{y}$ & $-B$ & $-C$ & $B+C$ & $0$ \\
$\mb{z}$ & $0$ & $0$ & $0$ & $D$ \\
\hline \end{tabular}}\,,
but we must point out the `anomaly' that $\mb{z}$ plays different r\^oles in the two cases --- it is a translation vector
for $+a1$, but only half a translation vector for $-a1$:
In consequence, the determinants are $D(AB+AC+BC)$ for $+a1$ and $4D(AB+AC+BC)$ for $-a1$.
The translation lattices $\m{T}$ are $\langle 2\x,\y,\z\rangle$ and  $\langle 2\x,2\y,2\z,\y+\z\rangle$,
and $|\m{N}/\m{T}|=2$ or 4 respectively (see Figure \ref{xfigamphic}).

\subsection*{The Amphidicosms $\pm a2^D_{A:B}$}

\begin{figure}[!htb]
\input{xfigamphid.pstex_t}
\caption{Defining vectors for the amphidicosms.}
\label{xfigamphid}
\end{figure}
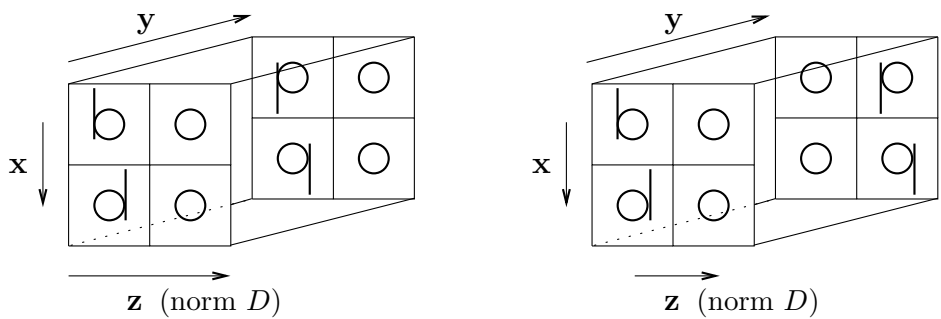

For the amphidicosms the Gram-matrix
{\small \begin{tabular}{c|ccc|}
& $\mb{x}$ & $\mb{y}$ & $\mb{z}$ \\
\hline
$\mb{x}$ & $A$ & $0$ & $0$ \\
$\mb{y}$ & $0$ & $B$ & $0$ \\
$\mb{z}$ & $0$ & $0$ & $D$ \\
\hline \end{tabular}}\,,
shows that the Naming lattice $\m{N}=\langle \mb{x},\mb{y},\mb{z}\rangle$ is a rectangular 3-dimensional
lattice $\Ld^D_{A\,B}\cong\Ld^{A\,B\,D}$.
We have the same anomaly about the meaning of $\mb{z}$, which leads to
differing determinants: $ABD$ for $+a1^D_{A:B}$ and
$4ABD$ for $-a1^D_{A:B}$.
The translation lattices $\m{T}$ are $\langle 2\x,2\y,\z\rangle$ and  $\langle 2\x,2\y,2\z\rangle$,
and $|\m{N}/\m{T}|=4$ or 8 respectively (see Figure \ref{xfigamphid}).


\subsection*{Which is which?}

It can be quite difficult to recognize which kind of amphicosm or amphidicosm one finds oneself in!
Thus it took 20 years before it was noticed that the four figures of amphichiral platycosms
in the celebrated book {\it The Shape of Space} \cite{We}, only define two inequivalent
manifolds!\ft{Figures 7.7 and 7.10 in that book are two first amphicosms with different metrics, while
Figures 7.11 and 7.12 are actually isometric first amphidicosms! To understand these equivalences, see
the remark on `variant forms' (\ref{bodo2}) in Section~\ref{spict}.}

There are several ways to answer this question.
One is that it is an {\it amphicosm} or {\it amphidicosm} according as its point group has size 2 or 4,
and {\it first} or {\it second} according as $|\m{N}/\m{T}|$ is 1 or 2 times this.
However, the safest rule we have found is

\smallskip

{\it An `amphi' is of first or second type if and only if it has a glide mirror
of first or second return, respectively.}

\smallskip

The reflecting plane of a glide reflection (a {\it glide mirror}) is of {\it first} or {\it second return} if and only if
the geodesic normal to its image at an arbitrary point closes exactly at its first or second return to the
surface respectively.

So to find which kind of amphichiral platycosm you are in, starting from a general point of some glide mirror,
walk perpendicularly to it until you first hit that mirror again. If this is always at the same point,
you are in a first (or positive) amphicosm or amphidicosm --- if always at a different point, in a second (or negative) one.
But beware: in an amphidicosm, there are also `ambiguous', glide mirrors that meet
most normals twice, but some only once. For more details, see the next section.

\section{Embedded Flat Surfaces}

In this section, ``embedded surface'' always means ``embedded compact flat 2-manifold''.
It will be either a torus or a Klein bottle, and either type may be embedded either
1-sidedly or 2-sidedly\ft{For surfaces embedded in Euclidean space, 2-sidedness is equivalent
to orientability, but in more general 3-manifolds this is no longer true.}.
Any such surface is the image of a plane in the universal cover with the property
that modulo $\Gamma$ it is compact.

Then the planes parallel to $\pi$ will have the same properties and so yield other embedded surfaces.
The configuration space of any such parallel family may be a circle, which we denote by $(2K)$ or $(2T)$
according as its elements (which are necessarily 2-sided) are Klein bottles or tori.
Alternatively, it may be an interval, in which case only the extreme members will be
1-sided; we will write $[1K\,\,(2T)\,\,1T]$ (say) for such an interval whose end surfaces
are a 1-sided Klein bottle and a 1-sided torus, and whose interior surfaces are 2-sided tori.

It is useful to classify 1-sided surfaces more closely, which we do by inserting
signs and the letters $g$ and $s$.
We write $+1$ for a first return (or `positive') surface,
meaning that the normal at any point $p$ intersects it only at $p$,
and $-1$ for a second return (or `negative') one,
when the normal at any $p$ intersects it in exactly two distinct points ($p$ and
the `negative' of $p$), and finally $\mp 1$ for an `ambiguous' one,
which intersects most normals twice but some only once.

Again if a surface is 1-sided then there is either a glide reflection or
a screw motion (or both) preserving it. We call it respectively a
{\it glide surface} or {\it screw surface} (or both), and indicate this by
the letters $g$ or $s$ (or $gs$).
So $+1gT$ refers to a positive (1-sided) glide torus,
$\mp1sK$ to an ambiguous screw Klein bottle.

Since embedded flat 2-manifolds are the images of planes in the universal cover, we must look for
such planes whose image modulo $\G$ is a torus or Klein bottle (the only compact possibilities).

Some considerations limit the search. The plane $\pi$ must not intersect any of its images under $\G$ other
than itself, which entails that all its images are parallel. Equivalently, $\G$ preserves the family of
parallel lines perpendicular to the plane.

There is a standard argument that usually proves $\pi$ is {\it basal}, that is, parallel to $\langle \mb{x},\mb{y}\rangle$,
or {\it perpendal}, that is, perpendicular to $\langle \mb{x},\mb{y}\rangle$. For if not, a unit
vector $\vv$ orthogonal to $\pi$ would resolve into basal and perpendal parts as $\vv=\mb{v_1}+\mb{v_2}$
with $\mb{v_1}\ne\mb{0},\mb{v_2}\ne\mb{0}$. Since the typical operation $\g\in\G$ takes this to
$\vv=\mb{v'_1}+\mb{v'_2}$ where $\mb{v'_2}=\pm\mb{v_2}$ we must have $\mb{v'_1}=\pm\mb{v_1}$ with the same sign.
However for every case other than $c1$ we can find a $\g$ for which this does not happen.

For the {\bf hexacosm} $c6$, {\bf tetracosm} $c4$ and {\bf tricosm} $c3$ the argument actually proves more, for
then the defining screw motion takes $\mb{v_1}$ to a vector $\mb{v'_1}$ which is not equal to $\pm\mb{v_1}$,
unless $\mb{v_1}=\mb{0}$, which forces $\vv=\mb{v_2}$ and $\pi$ to be basal. The basal
planes yield a circular family $(2T)$ of 2-sidedly embedded tori.

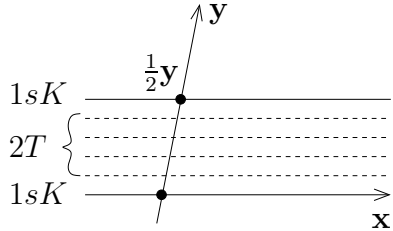
\begin{figure}[!ht]
\input{planesc2.pstex_t}
\caption{The extremes (through $0$ and $\frac 12 y$) are 1-sided Klein bottles while
all others are 2-sided tori.}
\label{figplanesc2}
\end{figure}

For the {\bf dicosm} $c2$ the same is true of the basal planes, but now there are also an infinite number
of families of perpendal planes, such as the family parallel to $\langle \mb{x},\mb{z}\rangle$.
These yield an interval $[1sK\,\,(2T)\,\,1sK]$ of embedded
surfaces as suggested in Figure \ref{figplanesc2}.
We obtain a similar interval of perpendal planes for each primitive vector $\mb{x'}$ (counted up to sign) of
the lattice $\langle \mb{x},\mb{y}\rangle$, since any such vector appears in some basis $\mb{x'},\mb{y'}$.
On the other hand, the standard argument shows that every embedded surface must be basal or perpendal,
since the defining screw motion takes $\mb{v_1}$ to $-\mb{v_1}$ but $\mb{v_2}$ to $+\mb{v_2}$.
This shows that the infinity of intervals $[1sK\,\,(2T)\,\,1sK]$ of the preceding paragraph
do indeed comprise all the non-basal embedded surfaces, agreeing with the entry
$(2T)^1$, $[1sK\,\,(2T)\,\,1sK]^{\infty}$ in Table~4.

For the {\bf torocosm} $c1$ the standard argument does not apply, but the answer is easy --- there is a
circular family $(2T)$ of
2-sided tori corresponding to each 2-dimensional section\ft{a {\it section} of a lattice is its full
intersection with some subspace.} of $\m{T}=\langle \x,\y,\z\rangle$,
or equivalently to each pair of primitive vectors $\pm\mb{v^*}$ of the dual lattice $\m{T}^*$, yielding
the entry $(2T)^{\infty}$ in Table~4.

\begin{figure}[!htb]
\centerline{\mbox{\includegraphics*[scale=.5]{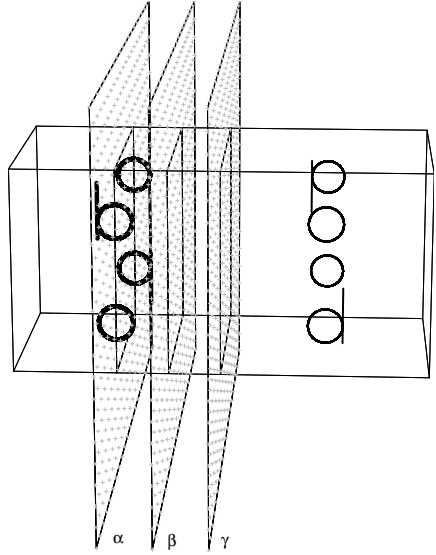}}}
\caption{The plane $\al$ is a Klein bottle, which is embedded 1-sidedly as shown
by the screw motion taking {\bf b} to {\bf q}, and the same is true for $\g$ as shown by
that taking {\bf b} to d. In view of these two screw motions, we need only
consider planes like $\be$ between $\al$ and $\g$  --- these are 2-sided tori.}
\label{figplanesc22}
\end{figure}

For the {\bf didicosm} $c22$ the standard argument can be applied in different directions, showing that $\pi$
must be parallel or perpendicular to $\langle \mb{y},\mb{z}\rangle$ as well as to $\langle \mb{x},\mb{y}\rangle$,
which forces it to be parallel to one of the three coordinate planes $\langle \mb{x},\mb{y}\rangle,\,
\langle \mb{x},\mb{z}\rangle,\, \langle \mb{y},\mb{z}\rangle$.
We obtain three intervals of type $[\mp1sK\,\,(2T)\,\,{\mp}1sK]$ (Figure \ref{figplanesc22}).

For the {\bf amphicosms} $\pm a1$ the standard argument applies since the glide reflection takes
$\mb{v_1}$ to $+\mb{v_1}$ but $\mb{v_2}$ to $-\mb{v_2}$. We first handle the perpendal planes.

We can see that the planes parallel to $\langle \mb{x},\mb{z}\rangle$
or $\langle \mb{y},\mb{z}\rangle$ yield circles of 2-sided Klein bottles and 2-sided tori respectively,
since $\mb{x}$ is a glide vector and $\mb{y}$ is a translation vector.\ft{In view of the symmetry
interchanging $\w$ and $\x$, the planes parallel to $\gen{\w}{\z}$ behave like those parallel to $\gen{\x}{\z}$.}
But any primitive vector of $\langle \mb{x},\mb{y}\rangle$ belongs to a superbase $\mb{w'},\mb{x'},\mb{y'}$ whose
vectors are respectively congruent modulo 2 to $\mb{w},\mb{x},\mb{y}$, which makes $\mb{w'},\mb{x'}$ be glide vectors
and $\mb{y'}$ a translation vector, and any such
superbase is equivalent to $\mb{w},\mb{x},\mb{y}$ by some isotopy of the lattice. So the perpendal planes yield
infinitely many circles of each type $(2K)$ or $(2T)$.
The basal planes form an interval of type $[+1gT\,\,(2T)\,\,{+}1gT]$ for $+a1$,
$[-1gT\,\,(2T)\,\,{-}1gT]$ for $-a1$.

\begin{table}[!htb]\label{tembsurf}
\begin{tabular}{c|rl}
platycosm & \multicolumn{2}{c}{families of surfaces} \\
\hline
$c1$ &  & $(2T)^\infty$ \\
$c2$ & $(2T)^1$; & $[1sK\,(2T)\,1sK]^\infty$ \\
$c3$ &  & $(2T)^1$ \\
$c4$ &  & $(2T)^1$ \\
$c6$ &  & $(2T)^1$ \\
$c22$ & & $[\mp1sK\,\,(2T)\,\,{\mp}1sK]^3$ \\
$+a1$ & $[+1gT\,\,(2T)\,\,{+}1gT]^1$; & $(2K)^\infty,\,\,\, (2T)^\infty$ \\
$-a1$ & $[-1gT\,\,(2T)\,\,{-}1gT]^1$; & $(2K)^\infty, \,\,\,(2T)^\infty$ \\
$+a2$ & $[+1gsK\,\,(2K)\,\,{+}1gsK]^1$; & $(2K)^1,\,\,\,\, [\mp1sK\,\,(2T)\,\,{\mp}1gT]^1$ \\
$-a2$ & $[-1gT\,\,(2T)\,\,{-}1sK]^1$; & $(2K)^1,\,\,\,\, [\mp1sK\,\,(2T)\,\,{\mp}1gT]^1$ \\
\end{tabular}
\bigskip
\caption{The parallel families of embedded surfaces, with the type symbol defined in the text
and the number of families of this type indicated by the superscript.
Those before a semicolon are images of basal planes, those after of perpendal ones.}
\end{table}

Finally, for the {\bf amphidicosms} $\pm a2$ the standard argument applies in both $\langle \mb{x},\mb{y}\rangle$
and $\langle \mb{x},\mb{z} \rangle$ directions, showing that $\pi$ must be parallel to a coordinate plane.
The perpendal planes parallel to
$\langle \mb{x},\mb{z}\rangle$ form a circle $(2K)$  and those parallel to $\langle \mb{y},\mb{z}\rangle$ an interval
$[\mp1sK\,\,(2T)\,\,{\mp}1gK]$ while the basal planes (parallel to $\gen{\mb{x}}{\mb{y}}$) form an
interval of type $[+1gsK\,\,(2K)\,\,{+}1gsK]$ for $+a2$ and $[-1gT\,\,(2T)\,\,{-}1sK]$ for $-a2$.

We summarize the results in Table~4.

\section{Infinite Platycosms}\label{sn-c}

Although elsewhere in this paper `platycosm' means `finite platycosm', in this section we briefly
discuss and name the infinite ones.
An {\it infinite platycosm} is a boundaryless flat manifold that has infinite volume, or equivalently
is not compact.

One kind of infinite platycosm, the `Product Space' or {\it Prospace}, is topologically the cartesian product of
a compact flat manifold of dimension $0,1,2$ by a Euclidean space of the complementary dimension $3,2,1$.
The cases are

\

\begin{tabular}{ccccl}
\begin{tabular}{c} compact \\ factor \end{tabular} & \begin{tabular}{c} complementing \\ (fiber) space \end{tabular} &
name & \begin{tabular}{c} notation and \\ parameters \end{tabular}\\
\hline
Point        & $\R^3$ & Euclidean Space     & $EUC$ \\
Circle       & $\R^2$ & Circular Prospace   & \quad $CPS_A(\theta)$ \\
Torus        & $\R^1$ & Toroidal Prospace   & \,\,\,\,\,\,$TPS_{A\,B\,C}$ \\
Klein Bottle & $\R^1$ & Kleinian  Prospace  & \,$KPS^A_B$ \\
\hline
\end{tabular}

\bigskip

The others are related to these in the way a M\"obius strip is to an annulus, so we call them
{\it M\"obius spaces}, or `Mospaces'.
Technically, they are fibrations whose base is a compact flat submanifold of dimension 1 or 2 and whose fiber
is the `complementary' Euclidean space of dimension 2 or 1, which is reflected when we traverse at
least one closed path in the base. This path must be homotopically non-trivial and may be orientation preserving ($+$)
or orientation-reversing ($-$).
The cases are

\

\begin{tabular}{ccccc}
\begin{tabular}{c} compact \\ base \end{tabular} & \begin{tabular}{c} compl. \\ fiber \end{tabular} &
\begin{tabular}{c} path \\ type \end{tabular} &  name & \begin{tabular}{c} notation and \\ parameters \end{tabular}\\
\hline
Circle & $\R^2$ & $+$ & Circular Mospace & $\,\,CMS_A$ \\[.05cm]
Torus & $\R^1$ & $+$ & Toroidal Mospace & $\quad\,\, TMS_{A\,B:C}$ \\[.05cm]
Klein Bottle & $\R^1$ & $-$ & chiral Kleinian Mospace & $+KMS^A_B$ \\[.05cm]
Klein Bottle & $\R^1$ & $+$ & achiral Kleinian Mospace & $-KMS^A_B$ \\
\hline
\end{tabular}

\bigskip

\begin{figure}[!htb]
\centerline{\mbox{\includegraphics*[scale=1]{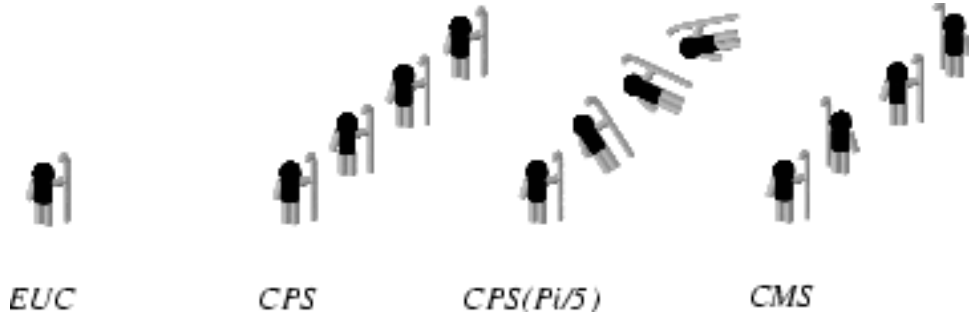}}}
\caption{One's images in the simplest infinite platycosms.}
\label{figcps}
\end{figure}

\begin{figure}[!htb]
\centerline{\mbox{\includegraphics*[scale=1]{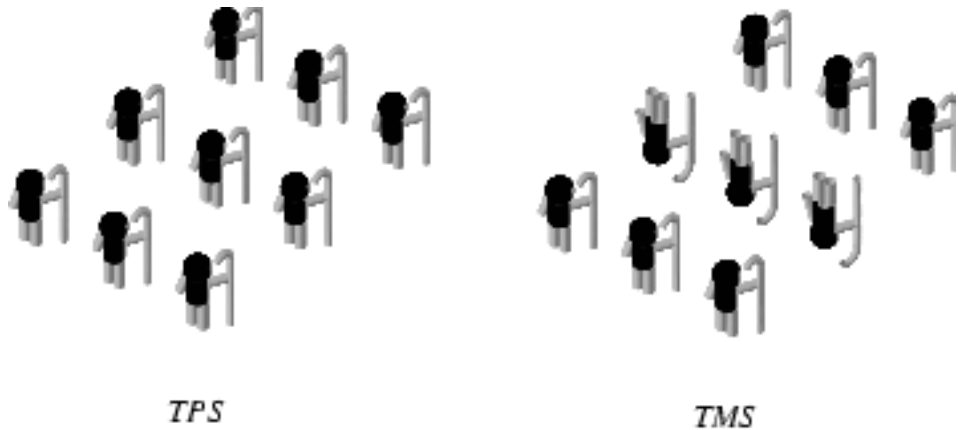}}}
\caption{The toroidal Product and M\"obius spaces.}
\label{figtps}
\end{figure}

\begin{figure}[!htb]
\centerline{\mbox{\includegraphics*[scale=1]{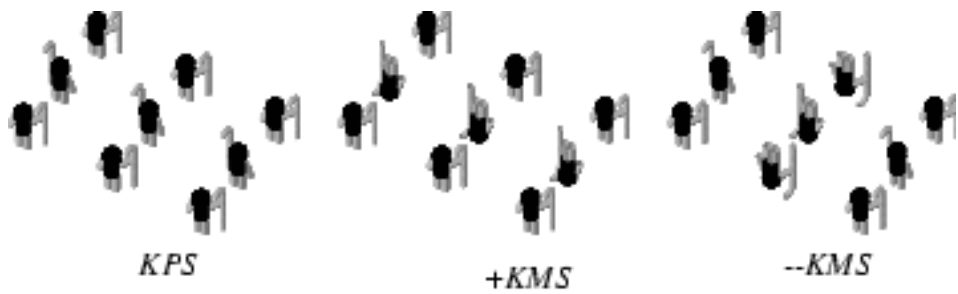}}}
\caption{Kleinian Product and M\"obius spaces.}
\label{figkps}
\end{figure}

Figures \ref{figcps}, \ref{figtps} and \ref{figkps} picture these manifolds (strictly, their universal covers).
The group for $CPS_A(\theta)$ is generated by a screw motion of angle $\theta$ and length $\sqrt{A}$.
In the {\it untwisted} case when $\theta=0$ this is the direct product of a circle and a
plane even in the metrical sense, in the twisted cases (when $\theta$ is not a multiple of $2\pi$) it is only topologically so.

In $CPS_A(\theta)$ ($=CPS\framebox{$_a$}(\theta)$) there is always a closed geodesic of length $a=\sqrt{A}$: this is the
only primitive closed geodesic unless $\frac {\theta}{2\pi}$ is rational, when the associated screw motion has
a finite period $N$, and all geodesics parallel to it are closed, of length $aN$.
For $CMS_A$ ($=CMS\framebox{$_a$}$) the group is generated by a glide reflection of length $a=\sqrt{A}$.
The manifolds described so far are illustrated in Figure~\ref{figcps}.

For $TPS_{A\,B\,C}$ the group is a 2-dimensional lattice of translations for which some superbase
$\mb{v_0},\mb{v_1},\mb{v_2}$ has conorms $A,B,C$.
In $TMS_{A\,B:C}$ the translations through $\mb{v_0}$ and $\mb{v_1}$ (so {\it not} $\mb{v_2}$) are replaced by the
glide reflections obtained by combining them with the reflection in the
base plane $\langle \mb{v_0},\mb{v_1},\mb{v_2}\rangle$.
These manifolds appear in Figure \ref{figtps}.

The group for $KPS^A_B$ ($=KPS\framebox{$^a _b$}$) is generated by a translation along $\mb{v_1}$ of length $a=\sqrt{A}$
and a glide reflection along a vector $\mb{v_2}$ of length $b=\sqrt{B}$ perpendicular to $\mb{v_1}$ whose
reflecting plane is orthogonal to the base plane $\langle \mb{v_1},\mb{v_2}\rangle$.
For $+KMS^A_B$ ($=+KMS\framebox{$^a _b$}$) we compose the latter with the reflection in the base plane, so that it becomes
a screw motion.
For $-KMS^A_B$ ($=-KMS\framebox{$^a_b$}$) we instead compose the former with that reflection, so making it a glide
reflection.
Figure \ref{figkps} shows these three Kleinian manifolds.

Beware: not every 3-dimensional product of flat manifolds is a `Prospace' --- we use that term
only for the product of a compact manifold and a Euclidean space.
The following little table shows that no fewer than 8 of the 18 platycosm types include non-trivial
direct products (Table~5).

\begin{table}[!htb]\label{tprodof}
\begin{tabular}{c|c|c|}
product of & line & circle \\
\hline
plane & $EUC$ & $CPS\,^*$ \\
cylinder & $CPS\,^*$ & $TPS\,^*$ \\
M\"obius cylinder & $CMS$ & $TMS\,^*$ \\
torus & $TPS$ & $c1\,^*$ \\
Klein bottle & $KPS$ & $+a1\,^*$ \\
\hline
\end{tabular}
\bigskip
\begin{caption} {An asterisk indicates that the metrical direct products of this kind do not yield all
parameter sets in this case.}
\end{caption}
\end{table}

\section{Fundamental Groups, Homology and Automorphisms}

The fundamental group of a platycosm is just the space group $\G$ we see in its universal cover.
We explain how to find a presentation, using the amphicosms and amphidicosms as our examples.
Let us mark any letter {\bf b}, {\bf d}, {\bf p}, {\bf q} in our diagrams with the element of $\G$ by which it can be
obtained from the
initial {\bf b}  --- this sets up a 1-1 correspondence, since these letters have no symmetry.\ft{Not even front-to-back!}
Then a presentation for $\G$ can be computed as follows.
First find generators for $\G$, and then generators for $\m{T}$ in terms of these.
Then a set of relations will be sufficient if and only if they imply the correct structures for

i) the subgroup $\m{T}$;

ii) the action of each generator on $\m{T}$;

iii) the point group $G=\G/\m{T}$.

\begin{figure}[!htb]
\input{xfigamphis.pstex_t}
\caption{Generators for the `amphis'.
In these, $X$ is always a glide reflection with vector $\x$ that takes the leading {\bf b} to the {\bf d} below it,
while $Z$ is the translation with vector $\z$ or $2\z$ that takes it to the next {\bf b} to its right.
For $+a1$, $W$ is a glide reflection with vector $\w$.
For $-a1$, it is a glide reflection with vector $\w+\z$ while the $W$'s
for $\pm a2$ are (different) screw motions, with vectors $\y$.}
\label{xfigamphis}
\end{figure}
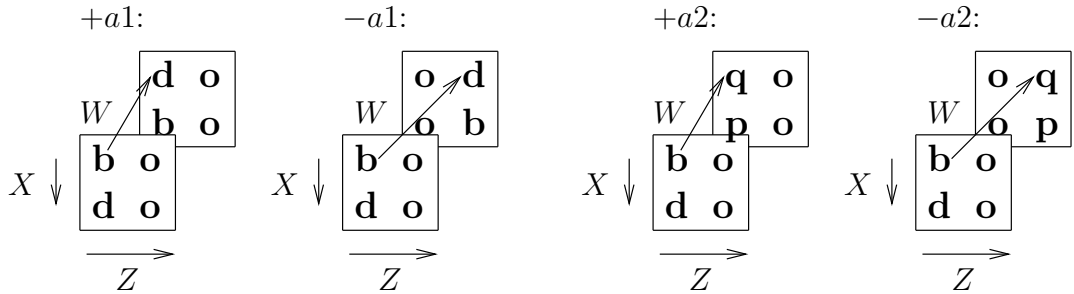

Since the two amphicosms $\pm a1$ are the hardest cases, we discuss them in detail.
It is easy to see that the operations $W,X,Z$ defined by Figures \ref{xfigamphic} and \ref{xfigamphis} are generators.
The vectors associated to these are $\mb{w},\mb{x}$ and $\mb{z'}=\mb{z}$ or $2\mb{z}$, and we see that the translation
lattice is generated in each case by $2\mb{w},2\mb{x},\mb{w}+\mb{x},\mb{z'}$, showing that the subgroup $\m{T}$ is
generated by $W^2,X^2,WX,Z$.

We obtain the correct structure for $\m{T}$ by demanding that these commute, and that
$(WX)^2=W^2X^2$ {\bf or} $W^2X^2Z$ in the two cases. Let us suppose in general that
we have found some generators $G_1,G_2\dots$ for $\G$, and certain products $w_1,w_2,\dots$
of them that generate $\m{T}$, and relations $R_1,R_2,\dots$ that define its structure.
Our next step is to define the correct action on $\m{T}$, for which it suffices
to express each $w_i^{G_j}=G_j^{-1}w_iG_j$ as a function of the $w_i$.
For the first amphicosm, we find that $W^2$ and $X^2$ are fixed by all three of
$W,X,Z$, while $Z$ is fixed by itself (obviously), but inverted by each of $W$ and $X$,
so we add these assertions as relations.

Finally, we must ensure that our relations imply the correct structure for the point group $G=\G/\m{T}$.
However, for the amphicosms they already do, since modulo $\m{T}$ we have $W^2\equiv X^2\equiv WX\equiv Z\equiv 1$,
which show that in the point group, $W$ and $X$ map to the same element of order two, while $Z$ maps to~$1$.

The relations one obtains in this way can almost always be simplified. The best way to
do this is to adjoin the shortest relations one can find and then delete any redundant ones.
For the first or second amphicosms, we find three simple relations, namely
$$Z^W=Z^X=Z^{-1} \quad\mathrm{and}\quad W^X=W\,\mb{or}\,ZW$$
in the two cases.

The fact that $W,X$ invert $Z$ implies that $W^2,X^2,WX$ commute with $Z$. But also, these three
relations imply that $W^2$ commutes with $X^2$  --- obviously for the first amphicosm, while for the negative one
$ZW^2Z^W=ZW^2Z^{-1}=W^2$.
Moreover, they imply that
$$XWXW=X^2W^XW=X^2(1\,\mathbf{or}\,Z)W^2,$$
which was the further relation needed for the structure of $\m{T}$.

Finally, our three simple relations specify the action of $W,X$ on $Z$ and each other, so in particular
on $W^2,X^2,WX,Z$, the generators of $\m{T}$. Rewriting the last relation as
$$[X,W]=1 \,\mathbf{or}\, Z,$$
we see that:

{\it The fundamental groups of the amphicosms have presentations:}
$$\begin{array}{ll}
\text{for } +a1\,\, & \langle W,X,Z : Z^X=Z^W=Z^{-1},[X,W]=1\rangle \\
\text{for } -a1\,\, & \langle W,X,Z : Z^X=Z^W=Z^{-1},[X,W]=Z\rangle.
\end{array}$$
Of course, their first homology groups are obtained by abelianizing these. They are
$$\begin{array}{cll}
\text{for } +a1 \,\, & \langle W,X,Z | \text{ \it abelian}, \,\, Z^2=1\rangle &\cong C_2\times C_\infty\times C_\infty \\
\text{for } -a1 \,\, & \langle W,X,Z | \text{ \it abelian}, \,\, Z=1\rangle  & \cong C_\infty\times C_\infty,
\end{array}$$
since in the $+a1$ case $Z$ was conjugate to $Z^{-1}$, so $Z^2$ maps to 1, while
for $-a1$ $Z$ was a commutator, so itself maps to 1.
In Table~6, we abbreviate these group structures to $2\cdot\infty^2$ and $\infty^2$.

\subsubsection*{Helicosms}

For the helicosms we have the presentations of the table, in which $Z$ is the
defining screw motion, while $X$
and $Y$ generate the 2-dimensional lattice perpendicular to it.

\subsubsection*{Didicosm}

For the didicosm $c22$ we find in this way that the translation subgroup is
generated by the squares of the
generating screw motions $X,Y,Z$, which invert the squares of each other, and
satisfy $XYZ=1$.
It turns out that this presentation reduces to $\langle X,Y\,|\,X=Y^2XY^2,
Y=X^2YX^2\rangle$ when we
define $Z=(XY)^{-1}$.
The abelianization is $C_4\times C_4$, which we denote by $4^2$.
This proves the well known fact that the Hantzsche-Wendt didicosm $c22$ is the only platycosm with finite homology,
or equivalently, with zero first Betti number.

\begin{table}[!htb]\label{tpres}
\begin{tabular}{cccc}
\smallskip
 & presentation & $H_1$  & transl. lattice \\
\hline
\bigskip
$c1$ & $X \rlh Y\rlh Z \rlh X$   &  ${\infty^3}^{\phantom N}$ & $X,Y,Z$  \\
\bigskip
$c2$ & $\begin{array}{c} X \rlh Y \\ Z:X\to X^{-1},Y\to Y^{-1} \end{array}$ & $2^2\cdot\infty$ & $X,Y,Z^2$  \\
\bigskip
$c3$ & $\begin{array}{c} X \rlh Y \\ Z:X\to Y\to (XY)^{-1} \end{array}$ & $3\cdot\infty$ & $X,Y,Z^3$ \\
\bigskip
$c4$ & $\begin{array}{c} X \rlh Y \\ Z:X\to Y\to X^{-1} \end{array}$ & $2\cdot\infty$ & $X,Y,Z^4$ \\
\bigskip
$c6$ & $\begin{array}{c} X \rlh Y \\ Z:X\to XY\to Y \end{array}$ & $\infty$ & $X,Y,Z^6$  \\
\bigskip
$c22$ & $\langle X,Y| X=Y^2XY^2, Y=X^2YX^2 \rangle$ & $4^2$ & $X^2,Y^2,Z^2=(XY)^{-2}$ \\
\bigskip
$+a1$ & $W,X:Z\to Z^{-1}, [X,W]=1$  & $2\cdot\infty^2$ & $W^2,X^2,WX,Z$  \\
\bigskip
$-a1$ & $W,X:Z\to Z^{-1}, [X,W]=Z$  & $\infty^2$ & $W^2,X^2,WX,Z$  \\
\bigskip
$+a2$ & $W,X:Z\to Z^{-1}, W:X\to X^{-1}$  & $2^2\cdot\infty$ & $W^2,X^2,Z$  \\
\bigskip
$-a2$ & $W,X:Z\to Z^{-1}, W:X\to X^{-1}Z$  & $4\cdot\infty$  & $W^2,X^2,Z$  \\
\end{tabular}
\caption{}
\end{table}

We remark that the fundamental groups of $c1,c2,+a1,+a2$ need three generators; each of
the others can be generated by two elements.

\subsection*{Automorphisms}

The largest group of `symmetries' of $\G$ is its {\it affine normalizer}\ft{so called because it is the normalizer of
$\G$ in the group of affine transformations of $\R^3$.}, which consists of the affine automorphisms of $\R^3$
that take $\G$ into itself. We discuss the four `parts' of this group\ft{Beware: these `parts' are not always subgroups
and not always disjoint.} illustrated in Fig.~\ref{figsymmetries}:

\begin{figure}[!htb]
\input{symmetries.pstex_t}
\caption{}
\label{figsymmetries}
\end{figure}

{\it Part I} is the {\it connected part} (the component of the identity), which
by continuity must fix $\G$ pointwise since the identity does. In other words, it is central.
For the platycosms, it consists only of translations and its dimension is the first Betti
number $\be_1$.
It is a subtle theorem (cf.\ \cite{ChV}) that factoring out part I yields precisely the automorphisms 
of $\G$.\ft{Our argument, although it refers to the way an automorphism acts on $\G$, does not in fact use this theorem.}

{\it Part II}, the {\it inner part}, is the space group $\G$, which acts trivially on $M$, but (of course) induces 
inner automorphisms on~$\G$. Since these induce the identity on $M$, they must be factored out to obtain the 
group of affinities of~$M$. 

The intersection of parts I and II is the translation subgroup of $\G$.
If we factor out both, we obtain the outer automorphism group of $\G$,
which is equally the discontinous part of the (affine) symmetry-group of $M$.
The elements of the connected and inner parts are isometric, but not all
elements of the outer automorphism group need be so.
We therefore separate it into:

{\it Part III}, the {\it rigidly isometric part}, the normal subgroup that contains precisely
those elements that are isometric for all the values of the parameters. Factoring these out, we
obtain:

{\it Part IV}, the {\it deformable part}. For example, the variant form (\ref{bodo2}) for $+a1$
shows that it has an automorphism interchanging $W$ and $X$ (see Fig.~\ref{xfigamphic}),
which is not always an isometry. It will however became an isometry if $B=C$. In general there
are certain finite subgroups, the {\it Bravais subgroups} $B_1,B_2,\dots$, of the deformable part
that just consist of the symmetries that
can become isometric for various choices of the parameters.
They are completely enumerated in Section~\ref{sbravo}.

Combining the results of that section and this, our paper describes all isometries
of all platycosms.

These groups are tabulated in Table~7.
We shall discuss only the two hardest cases, $c22$ and $-a1$.

\subsubsection*{The second amphicosm}

For $-a1$ it is easily checked that the most general inner automorphism $(m,n,\sigma)$, say, takes
$W\to WZ^m$, $X\to XZ^n$, $Z\to Z^{\sigma}$, where $\sigma=\pm 1$, $m-n=\frac{\sigma-1}2$, the ones
induced by $W,X,Z$ being $(0,+1,-)$, $(-1,0,-)$, $(2,2,+)$, respectively.

Now the glide planes of $-a1$ from a parallel series, say
$$\dots\overset{\pi_{-1}}{\big|}\quad\overset{\pi_0}{\big|}\quad\overset{\pi_1}{\big|}\quad\overset{\pi_2}{\big|}
\dots\big|\dots\big|\dots$$
and the inner automorphisms can move the $W$ and $X$ planes (say $\pi_0$ and $\pi_1$) any even number
of steps, while the further automorphism that takes $(W,X,Y)$ to $(WZ,XZ,Z^{-1})$ moves them just one step.

Modulo these we may suppose $\pi_0$ and $\pi_1$ are fixed or interchanged. The most general automorphism that
does this is
$$W\to W^a X^b,\quad X\to W^c X^d,\quad Z\to Z^\epsilon \quad (\epsilon=\pm 1),$$
where
\begin{eqnarray}
\text{if } \ep=+1,\,\, \left[\begin{array}{cc} a & b \\ c & d \end{array}\right]&\equiv&
\left[\begin{array}{cc}1 & 0 \\ 0 & 1 \end{array}\right]\,\, ({\rm mod}\,\, 2),  \quad  \text{while} \label{g2}
\\[.3cm]
\text{if } \ep=-1,\,\, \left[\begin{array}{cc} a & b \\ c & d \end{array}\right]&\equiv&
\left[\begin{array}{cc}0 & 1 \\ 1 & 0 \end{array}\right]\,\, ({\rm mod}\,\, 2). \label{g2b}
\end{eqnarray}
The group defined by (\ref{g2}) is often called $\G(2)$.
This explains the entry in Table~7.

{\small
{\setlength\arraycolsep{0pt}
\begin{table}[!htb]\label{tinner}
\begin{tabular}{|cc@{}c@{}c@{\hspace{-.4cm}}cc@{}c|}
\cline{2-7}
\multicolumn{1}{c|}{} & \multicolumn{1}{|c|}{I} & \multicolumn{1}{|c|}{II} & \multicolumn{2}{|c|}{III}
& \multicolumn{1}{|c|}{IV} & \multicolumn{1}{|c|}{V} \\
\cline{2-7}
\multicolumn{1}{c|}{} & \multicolumn{1}{|c|}{$\be_1$} &  \multicolumn{1}{|c|}{inner automorphisms}
& \multicolumn{3}{|c|}{outer automorphisms} & \multicolumn{1}{|c|}{$\#$ of} \\
\cline{4-6}
\multicolumn{1}{c|}{} & \multicolumn{1}{|c|}{} & \multicolumn{1}{|c|}{} & \multicolumn{2}{|l|}{rigidly isometric}
&  \multicolumn{1}{|r|}{deformable} & \multicolumn{1}{|c|}{$B_i$'s} \\
\hline
$c1$ & 3 & trivial & $X^{-1},Y^{-1},Z^{-1}$ & 2 & $\PGl_3(\Z)\,\mathrm{on}\,X,Y,Z$ & 14 \\
\hline
$c2$ & $1$  &
$\begin{array}{c} X^\epsilon,Y^\epsilon,X^a Y^b Z  \\ : a,b\in 2\Z \end{array}$  &
$\begin{array}{c} X,Y,XZ \\ X,Y,YZ \\ X,Y,Z^{-1} \end{array}$ & $2^3$
& $\PGl_2(\Z)\,\mathrm{on}\,X,Y$ & 5 \\
\hline
$c3$ & 1  &
$\begin{array}{c} X,Y,X^a Y^b Z  \\ Y,W,X^a Y^b Z  \\ W,X, X^a Y^b Z  \\
: a+b\in 3\Z \end{array}$
& $\begin{array}{c} X,Y,XZ \\ X^{-1},Y^{-1},Z \\ Y,X,Z^{-1} \end{array}$  & $2\times S_3$ & 1 & $1$ \\
\hline
$c4$ & 1  &
$\begin{array}{c}X,Y,X^a Y^b Z\\ Y,X^{-1},X^a Y^b Z\\ X^{-1},Y^{-1},X^a Y^b Z
\\ Y^{-1},X,X^a Y^b Z  \\ :a+b\in 2\Z \end{array}$
& $\begin{array}{c} X,Y,XZ  \\ Y,X,Z^{-1} \end{array}$  & $2^2$ & $1$ & $1$ \\
\hline
$c6$ & 1 &
$\begin{array}{c} X^\sigma,Y^\sigma,X^a Y^b Z \\ Y^\sigma,W^\sigma,X^a Y^b Z  \\ W^\sigma,X^\sigma, X^a Y^b Z
\\ :\sigma=\pm1,\, a+b\in\Z \end{array}$
& $Y,X,Z^{-1}$  & $2$ & $1$ & $1$ \\
\hline
$c22$ & 0 &
$\begin{array}{c} X^TY^b Z^c, Y^TX^a Z^c, Z^TX^a Y^b \\
: a,b,c\in 4\Z,\,\, T=1,X,Y,Z \end{array}$
& $\begin{array}{c} X,YX^2,ZX^2 \\  XY^2,Y,ZY^2 \\ XZ^2,YZ^2,Z \\ YZ^{-1},ZX^{-1},XY^{-1} \end{array}$ & $2^4$
& $S_3$  & $3$ \\
\hline

$+a1$ & 2  & $\begin{array}{c} WZ^{m}, XZ^{m}, Z^{\sigma} \\ :\sigma=\pm1,\, m\in2\Z \end{array}$ &  $WZ,XZ,Z^{-1}$ & $2$
& $\begin{array}{c} \G(2)\,\mathrm{on}\,W,X \\ X,W,Z^{-1} \end{array}$ & $5$ \\
\hline

$-a1$ & 2  &
$\begin{array}{c} WZ^m, XZ^n, Z^\sigma  \\ :\sigma=\pm1,\, m-n=\frac{\sigma-1}2\end{array}$ & $WZ,XZ,Z^{-1}$ & $2$
& $\begin{array}{c} \G(2)\,\mathrm{on}\,W,X \\ X,W,Z^{-1} \end{array}$ & $5$ \\
\hline

$+a2$ & 1 &
$\begin{array}{c} WX^{2a}Z^{2b}, X^\delta Z^{2b},Z^\epsilon \\
: \delta,\epsilon=\pm1, 2a\equiv\delta-\epsilon\mod 4 \end{array}$
& $\begin{array}{c} WZ,XZ,Z \\ W^{-1},X,Z \\ W,X^{-1},Z \end{array}$ & $2^3$ & $1$ & $1$ \\
\hline

$-a2$ & 1 &
$\begin{array}{c} WX^aZ^b, X^\delta Z^c,Z^\epsilon  \\
: \delta,\epsilon=\pm1, 2c+\delta\equiv1\mod 4 \\ a\equiv\delta-\epsilon\mod 4, 2b+1=2c+\epsilon\end{array}$
& $\begin{array}{c} WZ,XZ,Z \\ W^{-1},X,Z \\ W,X^{-1},Z \end{array}$ & $2^3$ & $1$ & $1$ \\
\hline
\end{tabular}
\bigskip
\caption{Columns I--IV concern the appropriate parts. Thus I gives the dimension of the connected component, II
the generic inner automorphisms, III generators and structure (modulo inner automorphisms) of the rigidly isometric
part, IV the deformable part. Column V gives the number of distinct types of Bravais subgroups $B_i$ (see
Section~\ref{sbravo}). Automorphisms are specified by giving images of $X,Y,Z$ for the chiral platycosms, $W,X,Z$
for the achiral ones.}
\end{table}
}}

\subsubsection*{The Didicosm}

To fix our ideas we take $A=B=C=1$.
Then the generators $X,Y$ are two screw motions of the minimal length~1 at the minimal distance~$1/2$.
The same must be true of their images, since if either length or distance were to be increased the resulting elements would
no longer generate. Now $X,Y$ determine $Z$ via $XYZ=1$. The lines of these three generators are edges of a unique
$\frac 12\times\frac 12\times\frac 12$ cube (or `{\it cubelet}') that constitutes one eighth of our defining `box'.
The general automorphism will take it to a similar cubelet, constituting an eighth of some other box of the tessellation.

We first show that there are enough automorphisms to take this standard cubelet to any of the eight
cubelets in any box. It may be taken to the corresponding one in any other box by a translation
through a typical vector $\al\x+\be\y+\g\z$ of the naming lattice $\m{N}=\langle \x,\y,\z\rangle$, which
achieves
$$(X,Y,Z)\longmapsto (XY^{2\be}Z^{2\g},YX^{2\al}Z^{2\g},ZX^{2\al}Y^{2\be}).$$
(This is inner just if $\al,\be,\g$ are even, since then the associated translation belongs
to $\m{T}=\langle 2\x,2\y,2\z\rangle$.)

The automorphisms we have just found effectively allow us to suppose that there is only one box, which contains
just eight cubelets. The inner automorphisms by $X,Y,Z$ take the standard cubelet to four of these.
(We have now found all inner automorphisms, since these elements represent the four cosets of~$\m{T}$ in~$\G$.)

So we need only find a further automorphism that takes the standard cubelet to one of the missing four. This is
$$(X,Y,Z)\longmapsto (YZ^{-1},ZX^{-1},XY^{-1}).$$
(the latter three elements are half turn screw motions about lines bisecting the other three faces of the
standard box.)

All the above automorphisms fix the three coordinate directions, so will be isometric even in
the general case ($A,B,C$ distinct). Modulo them, we can suppose $X,Y,Z$ go to some
permutation of $X^{\pm 1},Y^{\pm 1},Z^{\pm 1}$, since these are the only length~1 screw motions
associated with the initial cubelet. The condition $XYZ=1$ restricts us to the cyclic permutations
of $X,Y,Z$ or $Z^{-1},Y^{-1},X^{-1}$. The identity is the only one of these that is isometric
in the general case.

The conclusion is that the general inner automorphism takes $(X,Y,Z)$ to
$$(X^TY^{2\be}Z^{2\g},Y^TX^{2\al}Z^{2\g},Z^TX^{2\al}Y^{2\be}), \quad (\al,\be,\g \text{ even }, T=1,X,Y,Z);$$
while the full automorphism group is obtained by letting $\al,\be,\g$ be arbitrary integers, and adjoining
the maps taking $(X,Y,Z)$ to any even permutation of $(X,Y,Z)$
or $(YZ^{-1},ZX^{-1},XY^{-1})$ or any odd permutation of the inverses of those.
Note that the argument proves $c22^{A\,B\,C}$ to be invariant under all permutations of $A,B,C$
(not just the more obvious even ones).

Let $\theta_1,\theta_2,\theta_3,\phi$ be the automorphisms that take $(X,Y,Z)$ to
$$(X,YX^2,ZX^2),\,(XY^2,Y,ZY^2),\,(XZ^2,YZ^2,Z),\,(YZ^{-1},ZX^{-1},XY^{-1}).$$

Then it is easily checked that modulo inner automorphisms, $\theta_1,\theta_2,\theta_3,\phi$ are
multiplicatively commuting elements of order two.
The elements $\pi$ that take $(X,Y,Z)$ to even permutations of $(X,Y,Z)$ or odd permutations of $(X^{-1},Y^{-1},Z^{-1})$
form an $S_3$ that acts in the obvious way on  $\theta_1,\theta_2,\theta_3$ and takes
$\phi$ to  itself or $\theta_1\theta_2\theta_3\phi$, according as $\pi$ is even or odd.
This shows that the outer automorphism group is a split extension $2^4{:}S_3$.
The structure of this group was found (in a different form) by Zimmermann \cite{Zim}, who disproved
an earlier assertion (\cite{ChV}, \cite{Ch}) that the group was a split extension of form $2^3{:}(2\times S_3)$.

The outer automorphism groups of flat manifolds have often been studied (e.g., \cite{Ch}, \cite{MS}, \cite{Hi}),
the last of which used them to enumerate compact flat 4-manifolds. Our treatment avoids uses of subtle
theorems, and clearly distinguishes the isometric part.

\section{Double Covers}

A double cover of a manifold $M$ determines and is determined by a homomorphism from its fundamental
group $\pi_1(M)$ onto $\{\pm 1\}$, the image of a path being $-1$ just if it interchanges  the two sheets.
To enumerate such homomorphisms, we can obviously abelianize $\pi_1(M)$ (to $H_1(M)$) and then make
$H_1(M)$ have exponent 2, turning it into an elementary abelian 2-group.
If there are $r$ generators after this, the number of double covers will be $2^r-1$. We discuss the cases.

\subsection*{Torocosm}
Since $H_1(M)=\infty^3$, a torocosm has 7 double covers, themselves all torocosms.

\subsection*{Dicosm}
Now $H_1(M)=2^2\cdot\infty$, so the dicosm also has 7 double covers, of which just one is a torocosm. Each
of the other 6 is a dicosm, since some period 2 screw motion survives in the kernel of its defining homomorphism.

\subsection*{Tricosm}
Since $W,X,Y$ are conjugate in $\pi_1(M)$, they must map to the same sign, which must be $+$, since $WXY=1$.
There is therefore just one double cover (given by $W\mapsto +$, $X\mapsto +$, $Y\mapsto +$, $Z\mapsto -$)
which is another tricosm $c3$.

\subsection*{Hexacosm} A similar argument shows that the hexacosm $c6$ also has a single double cover, a tricosm $c3$.

\subsection*{Tetracosm}
$X$ and $Y$ are again conjugate, but now can both map to $+$ or both to $-$. There are therefore three double covers:
$(X,Y,Z)\mapsto ++-$ (type $c2$) or $--+$ or $---$ (both of type $c4$).

\subsection*{Didicosm}
The relation $XYZ=1$ shows that two of $X,Y,Z$ must map to $-$ and the third to $+$. So we have three double covers,
all dicosms $c2$.

\subsection*{Amphicosms and Amphidicosms}
The abelianized groups show that the positive cases each have seven double covers and the negative ones three.
The types are

for $+a1$: 1 of type $c1$, 4 of type $+a1$, 2 of type $-a1$;

for $-a1$: 1 of type $c1$, 2 of type $+a1$;

for $+a2$: 1 of type $c2$, 2 of type $+a1$, 2 of type $+a2$, 2 of type $-a2$;

for $-a2$: 1 of type $c2$, 2 of type $+a1$.

\begin{table}[!hbt]\label{tdoublecovers}
\begin{tabular}{|cc@{}c@{}c@{}cccc|}
\hline
platycosm & \multicolumn{4}{c}{generators}  & double covering & $\#$ & total $\#$ \\
  &  $W$ & $X$ & $Y$ & $Z$ & & & \\
\hline
$c1$ & & & & & $c1$ & 7 & 7 \\
\hline
$c2^D_{A\,B\,C}$ & $(+$ & $+$ & $+)$ & $-$ & $c1^{4D}_{A\,B\,C}$  & 1 & \\
                 & $(+$ & $-$ & $-)$ & $\pm$ & $c2^{D}_{[A]\,B\,C}$, etc. & 2 of each  & 7 \\
\hline
$c3^D_{A\,A\,A}$ & $(+$ & $+$ & $+)$ & $-$ & $c3^{4D}_{A\,A\,A}$ & 1 & 1 \\
\hline
$c4^D_{A\,A}$    &  & $(+$ & $+)$ & $-$ & $c2^{4D}_{A\,A}$ & 1 & \\
            & & $(-$ & $-)$ & $\pm$ & $c4^D_{2A\,2A}$ & 1 of each & 3 \\
\hline
$c6^D_{A\,A\,A}$ & $(+$ & $+$ & $+)$ & $-$ & $c3^{4D}_{A\,A\,A}$ & 1 & 1 \\
\hline
$c22^{A\,B\,C}$ &  & $(+$ & $-$ & $-)$ & $c2^A_{4B\,4C}$, etc. & 1 of each & 3 \\
\hline
\end{tabular}
\bigbreak
\begin{tabular}{|ccc|cccc|}
\hline
$X$ & $Y$ & $Z$ & $+a1^D_{A:B\,C}$ & $-a1^D_{A:B\,C}$ & $+a2^D_{A:B}$ & $-a2^D_{A:B}$ \\
\hline
$-$&$+$&$+$ & ${c1^D_{[A]\,B\,C}}^{\phantom N}$ & $c1$ see (\ref{fourforms})  & $+a1^{4A}_{B:D}$
& $+a1^{4A}_{[0]:B\,D}$ \\[.2cm]
$+$&$-$&$+$ & $+a1^D_{[B]:A\,C}$ & $+a1^{4D}_{[B]:A\,C}$ & $+a1^D_{A:4B}$ & $+a1^{4D}_{A:4B}$ \\[.2cm]
$-$&$-$&$+$ & $+a1^D_{[C]:A\,B}$ & $+a1^{4D}_{[C]:A\,B}$ & $c2^B_{4A\,D}$ & $c2^B_{4A\,4D}$ \\[.2cm]
$\pm$&$+$&$-$ & $+a1^{4D}_{A:B\,C}$ & ---  & $+a2^{4D}_{A:B}$ & --- \\[.2cm]
$\pm$&$-$&$-$ & $-a1^D_{A:B\,C}$ & --- & $-a2^D_{A:B}$ & --- \\[.2cm]
\multicolumn{3}{|c|}{total $\#$} & 7 & 3 & 7 & 3 \\[.1cm]
\hline
\end{tabular}
\bigskip
\caption{Parameters for the double covers.}
\end{table}
\subsection*{Parameters}
Tables~8 give the types and parameters of the double covers corresponding to
all homomorphisms : $\pi_1(M)\to\{\pm{1}\}$.
The parameters for the orientable double covers ($c1$) of $-a1$ take one of the four forms
\begin{equation}\label{fourforms}
\begin{array}{ccc}
c1^{2D\,\,2D}_{2C\,\,2C+4A\,\,B-C-D} & c1^{2D\,\,2D}_{2B\,\,2B+4A\,\,C-B-D} &
c1^{2B\,\,2B\,\,4A}_{2C\,\,2C\,\,D-B-C}  \\[.1cm]
(\text{if } B\ge C+D) & (\text{if } C\ge B+D) & (\text{if } D\ge B+C)  \\[.1cm]
 & c1^{C+D-B\,\,B+D-C\,\,B+C-D}_{C+D-B\,\,B+D-C\,\,B+C-D+4A} & \\[.1cm]
 & (\text{otherwise}) &
\end{array}
\end{equation}
and in several other cases, there are two forms, which we abbreviate using the notation
\begin{equation}\label{notation}
[u]{:}v\,w = \left\{ \begin{array}{ll} w-v:2v\,\,2v+4u &  \text{ if } v\le w; \\
v-w:2w\,\,2w+4u & \text{ if } v\ge w, \end{array}\right.
\end{equation}
and the similar notation without a colon.

\section{Diameters and Injectivity Radii}\label{sdiam}

The {\it covering radius} of a lattice $\m{T}$ is the minimal radius for which the closed balls
of that radius centered at points of $\m{T}$ cover the space.
It coincides with the {\it diameter} of the torus $R^n/{\m{T}}$,
where the diameter of a Riemannian manifold is defined in general to be the maximal distance between any pair
of points in it.

The {\it packing radius} of $\m{T}$ is the maximal radius
for which the open balls centered at points of $\m{T}$ are disjoint.
In the torus $R^n/{\m{T}}$, it becomes the {\it injectivity radius}, which is defined
for every Riemannian manifold to be the maximal radius for
which balls of that radius centered at any point of the manifold will embed in the manifold.
These are easier to determine:

\subsection*{Injectivity Radii}

\begin{proposition} The injectivity radius of a platycosm is bounded below by the packing radius of its naming lattice $\m{N}$.
Moreover, these two concepts coincide for 8 of the 10 platycosms.
\end{proposition}
\proof The injectivity radius equals one half of the length of the shortest closed geodesic in the manifold.
Every closed geodesic corresponds to a translation, screw motion or glide reflection
and its length is exactly the length of the translation vector, screw vector or glide vector, respectively,
and all these are by definition in the naming lattice $\m{N}$. This proves the first assertion.

On the other hand, one can check that the shortest possible vectors of
$\m{N}$ do correspond to such `geodesic' vectors in 8 out of the 10 cases, the exceptions being the two
`negatives' $-a1$ and $-a2$.
\endproof

We now discuss these two exceptions.

\subsubsection*{Injectivity radius of $-a1$}
Recall that $\m{N}=\langle \w,\x,\z\rangle$, whose conorms and vonorms
we display:
\begin{figure}[!htb]
\input{con-von+a1.pstex_t}
\caption{}
\label{figcon-von+a1}
\end{figure}

We now classify the geodesic vectors:

\

\begin{tabular}{ccc}
kind of vector &  vectors & norm of minimal vector \\
\hline
glide vector       & $odd\cdot \w+even\cdot \x$     & $A+B=N(\w)$ \\
                   & $even\cdot \w+odd\cdot \x$      & $A+C=N(\x)$ \\
\hline
translation vec. & $odd\cdot \w+odd\cdot \x+odd\cdot \z$  & $B+C+D=N(\w+\x+\z)$ \\
                    & $even\cdot \w+even\cdot \x+even\cdot \z$  & $4\cdot \min(\text{\it all vonorms})$
\end{tabular}

\

\noindent seeing that they lie in four cosets of $\m{N}$ in $2\m{N}$.
But by vonorm theory (see Appendix II), the minimal norms in these cosets are the numbers in the right column.

So the squared injectivity radius is the minimum of the ten numbers
$$\begin{array}{c} A+B,A+C,A+B+D,4D,4(B+C) \\
4(A+B),4(A+C),4(A+B+D),4(A+C+D),4(B+C+D), \end{array}
$$
of which those in the second line are dominated. Any of the five numbers in the
first line can be the minimum!

\subsubsection*{Injectivity radius of $-a2$}
The vectors corresponding to a letter `{\bf d}' in the figure are of the form $odd\cdot \x+even\cdot \w$.
In this way we obtain

\

\begin{tabular}{cccc}
letter & vector & minimal vector & norm \\
\hline
{\bf d} & $odd\cdot \x+even\cdot \y$ & $\x$ & $A$ \\
{\bf q} & $odd\cdot \y$ & $\y$ & $B$ \\
{\bf p} & $odd\cdot \y+odd\cdot \z$ & $\y+\z$ & $B+D$ \\
{\bf b} & $even\cdot \x+even\cdot \y+even\cdot \z$ & \text{(in $2\m{N}$)} & $4\cdot\min(\text{\it vonorms})$
\end{tabular}

\

Hence the answer is $\min(A,B,4D)$, since the other numbers are dominated by $A$ or $B$.

For $c1$ the answer is the minimal vonorm --- see the end of this section.

For $cN^D_{A\,B\,C}$ ($N>1$), it is generically $\min(B+C,C+A,A+B,D)$.
The particular cases are listed in Table~9.

\subsection*{Diameters}

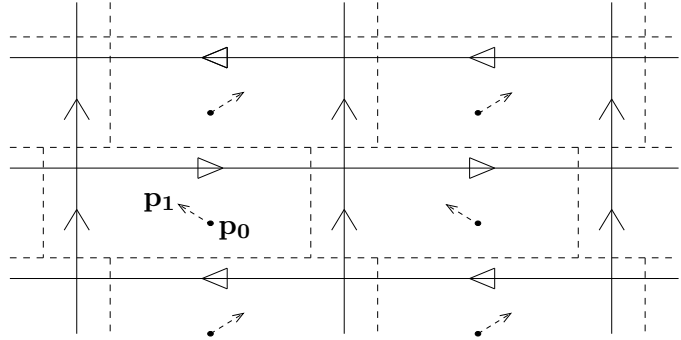
\begin{figure}[!htb]
\input{kb.pstex_t}
\caption{Translating the Voronoi cells.}
\label{figkb}
\end{figure}

We illustrate our method for finding diameters by considering that of the Klein bottle.
The solid rectangles in Figure \ref{figkb} are the Voronoi cells (see Section \ref{sbravo}) of the
images of $\mb{p_0}$ under $\G$.
This is an orbit which is also a lattice, so that the rectangles form
the {\it Voronoi tiling} of an {\it orbit lattice}.

\noindent Moreover, this tiling is {\it translatable} in the sense that
it remains a tiling if we move its tiles (say as indicated by the dashed lines)
so that their centers form the orbit of an arbitrary other point $\mb{p_1}$
(although the tiles might no longer be Voronoi cells).

It follows that the diameter of the Klein bottle is exactly the circumradius of these rectangles,
since the corners of the original rectangles are exactly this distance from $\mb{p_0}$, and clearly
no point can be further than this distance from the generic point $\mb{p_1}$.
The argument proves in general

\begin{proposition} The diameter of a platycosm $P$ is bounded below by the covering radius of
any orbit lattice.
If $P$ has an orbit lattice whose Voronoi tiling is translatable then
its diameter equals the covering radius of that orbit lattice. \qed
\end{proposition}

The {\it orbit lattice bound} $R_0$ for the diameter of a platycosm $P$ is the
maximum of the covering radii of all orbit lattices for $P$, and we say $P$ has the
{\it orbit lattice property} if its diameter equals this lower bound. We believe

\newtheorem*{conj}{(The orbit lattice conjecture for platycosms)}
\begin{conj} All 10 platycosms have the orbit lattice property.
\end{conj}

The proposition establishes this for 7 of the 10 since they have  orbit lattices whose Voronoi tilings are
translatable. The authors have worked out the orbit lattice bound for the three exceptions, and now
summarize the results.

\subsubsection*{Didicosm}

The didicosm $c22^{A\,B\,C}$ has 8 orbit lattices, of just 4 distinct shapes, for which the values
of $4R^2$ are:
$$\al=A+\frac{(B+C)^2}{\max(B,C)},\quad \be=B+\frac{(C+A)^2}{\max(C,A)},\quad \g=C+\frac{(A+B)^2}{\max(A,B)}$$
and
$$\delta=A+B+C+m^2\cdot\max\left(\frac 1A+\frac 1B+\frac 1C-\frac 2m,0\right),\,\,\,\text{where } m=\min(A,B,C).$$

Hence the orbit lattice bound is given by $4R_0^2=\max(\al,\be,\g)$,
since it can be shown that $\delta$ is dominated by any of $\al,\be,\g$.

\begin{table}[!hbt]\label{tInjDiam}
$$\begin{array}{rlc}
\text{platycosm} & \text{squared injectivity radius} & \text{squared diameter} \\
\hline
\phantom{{c^D}^S}
c1_{A\,B\,C}^{D\,E\,F} & \text{the minimal vonorm --- see below} & \text{see Section \ref{sformtor}} \\[.2cm]
c2_{A\,B\,C}^D  & \min(B+C,C+A,A+B,D) & \frac{(B+C)(C+A)(A+B)}{4(BC+CA+AB)} + \frac D4 \\[.2cm]
c3_{A\,A\,A}^D & \min(2A,D)    & \frac 23 A+\frac D4 \\[.2cm]
c4_{A\,A}^D & \min(A,D) &  \frac 12 A +\frac D4 \\[.2cm]
c6_{A\,A\,A}^D & \min(2A,D) & \frac 23 A+\frac D4 \\[.2cm]
c22^{A\,B\,C} & \min(A,B,C) & \ge \frac 14 \max(\al,\be,\g) \\[.2cm]
+a1^D_{A:B\,C} & \min(A+B,B+C,A+C,D) & \frac{(B+C)(C+A)(A+B)}{4(BC+CA+AB)} + \frac D4 \\[.2cm]
-a1^D_{A:B\,C} & \min(A+B,A+C,B+C+D,4D,4(B+C)) & \text{see text} \\[.2cm]
+a2^C_{A:B} & \min(A,B,C)  & \frac 14 (A+B+C) \\[.2cm]
-a2^C_{A:B} & \min(A,B,4C) & \ge \frac 14\max(\be,\g)
\end{array}$$
\bigskip
\caption{Injectivity radii and diameters. Here \newline
$\al=A+\frac{(B+C)^2}{\max(B,C)}$, $\be=B+\frac{(C+A)^2}{\max(C,A)}$, $\g=C+\frac{(A+B)^2}{\max(A,B)}$.}
\end{table}

\subsubsection*{Second amphidicosm}
The second amphidicosm $-a2^C_{A:B}$ has infinitely many orbit lattices, of just three different shapes,
for which the values of $4R^2$ are the above $\be,\g,\delta$. The
orbit lattice bound is therefore given by $4R_0^2=\max(\be,\g)$.

\subsubsection*{Second amphicosm}
The hardest case is $-a1^D_{A:B\,C}$, which has just 2 distinct shapes of orbit lattices.
For the first of these,

\begin{equation}\label{maxs}
4R^2=\left\{ \begin{array}{ll}
I & \textrm{if  $C\ge D$}, \\
\max(II,V) & \textrm{if  $C\le D\le A+C,B+C$}, \\
\max(III,V) & \textrm{if  $B+C\le D$, $B\le A$}, \\
\max(IV,V) & \textrm{if  $B+C\le D$, $B\ge A$},
\end{array}\right.
\end{equation}
where
\begin{eqnarray*}
I&=&A+B+4D+(A+B)(D-C)^2/\Omega \\
II&=&A+B+4C+(A+B+4C)(D-C)^2/\Omega \\
III&=&2B+2C+D+(B+C)^2/D+(B+C)(B-A)^2/\Omega \\
IV&=&2A+2C+D+(A+C)^2/D+(A+C)(B-A)^2/\Omega \\
V&=&A+B+2C+C^2/D
\end{eqnarray*}
where $\,\Omega=BC+CA+AB$.
Equivalently, it can be shown to be the minimum of the four expressions in (\ref{maxs}),
which may be rewritten as
$$\max\Big(\min(I,II,III,IV),\min(I,V)\Big).$$
For the second orbit lattice, $4R^2$ is given by the similar expression
$$\max\Big(\min(I',II',III',IV'),\min(I',V')\Big).$$
found by interchanging $B$ and $C$.

Therefore the orbit lattice bound is given by
$$4R_0^2=\max\Big(\min(I,II,III,IV),\min(I,V),\min(I',II',III',IV'),\min(I',V')\Big).$$

\section{Some Formulae for Torocosms and Lattices}\label{sformtor}

Many functions of a 3-dimensional lattice (or its associated torocosm) can
be expressed as symmetrical functions of all 7 conorms.
For instance, the squared packing or injectivity radius is the minimal vonorm,
and of course the typical conorm is just the sum of conorms not on a given line,
so the lattice $\Ld_{A\,B\,C}^{D\,E\,F\,G}$ has the conorms and vonorms of Figure \ref{figcon-von}.
\begin{figure}[!htb]
\hskip-1cm
\input{con-von.pstex_t}
\caption{}
\label{figcon-von}
\end{figure}
The $G$ conorm is not necessarily the zero one. If it is, the squared injectivity radius (of $c1_{A\,B\,C}^{D\,E\,F}$)
simplifies to
$$\min\big(B+E+C+F,C+F+A+D,A+D+B+E,D+B+C,A+E+C,A+B+F,D+E+F\big).$$

The determinant of this lattice (which is the squared volume of the corresponding torocosm) is the sum of the conorm
products over all 28 triangles in the conorm plane. When $G$ is zero, 12 of those products vanish,
and the determinant simplifies to
$$\Delta=(A+D)(BE+CF)+(B+E)(CF+AD)+(C+F)(AD+BE)+DBC+AEC+ABF+DEF.$$

The actual length of the edge in Figure \ref{figsensual} marked with a given conorm $x$ is
\begin{equation}\label{length}
\frac x{\Delta} \sqrt{\frac {\partial \Delta}{\partial x}}
\end{equation}
where $\frac {\partial \Delta}{\partial x}$
denotes the formal derivative of the above expression for $\Delta$ with respect to
the variable $x$. It is easy to check that (\ref{length}) can vanish only when $x$ does.

Another useful expression is $\Delta'$, the sum of conorm products over the complements
of all triangles; when $G=0$ this simplifies to
$$\Delta'=AD(B+E)(C+F)+BE(C+F)(A+D)+CF(A+D)(B+E),$$
since now 16 of the products vanish.
It can be shown that the covering radius is given (when $G=0$) by
\begin{equation}\label{diam}
4R^2=A+B+C+D+E+F- \frac{\Delta'+ 4\cdot\min (BECF,CFAD,ADBE)}{\Delta}
\end{equation}

\

Finally, we remark that the lattice $\Ld_{a\,b\,c}^{d\,e\,f\,g}$ dual to $\Ld_{A\,B\,C}^{D\,E\,F}$ has
conorms given by
$$\begin{array}{ll}
a\Delta=AD-m, & d\Delta=AE+AF+EF+m, \\
b\Delta=BE-m, & e\Delta=BD+BF+DF+m, \\
c\Delta=CF-m, & f\Delta=CD+CE+DE+m, \\
& g\Delta=AB+AC+BC+m, \\
\end{array}$$
where $m=\min(AD,BE,CF)$.

\section{The Bravais-Voronoi Classes}\label{sbravo}

\subsection*{Voronoi classification of lattices}

The points of space that are closer to one particular lattice point than to any other constitute
the {\it Voronoi cell} of that point. In \cite{sensualform} conorms were used to give a simple
derivation of Voronoi's classification of lattices by the topological type of their Voronoi cells.

\begin{figure}[!htb]
\hskip-1cm
\input{sensual.pstex_t}
\caption{The shape of the Voronoi cell.}
\label{figsensual}
\end{figure}
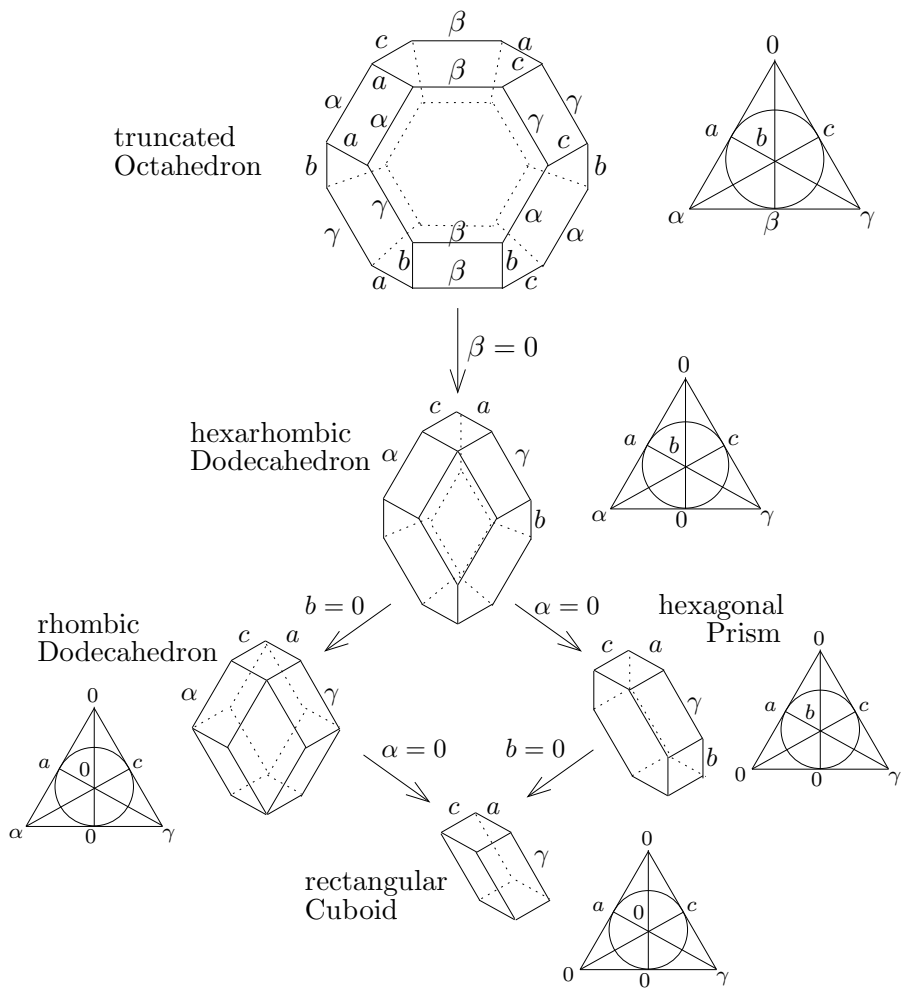

The generic Voronoi cell is a truncated octahedron as in Figure \ref{figsensual}, adapted from \cite{sensualform}.
The six families of parallel edges
correspond to its six non-zero conorms. When any conorm vanishes, the corresponding edges shrink
to points, changing the Voronoi cell into one of the four other shapes:
a hexarhombic dodecahedron, a rhombic dodecahedron, a hexagonal prism, or a rectangular cuboid (see Figure~\ref{figsensual}).

Each of these figures has antipodal symmetry. Any other topological symmetry will become a metrical one
provided certain conorms are equal.
Since the possible topological symmetries depend on which conorms vanish, we have found it convenient
to adjust the lettering to display this more clearly. The cases are:

\

{\setlength\arraycolsep{1pt}
\begin{tabular}{ccl}
lattice & (abbreviating) & \quad effect of symmetries on the conorms \\
\hline
$\Ld^{D\,E\,F}_{A\,B\,C}$ & &
\!\!\!\!\begin{tabular}{l}
Make two of the interchanges $\! ^{\phantom n}\!\! A\leftrightarrow D$, $B\leftrightarrow E$, $C\leftrightarrow F$.
\\ Bodily permute the three columns. \end{tabular} \\[.3cm]
$\Ld^{A\,B}_{C\,D:E}$ & $\Ld^{A\,B\,\text{O}}_{C\,D\,E}$
& \!\!\!\!\begin{tabular}{l}Make either of the interchanges $A\leftrightarrow B$, $C\leftrightarrow D$.
\\ Bodily permute the these two couples.\end{tabular} \\[.3cm]
$\Ld^{A\,B}_{C\,D}$ & $\Ld^{A\,B\,\text{O}}_{C\,D\,\text{O}}$
&  Permute $A,B,C,D$ in any way. \\[.3cm]
$\Ld^{D}_{A\,B\,C}$ & $\Ld^{D\,\text{O}\,\text{O}}_{A\,B\,C}$ &  Permute $A,B,C$ in any way. \\[.3cm]
$\Ld^{A\,B\,C}_{\phantom{\text{O}\,O\,O}}$ & $\Ld^{A\,B\,C}_{\text{O}\,\text{O}\,\text{O}}$ &  Permute $A,B,C$ in any way.
\end{tabular} }

\

\subsection*{Deducing the Bravai classification}
The arrows in the following figures show how certain equalities between conorms permit additional metrical
symmetries (the arrows being downward or rightward according as the symmetries move numbers within their columns
or across columns).
They show that 3-dimensional lattices belong to just 24 `BraVo' (Bravais-Voronoi) classes, where we say that
two lattices are in the same BraVo class just if each can be continuously deformed into the other while
keeping all its symmetries and without changing the topological shape of its Voronoi cell.

Below each conorm array is a naming letter A,B,\dots,X, and the orbifold notation for
the corresponding point group.

\underline{Truncated Octahedron} (tO).
\begin{displaymath}
\begin{array}{ccccc}
{\underset{\mathrm{A}}{\phantom {\begin{array}{c} z\\ z\end{array}} }
\framebox{$\begin{array}{ccc} a&b&c\\ \alpha&\beta&\g \end{array}$}}_{\times} & \longrightarrow
&
{\underset{\mathrm{D}}{\phantom {\begin{array}{c} z\\ z\end{array}} }
\framebox{$\begin{array}{ccc} a&a&b\\ \alpha&\alpha&\beta \end{array}$}}_{2*}
&  \longrightarrow &
{\underset{\mathrm{G}}{\phantom {\begin{array}{c} z\\ z\end{array}} }
\framebox{$\begin{array}{ccc} a&a&a\\ \alpha&\alpha&\alpha \end{array}$}}_{2{*}3}
\\
\big\downarrow & & \big\downarrow &
&
\raisebox{-.8cm}[0pt][0pt]{\begin{tabular}{c} $|$ \\[-.1cm] $|$ \\[-.1cm] $|$ \\[-.1cm] $|$ \\[-.1cm] $|$
\\[-.7cm] \\ \big\downarrow \end{tabular}}
\\
{\underset{\mathrm{B}}{\phantom {\begin{array}{c} z\\ z\end{array}} }
\framebox{$\begin{array}{ccc} a&b&c\\ a&b&\g \end{array}$}}_{2*} & \longrightarrow
&
{\underset{\mathrm{E}}{\phantom {\begin{array}{c} z\\ z\end{array}} }
\framebox{$\begin{array}{ccc} a&a&b\\ a&a&\beta \end{array}$}}_{*222} &
\hspace{-4cm} $--------------------\!\!\!\!$ \hspace{-5.6cm}
&
\\
\big\downarrow & & \big\downarrow & &
\\
{\underset{\mathrm{C}}{\phantom {\begin{array}{c} z\\ z\end{array}} }
\framebox{$\begin{array}{ccc} a&b&c\\ a&b&c \end{array}$}}_{*222} & \longrightarrow &
{\underset{\mathrm{F}}{\phantom {\begin{array}{c} z\\ z\end{array}} }
\framebox{$\begin{array}{ccc} a&a&b\\ a&a&b \end{array}$}}_{*422} &  \longrightarrow &
{\underset{\mathrm{H}}{\phantom {\begin{array}{c} z\\ z\end{array}} }
\framebox{$\begin{array}{ccc} a&a&a\\ a&a&a \end{array}$}}_{*432}
\end{array}
\end{displaymath}

\underline{Hexarhombic Dodecahedron} (hD).
\begin{displaymath}
\begin{array}{ccccc}
{\underset{\mathrm{I}}{\phantom {\begin{array}{c} z\\ z\end{array}} }
\framebox{$\begin{array}{ccc}a&b&\\ \alpha&\beta&:\g \end{array}$}}_{\times} & \longrightarrow &
{\underset{\mathrm{K}}{\phantom {\begin{array}{c} z\\ z\end{array}} }
\framebox{$\begin{array}{ccc} a&a&\\ \alpha&\beta&:\g \end{array}$}}_{2{*}} & \longrightarrow &
{\underset{\mathrm{L}}{\phantom {\begin{array}{c} z\\ z\end{array}} }
\framebox{$\begin{array}{ccc} a&a&\\ \alpha&\alpha&:\g \end{array}$}}_{*222} \\
$\big\downarrow$ & &
\raisebox{-.35cm}[0pt][0pt]{\begin{tabular}{c} $|$  \\[-.2cm] $|$  \\[-.2cm] $|$ \end{tabular}}
& & $\big\downarrow$  \\
{\underset{\mathrm{J}}{\phantom {\begin{array}{c} z\\ z\end{array}} }
\framebox{$\begin{array}{ccc} a&b&\\ a&b&:\g \end{array}$}}_{2{*}} &
$------$ &
\hspace{-1cm}$---------------------------------$\hspace{-.7cm}
& \raisebox{.020cm}[0pt][0pt]{$\longrightarrow $}
&
{\underset{\mathrm{M}}{\phantom {\begin{array}{c} z\\ z\end{array}} }
\framebox{$\begin{array}{ccc} a&a&\\ a&a&:\g \end{array}$}}_{*422}
\end{array}
\end{displaymath}

\underline{Rhombic Dodecahedron} (rD).
\begin{displaymath}
\begin{array}{ccccc}
{\underset{\mathrm{N}}{\phantom {\begin{array}{c} z\\ z\end{array}} }
\framebox{$\begin{array}{cc}a&b\\ \alpha&\beta \end{array}$}}_{\times} & \longrightarrow &
{\underset{\mathrm{O}}{\phantom {\begin{array}{c} z\\ z\end{array}} }
\framebox{$\begin{array}{cc} a&a\\ \alpha&\beta \end{array}$}}_{2*} & \longrightarrow &
{\underset{\mathrm{Q}}{\phantom {\begin{array}{c} z\\ z\end{array}} }
\framebox{$\begin{array}{cc} a&a\\ \alpha&\alpha \end{array}$}}_{*222} \\
& & $\big\downarrow$ & & $\big\downarrow$  \\
& &
{\underset{\mathrm{P}}{\phantom {\begin{array}{c} z\\ z\end{array}} }
\framebox{$\begin{array}{cc} a&a\\ a&\alpha \end{array}$}}_{2{*}3} & \longrightarrow &
{\underset{\mathrm{R}}{\phantom {\begin{array}{c} z\\ z\end{array}} }
\framebox{$\begin{array}{cc} a &a\\ a&a \end{array}$}}_{*432}
\end{array}
\end{displaymath}

\underline{Hexagonal Prism} (hP).
\begin{displaymath}
\begin{array}{ccccc}
{\underset{\mathrm{S}}{\phantom {\begin{array}{c} z\\ z\end{array}} }
\framebox{$\begin{array}{ccc}d\\ \alpha&\beta&\g \end{array}$}}_{2*} & \longrightarrow &
{\underset{\mathrm{T}}{\phantom {\begin{array}{c} z\\ z\end{array}} }
\framebox{$\begin{array}{ccc} d\\ \alpha&\alpha&\beta \end{array}$}}_{*222} & \longrightarrow &
{\underset{\mathrm{U}}{\phantom {\begin{array}{c} z\\ z\end{array}} }
\framebox{$\begin{array}{ccc} d\\ \alpha&\alpha&\alpha \end{array}$}}_{*622}
\end{array}
\end{displaymath}

\underline{Rectangular Cuboid} (rC).
\begin{displaymath}
\begin{array}{ccccc}
{\underset{\mathrm{V}}{\phantom {\begin{array}{c} z\\ z\end{array}} }
\framebox{$\begin{array}{ccc} a&b&c \\ \phantom{a} \end{array}$}}_{*222} & \longrightarrow &
{\underset{\mathrm{W}}{\phantom {\begin{array}{c} z\\ z\end{array}} }
\framebox{$\begin{array}{ccc} a&a&b\\ \phantom{a}\end{array}$}}_{*422}
& \longrightarrow &
{\underset{\mathrm{X}}{\phantom {\begin{array}{c} z\\ z\end{array}} }
\framebox{$\begin{array}{ccc} a&a&a\\ \phantom{a}\end{array}$}}_{*432}
\end{array}
\end{displaymath}

The easiest proof that there are just 14 Bravais classes is to see how certain continuous variations link
the 24 BraVo classes. For example, take the lattice
$\Ld ^{a\,a\,a}_{\al\,\al\,\al}$ (case G)
and vary $\al$ until it becomes a small negative number $-\epsilon$, forcing us to renormalize the conorms
as in Figure~\ref{figconbravo}.

\begin{figure}[!ht]
\input{conbravo.pstex_t}
\caption{}
\label{figconbravo}
\end{figure}

Then the new lattice has shape
$\Ld^{a\,a}_{a\,b}$ (case P),
showing that  cases G and P are in the same Bravais class.
To obtain the exact partition into Bravais classes one merely has to vary the parameters like this in all
possible ways for which the new lattice has no more symmetries than the
old one\ft{although at isolated intermediate times the symmetry might have been larger, as it was
in this case when $\alpha=0$}.
Similarly, making $b$ negative in case F leads to case M. All other identifications can be achieved just by 
making some conorms vanish, as in Figure \ref{figw}.
\begin{figure}[!ht]
\input{figw.pstex_t}
\caption{}
\label{figw}
\end{figure}

This gives Table~10, whose last column gives the ``symmetry factors'' by which the sizes of the isometry
groups have been increased.

\begin{table}[!htb]\label{tbravo}
\begin{tabular}{ccccclc}
\multicolumn{5}{c}{\it Voronoi classes} & {\it  Bravais classes} & \\
tO & hD & rD & hP & rC & crystallographic name & symmetry factor \\
\hline
A & I & N & & & Triclinic & 1 \\
B\,D & J\,K & O & & & base-centered Monoclinic  & 2 \\
C & L & Q & & & body-centered Orthorhombic  & 4 \\
E &  &  & & & face-centered Orthorhombic & 4  \\
F & M &  & & & body-centered Tetragonal & 8 \\
G &  & P & & & Rhombohedral (or Trigonal) & 6 \\
H & & & & & body-centered Cubic (bcc) & 24 \\
 & & R & & &  face-centered Cubic (fcc) & 24 \\
 & & & S & &  primitive Monocline & 2 \\
 & & & T & &  base-centered Orthorhombic & 4 \\
 & & & U & &  Hexagonal & 12 \\
 & & & & V &  Orthorhombic & 4 \\
 & & & & W & primitive Tetragonal & 8 \\
 & & & & X & primitive Cubic & 24 \\
\end{tabular}
\bigskip
\caption{How the 24 BraVo classes correspond to the 5 Voronoi and 14 Bravais classes.}
\end{table}

To prove this list is complete, it suffices to show that any other cases  with the same group are in distinct
Bravais classes, which follows from the fact that each of E,H,R,S,T,U,V,W,X is in a Bravais class of its own.
Why is this? For H,R,V,W,X no parameter can pass through $0$ without the lattice degenerating. The same
is true for the parameter $a$ in E and $d$ in S,T,U.
For E this leaves essentially only the normalization of
$\left(\begin{array}{ccc} a & a & \epsilon \\ a & a & \beta \end{array}\right)$
to
$\left(\begin{array}{ccc} a-\epsilon & a-\epsilon & \epsilon \\ a-\epsilon & a-\epsilon & \beta+\epsilon \end{array}\right)$,
which is again in case E.
The argument for S,T,U is even easier, and completes the proof of the Bravais classification.

\subsection*{Bravais types for the other platycosms}
For the remaining helicosms $cN_{A\,B\,C}^D$ the Bravais classification reduces to that of the 2-dimensional
lattice $\Ld_{A\,B\,C}$, since $D$ can be varied independently of $A,B,C$. Thus $c2$ has five Bravais types:
$$c2_{A\,B\,C}^D \longrightarrow c2_{A\,A\,B}^D \longrightarrow c2_{A\,A\,A}^D,\quad
c2_{A\,B}^D\longrightarrow c2_{A\,A}^D$$
where $c2^D_{A\,B}$ means $c2^D_{A\,B\,0}$, and $c3,c4,c6$ just one each
$$c3_{A\,A\,A}^D,\,\,c4_{A\,A}^D,\,\,c6_{A\,A\,A}^D.$$
For the didicosm $c22$, the classification is just that of the rectangular 3-dimensional lattice $\Ld^{A\,B\,C}$:
$$c22^{A\,B\,C} \longrightarrow c22^{A\,A\,B} \longrightarrow c22^{A\,A\,A}.$$
For the two amphicosms $\pm a1_{A:B\,C}^D$, the Bravais classification reduces to that of
the 2-dimensional lattice $\Ld_{A\,B\,C}$,
but taking into account the distinguished r\^ole of $A$:
$$+a1_{A:B\,C}^D \longrightarrow +a1_{A:B\,B}^D,\quad +a1_{:B\,C}^D\longrightarrow +a1_{:B\,B}^D,\quad +a1_{A:B}^D,$$
$$-a1_{A:B\,C}^D \longrightarrow -a1_{A:B\,B}^D,\quad -a1^D_{:B\,C}\longrightarrow -a1^D_{:B\,B},\quad -a1_{A:B}^D,$$
where, of course,  $\pm a1_{A:B}^D$ and $\pm a1_{:B\,C}^D$ abbreviate $\pm a1_{A:B\,0}^D$ and $\pm a1_{0:B\,C}^D$
respectively.

Finally, each of the amphidicosms $\pm a2^D_{A:B}$ forms a single Bravais type since all three parameters
are distinguishable. We conclude:

\smallbreak

{\it The numbers of Bravais types of platycosms are}

$\begin{array}{cccccccccccc}
 &14 & 5 & 1 & 1& 1 & 3 & 5 & 5 & 1 & 1 & (\mathit{total}\,\, 37) \\
\mathit{for} & c1 & c2 & c3& c4 & c6 & c22 & +a1 & -a1 & +a2 & -a2 &
\end{array}
$

\section*{Appendix I: Why there are just 10 platycosms}

Why are there exactly 10 compact platycosms?
We sketch a proof of this to make our paper complete.
The proof is outlined in Figure \ref{figtelegraphic} and the accompanying explanations (i)--(vi).
We quote Bieberbach's theorem that the subgroup of translations has finite index.

(i) When there are screw motions there is a
\newtheorem*{splitlemma}{Splitting Lemma}
\begin{splitlemma}
If $\vv$ is the smallest screw vector in a given direction, say vertical, then the
translation lattice $\m{T}$ decomposes as $\langle N\vv\rangle\oplus \m{T}'$, where $N$ is
the period of the corresponding screw motion~$\sigma$.
\end{splitlemma}
\begin{proof} It suffices to prove that there is no translation $\tau$ whose vertical part $c\vv$ is
strictly between $0$ and~$N\vv$.
But then $\sigma\tau^{-1}$ if $0<c<1$, or $\tau\sigma^{-m}$ if $m<c<m+1$ ($m=1,\dots,N-1$)
would be a vertical screw motion with shorter vector than~$\vv$.
\end{proof}
(ii) If all screw vectors are parallel we have a cyclic point group whose order $N$ is the least common
multiple of their periods and which acts on the 2-dimensional lattice $\m{T}'$.
Obviously $\G$ is generated by the shortest screw motion of this period together with the translations
of $\m{T}'$. The identification with $c2,c3,c4$ or $c6$ follows from the well known lemma:
\newtheorem*{Barlow}{Barlow's Lemma}
\begin{Barlow}
If a rotation of order $N$ fixes a 2-dimensional lattice $\m{T}'$, then $N=1,2,3,4,6$.
\end{Barlow}
\begin{proof}
Apply the rotation to a minimal non-zero vector of the lattice. We
obtain a `star' of $N$ vectors,
\begin{figure}[!ht]
\centerline{\mbox{\includegraphics*[scale=1]{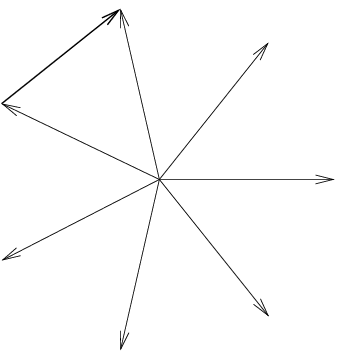}}}
\label{figbarlow}
\end{figure}
the difference of adjacent members of which will be a shorter vector if $N\ge 7$.
A rotation of order 5 would combine with negation to produce one of order 10.
\end{proof}

(iii) Two non parallel screw motions $\sigma_1$ and $\sigma_2$ would make $\m{T}$ decompose
in two ways, say $\langle \mb{v_1}\rangle\oplus\m{T}_1\cong\langle \mb{v_2}\rangle\oplus\m{T}_2$, showing
immediately that $\mb{v_1}$ and $\mb{v_2}$ are orthogonal. If $\sigma_i$ had a period other
than 2, then $\sigma_j$ and $\sigma_j^{\sigma_i}$ would not be orthogonal.
Finally, if $\sigma_1$ and $\sigma_2$  are orthogonal and of period 2, then
$\sigma_1\sigma_2$ is a screw motion in the direction of $\mb{v_3}$ orthogonal to $\mb{v_1}$ and $\mb{v_2}$.
These together with the translations must generate $\G$, and can be identified with the $X,Y,Z$ of~$c22$.

\begin{figure}[!ht]
\input{telegraphic.pstex_t}
\caption{Guide to the proof.}
\label{figtelegraphic}
\end{figure}
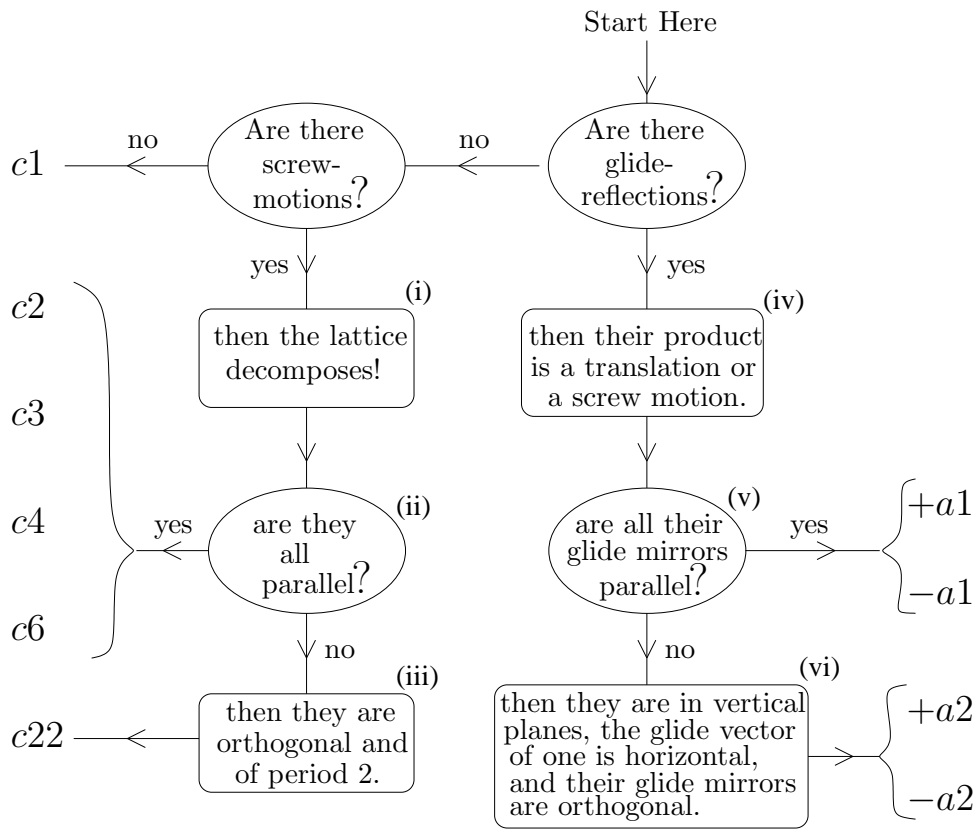

(iv) Any glide reflection maps to a reflection through the origin when ignore translations.
The product of two of them is therefore a translation or a screw motion since the
product of two reflections through {\bf o} is the identity or a non-trivial rotation according as
their planes are parallel or not.

(v) If all glide reflections are in parallel planes --- call them basal --- the point group has
order 2, so $\G$ is generated by a single glide reflection, together with its translations.
We resolve any translation vector into basal and perpendal parts $\mb{v_1}+\mb{v_2}$, and see that since
$\tau$ takes this to $\mb{v_1}-\mb{v_2}$, $2\mb{v_2}$ is also a translation vector. If $\mb{v_2}$ itself (and so $\mb{v_1}$)
is a translation vector, then the translation lattice $\m{T}$ decomposes into its basal
and perpendal parts, leading to $+a1$. Otherwise, adjoining the translation through $\mb{v_2}$
embeds our platycosm in a copy of $+a1$, and easily identifies the original platycosm with $-a1$.

(vi) If two glide reflections  $g,g'$  have non-parallel planes intersecting in a line we will
call vertical, then their product will be a screw motion $s$ whose screw vector $\mb{v}$ will be
vertical. We choose things so that $\mb{v}$ is as short as possible. Then we can multiply by
a power of $s$ to reduce $g$ until the vertical component
of its glide vector $\ld\mb{v}$ has $|\ld|\le \frac 12$. But then $\ld$ must be $0$
since $g^2$ is a translation and by the Splitting Lemma its vertical component must be a multiple
of $N\mb{v}$, where $N\ge 2$ is the period of $s$.

The glide vector of the new $g'$ (defined by $gg'=s$) will have the same vertical component $\mb{v}$
as $s$. But now the translation vector ${(g')}^2$ has vertical component $2\mb{v}$, showing (by the Splitting Lemma)
that $s$ must have period 2, from which it follows that the planes of $g$ and $g'$ are orthogonal.

We can now see that the translation lattice is generated by 3 orthogonal translations, say
$(x,y,z)\mapsto (x+a,y,z),\,(x,y+b,z),\,(x,y,z+c)$.
We can also see that the point group has order 4, i.e. that $g,g'$, together with translations, must
generate $\G$. For if not, there would be a glide mirror not parallel to either of those of $g,g'$,
and so perpendicular to both, by the previous paragraph. But then the product of the corresponding glide reflection
with $s$ would have the form
$$(x,y,z)\mapsto (\mathrm{constant}-x,\mathrm{constant}-y,\mathrm{constant}-z),$$
which fixes a point.
In these coordinates
\begin{eqnarray*}
s&:&(x,y,z)\to(-x,-y,z+\frac 12 c) \\
g&:&(x,y,z)\to(x+\frac 12 a,\ld b-y,z)
\end{eqnarray*}
where we may take $0\le\ld<1$ by compounding $g$ with translations.
From these we find
\begin{eqnarray*}
g'=g^{-1}s&:&(x,y,z)\to(\frac 12 a-x,y-\ld b,z+\frac 12 c) \\
{(g')}^2&:&(x,y,z)\to(x,y-2\ld b,z+c)
\end{eqnarray*}
showing that $2\ld$ must be an integer, so either $\ld=0$, which gives $+a2$; or $\ld=\frac 12$, which gives $-a2$.
In the former case, our screw vector $(0,0,\frac 12 c)$ is also a glide vector (of $g'$)  --- in the latter
case, no screw vector is a glide vector.

\section*{Appendix II: Conorms of Lattices}\label{sconorms}

\subsection*{Introduction}
The conorms of a lattice ${\m L}$ are certain numbers that are determined by ${\m L}$ (up to equivalence)
and return the compliment by determining ${\m L}$, at least in low dimensions.
More precisely ${\m L}$ has a {\it conorm function} defined on {\it conorm space}, which is a
finite set that has the structure of a projective $(n-1)$-dimensional space over the field of order two.
At least for $n\le 4$, two $n$-dimensional lattices ${\m L}$ and ${\m L}'$ are isometric if and only if there is an isomorphism between their conorm spaces that takes one conorm function to the other.

The theory in low dimensions is greatly simplified by the observation that for $n\le 3$ every lattice has
an {\it obtuse superbase}, which implies in particular that the conorms are $\ge 0$.
A {\it superbase} for an $n$-dimensional lattice ${\m L}$ is an $(n+1)$-tuple $\{\vv_0,\vv_1,\dots,\vv_n\}$ of vectors that
generate ${\m L}$ and sum to zero. It is {\it obtuse} if all inner products $\vv_i\cdot \vv_j$ of distinct vectors are non-positive
(and {\it strictly obtuse} if they are strictly negative).  For a lattice with an obtuse superbase, it can be shown that the
conorms are the negatives of the inner products of pairs of distinct superbase vectors, supplemented by zeros.
(If $n\ge 4$ other things can happen; a lattice may not have an obtuse superbase, and some conorms may be strictly negative.)

{\it For $n=0$, conorm space is empty}, so there are no conorms.

{\it For $n=1$, conorm space is a single point}, and there is one strictly positive conorm $A$.
We represent this by the picture $\,\stackrel{A}{\bullet}\,$; it means that the lattice has an obtuse suberbase
$\{\vv,-\vv\}$ with Gram-matrix $\left[\begin{array}{cc} A & -A \\ -A & A \end{array}\right]$;
equivalently a base $\{\vv\}$ with Gram-matrix $\left[ A \right]$.

{\it For $n=2$, conorm space is a 3-point projective line}, and the 3 conorms $A,\,B,\,C$
are non-negative.
We represent this by the picture
$\,\,\stackrel{A}{\bullet}\!\!\!\!-\!\!\!-\!\!\!-\!\!\!\!\stackrel{B}{\bullet}\!\!\!\!-\!\!\!-\!\!\!-\!\!\!\!
\stackrel{C}{\bullet}$~;
it means that the lattice has an obtuse superbase $\{\vv_0,\vv_1,\vv_2\}$ with matrix
$\left[\begin{array}{ccc} B+C & -C & -B \\ -C & C+A & -A \\ -B & -A & A+B \end{array}\right]$,
or equivalently a base $\{\vv_0,\vv_1\}$ with matrix
$\left[\begin{array}{cc} B+C & -C  \\ -C & C+A \end{array}\right]$.

\begin{figure}[!ht]
\centerline{\input{fig1.pstex_t}}
\caption{A Fano plane with $p_{01}, p_{02}, p_{03}, p_{12}, p_{13}, p_{23}, 0$ in the usual arrangement.}
\label{fig1}
\end{figure}
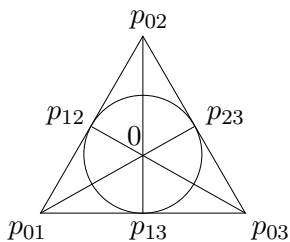

{\it For $n=3$, conorm space is the 7-point ``Fano plane''},
represented by Figure~\ref{fig1}.
The superbase has Gram-matrix
$$\left[\begin{array}{cccc} p_{0|123} & -p_{12} & -p_{13} & -p_{01} \\
-p_{12} & p_{1|023} & -p_{23} & -p_{02} \\
-p_{13} & -p_{23} & p_{2|013} & -p_{03} \\
-p_{01} & -p_{02} & -p_{03} & p_{3|012}
\end{array}\right],$$
where $p_{ij}=p_{ji}\ge 0$, $p_{i|jkl}:=p_{ij}+p_{ik}+p_{il}$, and
the minimal conorm is~0.

The sum of the four conorms not on a given line in Fig.~\ref{fig1} is a {\it vonorm}, which can
be defined either as the norm of a Voronoi vector, or as the minimal norm of any vector in a
non-trivial coset of $\m{L}$ in $2\m{L}$.

\subsection*{Putative conorms and the reduction algorithm}

A 3-dimensional lattice has many systems of ``putative conorms'' in addition to its unique system of  actual
conorms, for instance the numbers obtained by arranging $0$ and the negatives $p_{ij}:=-v_i\cdot v_j$ of the
inner products of distinct members of {\sl any} superbase on a Fano plane in the manner of Figure \ref{fig1}.
If the putative conorms are all non-negative, they will be the actual conorms.
Otherwise, the following algorithm quoted from \cite{sensualform} and \cite{CSVI} will produce the latter.

Select a `working line' that contains both a $0$ conorm and a negative one, say $-\epsilon$.
Then we transform to an improved system of putative conorms by adding $\epsilon$ to the 3 conorms
on the working line, and subtracting $\epsilon$ from the 4 conorms off this line.
If the improved system still has a negative conorm, we can define a new working line and repeat
the procedure. A finite number of repetitions will suffice to produce the actual conorms.

As an example Figure \ref{fig0} finds the conorms for the lattice whose Gram-matrix
with respect to a suitable base $\vv_1,\vv_2,\vv_3$ is
$\left[\begin{array}{ccc} 2&1&1 \\ 1&3&1 \\ 1&1&4 \end{array}\right]$
The Gram-matrix for the superbase $\vv_1,\vv_2,\vv_3,-\vv_1-\vv_2-\vv_3$ is
$\left[\begin{array}{cccc} 2&1&1&-4 \\ 1&3&1&-5 \\ 1&1&4&-6 \\ -4&-5&-6&15 \end{array}\right]$
(found by making each row and column sum to $0$),
which leads to the putative conorms of Figure \ref{fig0}(a).
\begin{figure}[!ht]
\input{fig0.pstex_t}
\caption{}
\label{fig0}
\end{figure}
Transforming this using the working line indicated on it, we obtain Figure \ref{fig0}(b), and we proceed from
this in a similar way to Figures \ref{fig0}(c), \ref{fig0}(d), \ref{fig0}(e),
the last of which gives the actual conorms.
The reader might like to verify that another choice of working line yields the same conorm function.

For a 2-dimensional lattice with putative conorms $A,\,B,\,C$ there is a similar algorithm. We select a
`working point' at which there is a negative conorm, say $-\epsilon$, and transform by adding $2\epsilon$
to this conorm and subtracting $2\epsilon$ from the other two conorms.
Thus
$\,\,\stackrel{-3}{\bullet}\!\!\!\!\!-\!\!\!-\!\!\!-\!\!\!\stackrel{5}{\bullet}\!\!\!-\!\!\!-\!\!\!-\!\!\!
\stackrel{10}{\bullet}\,\,$ transforms through
$\,\,\stackrel{3}{\bullet}\!\!\!-\!\!\!-\!\!\!-\!\!\!\!\!\stackrel{-1}{\bullet}\!\!\!\!\!-\!\!\!-\!\!\!-\!\!\!
\stackrel{4}{\bullet}\,\,$ to
$\,\,\stackrel{1}{\bullet}\!\!\!-\!\!\!-\!\!\!-\!\!\!\stackrel{1}{\bullet}\!\!\!-\!\!\!-\!\!\!-\!\!\!
\stackrel{2}{\bullet}\,$.

\section*{Appendix III: Dictionary of Names and Notations}\label{sdictionary}

The purpose of this appendix is to give to readers who have met with a platycosm
in some other notation a way to recognize it quickly.
Some of the notations are for the corresponding space groups.
We take the translations of these in the form
$$(x,y,z)\to (x+l,y+m,z+n)\quad (l,m,n\in\Z)$$
and use the additional generators given in the last column (in which the coordinates are not necessarily orthonormal).
These are taken from CARAT (cf. \cite{CS}), which has such information up to dimension~6.

\begin{table}[!htb]\label{tnamesplaty}
\begin{tabular}{|c|c|c|c|}
\hline
our name               & symbol & other names               & Wolf \\
\hline
torocosm               & $c1$   & 3-torus                   &  $\mathcal{G}_1$  \\
dicosm                 & $c2$   & half turn space        & $\mathcal{G}_2$   \\
tricosm                & $c3$   & one-third turn space   & $\mathcal{G}_3$   \\
tetracosm              & $c4$   & quarter turn space     & $\mathcal{G}_4$   \\
hexacosm               & $c6$   & one-sixth turn space   & $\mathcal{G}_5$   \\
didicosm               & $c22$  & Hantzsche-Wendt space  & $\mathcal{G}_6$   \\
first amphicosm        & $+a1$  & Klein bottle times circle & $\mathcal{B}_1$   \\
second amphicosm       & $-a1$  &                        & $\mathcal{B}_2$   \\
first amphidicosm      & $+a2$  &                        & $\mathcal{B}_3$   \\
second amphidicosm     & $-a2$  &                        & $\mathcal{B}_4$   \\
\hline
\end{tabular}
\bigskip
\caption{Names and notations for platycosms.}
\end{table}

\begin{table}[!htb]\label{tnotationsgroups}
\begin{tabular}{|c|c|r@{.  }l|c|}
\hline
symbol & \cite{CDHT}     &  \multicolumn{2}{c|}{\begin{tabular}{l} internatl.\\ no.\, name\end{tabular} }
& \begin{tabular}{c} non-translation \\ generators\end{tabular} \\
\hline
$c1$ & $(\circ)$         &  1 & $P1$  & --- \\
\hline
$c2$ & $(2_12_12_12_1)=(\bar\times\bar\times)$    & 4 & $P2_1$ & $(-x,-y,z+1/2)$ \\
\hline
$c3$ & $(3_13_13_1)$  &
\multicolumn{2}{@{\hspace{-.7cm}}c@{\hspace{-.0cm}}|}{ \begin{tabular}{r@{. }l} 144 & $P3_1$ \\ 145 & $P3_2$ \end{tabular} }
& $(-x-y,x,z+1/3)$ \\
\hline
$c4$ & $(4_14_12_1)$     &
\multicolumn{2}{@{\hspace{-.4cm}}c|}{ \begin{tabular}{r@{. }l} 76 & $P4_1$ \\ 78 & $P4_3$ \end{tabular} }
& $(-y,x,z+1/4)$ \\
\hline
$c6$ & $(6_13_12_1)$     &
\multicolumn{2}{@{\hspace{-.7cm}}c@{\hspace{-.0cm}}|}{ \begin{tabular}{r@{. }l} 169 & $P6_1$ \\ 170 & $P6_5$ \end{tabular} }
& $(x+y,-x,z+1/6)$ \\
\hline
$c22$ & $(2_12_1\bar\times)$  & 19 & $P2_12_12_1$ & $\begin{array}{c} (-x,y+1/2,-z+1/2) \\ (x+1/2,-y,-z) \end{array}$
\\
\hline
$+a1$ & $(\bar\circ_0)=(*{:}{*}{:})=(\times\times_0)$    & 7 & $Pc$ & $(x+1/2,y,-z)$ \\
\hline
$-a1$ & $(\bar\circ_1)=(*{:}\times)=(\times\times_1)$ & 9 & $Cc$ & $(x+1/2,z,y)$  \\
\hline
$+a2$ & $(2_12_1{*}{:})=(\bar*{:}\bar*{:})=(\bar\times\times_0)$ & 29 & $Pca2_1$ &
$\begin{array}{c} (-x,y,z+1/2) \\ (x+1/2,-y,z) \end{array}$
\\
\hline
$-a2$ & $(2_12_1\times)=(*{:}\bar\times)=(\bar\times\times_1)$ & 33 & $Pa2_1$ &
$\begin{array}{c} (-x,y+1/2,z+1/2) \\ (x+1/2,-y,z) \end{array}$
\\
\hline
\end{tabular}
\bigskip
\caption{Notations for space groups.}
\end{table}

The two international numbers and names given in the three metachiral cases correspond to the two
orientations.

The infinite --- or non-compact --- platycosms were systematically treated probably for the first time
in \cite{Wo} p.123, with the symbols indicated in the corresponding column of the following table:

\begin{table}[!htb]\label{tinfiniteplaty}
\begin{tabular}{|c|r@{}l|c|}
\hline
our name            & \multicolumn{2}{|c|}{symbol}    & Wolf \\
\hline
Euclidean Space     &  & $EUC$            & $\m{E}$      \\
Circular Product space   &  & $CPS_A(\theta)$  & $\m{S}^\theta_1$ \\
Circular M\"obius space    &  & $CMS_A$          & $\m{S}_2$    \\
Toroidal Product space   &  & $TPS_{A\,B\,C}$  & $\m{T}_1$    \\
Toroidal M\"obius space    &  & $TMS_{A\,B:C}$   & $\m{T}_2$    \\
Kleinian  Product space  &  & $KPS^A_B$      & $\m{K}_2$    \\
chiral Kleinian M\"obius space & $+$ & $KMS^A_B$   & $\m{K}_1$    \\[.05cm]
achiral Kleinian M\"obius space & $-$ & $KMS^A_B$   & $\m{K}_3$    \\
\hline
\end{tabular}
\bigskip
\caption{Infinite platycosms.}
\end{table}

\

\end{document}

%% file: torklein.pstex_t
\begin{picture}(0,0)%
\includegraphics{torklein.pstex}%
\end{picture}%
\setlength{\unitlength}{1973sp}%
\begingroup\makeatletter\ifx\SetFigFont\undefined%
\gdef\SetFigFont#1#2#3#4#5{%
  \reset@font\fontsize{#1}{#2pt}%
  \fontfamily{#3}\fontseries{#4}\fontshape{#5}%
  \selectfont}%
\fi\endgroup%
\begin{picture}(11349,5892)(-611,-1819)
\put(301,2564){\makebox(0,0)[lb]{\smash{\SetFigFont{12}{14.4}{\rmdefault}{\mddefault}{\updefault}$\x$}}}
\put(2101,-1711){\makebox(0,0)[lb]{\smash{\SetFigFont{10}{12.0}{\rmdefault}{\mddefault}{\updefault}$B=-\w\cdot\y$}}}
\put(-599,314){\makebox(0,0)[lb]{\smash{\SetFigFont{10}{12.0}{\rmdefault}{\mddefault}{\updefault}$C=-\x\cdot\y$}}}
\put(1651,2564){\makebox(0,0)[lb]{\smash{\SetFigFont{10}{12.0}{\rmdefault}{\mddefault}{\updefault}$A=-\w\cdot\x$}}}
\put(1126,-1711){\makebox(0,0)[lb]{\smash{\SetFigFont{12}{14.4}{\rmdefault}{\mddefault}{\updefault}$\y$}}}
\put(2701,314){\makebox(0,0)[lb]{\smash{\SetFigFont{12}{14.4}{\rmdefault}{\mddefault}{\updefault}$\w$}}}
\put(9601,314){\makebox(0,0)[lb]{\smash{\SetFigFont{12}{14.4}{\rmdefault}{\mddefault}{\updefault}$\x$}}}
\put(9601, 14){\makebox(0,0)[lb]{\smash{\SetFigFont{10}{12.0}{\rmdefault}{\mddefault}{\updefault}$N(\x)=A$}}}
\put(8851,2564){\makebox(0,0)[lb]{\smash{\SetFigFont{10}{12.0}{\rmdefault}{\mddefault}{\updefault}$N(\y)=B$}}}
\put(8101,2564){\makebox(0,0)[lb]{\smash{\SetFigFont{12}{14.4}{\rmdefault}{\mddefault}{\updefault}$\y$}}}
\end{picture}

%% file: contor.pstex_t
\begin{picture}(0,0)%
\includegraphics{contor.pstex}%
\end{picture}%
\setlength{\unitlength}{1973sp}%
\begingroup\makeatletter\ifx\SetFigFont\undefined%
\gdef\SetFigFont#1#2#3#4#5{%
  \reset@font\fontsize{#1}{#2pt}%
  \fontfamily{#3}\fontseries{#4}\fontshape{#5}%
  \selectfont}%
\fi\endgroup%
\begin{picture}(12075,2701)(526,-2294)
\put(8101,-961){\makebox(0,0)[lb]{\smash{\SetFigFont{8}{9.6}{\rmdefault}{\mddefault}{\updefault}$0$}}}
\put(7351,-2236){\makebox(0,0)[lb]{\smash{\SetFigFont{8}{9.6}{\rmdefault}{\mddefault}{\updefault}$B$}}}
\put(8701,-2236){\makebox(0,0)[lb]{\smash{\SetFigFont{8}{9.6}{\rmdefault}{\mddefault}{\updefault}$C$}}}
\put(7261,-1181){\makebox(0,0)[lb]{\smash{\SetFigFont{8}{9.6}{\rmdefault}{\mddefault}{\updefault}$0$}}}
\put(7366,164){\makebox(0,0)[lb]{\smash{\SetFigFont{8}{9.6}{\rmdefault}{\mddefault}{\updefault}$0$}}}
\put(6001,239){\makebox(0,0)[lb]{\smash{\SetFigFont{8}{9.6}{\familydefault}{\mddefault}{\updefault}(ii)}}}
\put(12001,-961){\makebox(0,0)[lb]{\smash{\SetFigFont{8}{9.6}{\rmdefault}{\mddefault}{\updefault}$0$}}}
\put(11251,-2236){\makebox(0,0)[lb]{\smash{\SetFigFont{8}{9.6}{\rmdefault}{\mddefault}{\updefault}$B$}}}
\put(12601,-2236){\makebox(0,0)[lb]{\smash{\SetFigFont{8}{9.6}{\rmdefault}{\mddefault}{\updefault}$C$}}}
\put(10501,-961){\makebox(0,0)[lb]{\smash{\SetFigFont{8}{9.6}{\rmdefault}{\mddefault}{\updefault}$0$}}}
\put(11161,-1181){\makebox(0,0)[lb]{\smash{\SetFigFont{8}{9.6}{\rmdefault}{\mddefault}{\updefault}$0$}}}
\put(11266,164){\makebox(0,0)[lb]{\smash{\SetFigFont{8}{9.6}{\rmdefault}{\mddefault}{\updefault}$D$}}}
\put(9901,239){\makebox(0,0)[lb]{\smash{\SetFigFont{8}{9.6}{\familydefault}{\mddefault}{\updefault}(iii)}}}
\put(9226,-961){\makebox(0,0)[lb]{\smash{\SetFigFont{10}{12.0}{\rmdefault}{\mddefault}{\updefault}$\cong$}}}
\put(2701,-961){\makebox(0,0)[lb]{\smash{\SetFigFont{8}{9.6}{\rmdefault}{\mddefault}{\updefault}$F$}}}
\put(1951,-2236){\makebox(0,0)[lb]{\smash{\SetFigFont{8}{9.6}{\rmdefault}{\mddefault}{\updefault}$B$}}}
\put(3301,-2236){\makebox(0,0)[lb]{\smash{\SetFigFont{8}{9.6}{\rmdefault}{\mddefault}{\updefault}$C$}}}
\put(1966,164){\makebox(0,0)[lb]{\smash{\SetFigFont{8}{9.6}{\rmdefault}{\mddefault}{\updefault}$0$}}}
\put(601,239){\makebox(0,0)[lb]{\smash{\SetFigFont{8}{9.6}{\familydefault}{\mddefault}{\updefault}(i)}}}
\put(1126,-961){\makebox(0,0)[lb]{\smash{\SetFigFont{8}{9.6}{\rmdefault}{\mddefault}{\updefault}$D$}}}
\put(1801,-1186){\makebox(0,0)[lb]{\smash{\SetFigFont{8}{9.6}{\rmdefault}{\mddefault}{\updefault}$E$}}}
\put(6526,-961){\makebox(0,0)[lb]{\smash{\SetFigFont{8}{9.6}{\rmdefault}{\mddefault}{\updefault}$D$}}}
\put(526,-2236){\makebox(0,0)[lb]{\smash{\SetFigFont{8}{9.6}{\rmdefault}{\mddefault}{\updefault}$A$}}}
\put(5926,-2236){\makebox(0,0)[lb]{\smash{\SetFigFont{8}{9.6}{\rmdefault}{\mddefault}{\updefault}$A$}}}
\put(9826,-2236){\makebox(0,0)[lb]{\smash{\SetFigFont{8}{9.6}{\rmdefault}{\mddefault}{\updefault}$A$}}}
\end{picture}

%% file: triang.pstex_t
\begin{picture}(0,0)%
\includegraphics{triang.pstex}%
\end{picture}%
\setlength{\unitlength}{1973sp}%
\begingroup\makeatletter\ifx\SetFigFont\undefined%
\gdef\SetFigFont#1#2#3#4#5{%
  \reset@font\fontsize{#1}{#2pt}%
  \fontfamily{#3}\fontseries{#4}\fontshape{#5}%
  \selectfont}%
\fi\endgroup%
\begin{picture}(10858,3570)(1,-3718)
\put(  1,-3436){\makebox(0,0)[lb]{\smash{\SetFigFont{12}{14.4}{\rmdefault}{\mddefault}{\updefault}$\Ld_{A\,A\,A}$}}}
\put(1276,-586){\makebox(0,0)[lb]{\smash{\SetFigFont{12}{14.4}{\rmdefault}{\mddefault}{\updefault}$\x$}}}
\put(1351,-2611){\makebox(0,0)[lb]{\smash{\SetFigFont{12}{14.4}{\rmdefault}{\mddefault}{\updefault}$\y$}}}
\put(7876,-436){\makebox(0,0)[lb]{\smash{\SetFigFont{12}{14.4}{\rmdefault}{\mddefault}{\updefault}$\y$}}}
\put(9526,-1861){\makebox(0,0)[lb]{\smash{\SetFigFont{12}{14.4}{\rmdefault}{\mddefault}{\updefault}$\x$}}}
\put(3526,-1861){\makebox(0,0)[lb]{\smash{\SetFigFont{12}{14.4}{\rmdefault}{\mddefault}{\updefault}$\w$}}}
\put(7726,-2911){\makebox(0,0)[lb]{\smash{\SetFigFont{12}{14.4}{\rmdefault}{\mddefault}{\updefault}$-\y$}}}
\put(7051,-1861){\makebox(0,0)[lb]{\smash{\SetFigFont{12}{14.4}{\rmdefault}{\mddefault}{\updefault}$-\x$}}}
\put(2251,-2611){\makebox(0,0)[lb]{\smash{\SetFigFont{12}{14.4}{\rmdefault}{\mddefault}{\updefault}$-\x$}}}
\put(2251,-586){\makebox(0,0)[lb]{\smash{\SetFigFont{12}{14.4}{\rmdefault}{\mddefault}{\updefault}$-\y$}}}
\put(901,-1861){\makebox(0,0)[lb]{\smash{\SetFigFont{12}{14.4}{\rmdefault}{\mddefault}{\updefault}$-\w$}}}
\put(6901,-3511){\makebox(0,0)[lb]{\smash{\SetFigFont{12}{14.4}{\rmdefault}{\mddefault}{\updefault}$\Ld_{A\,A}$}}}
\end{picture}

%% file: xfigamphic.pstex_t
\begin{picture}(0,0)%
\includegraphics{xfigamphic.pstex}%
\end{picture}%
\setlength{\unitlength}{3947sp}%
\begingroup\makeatletter\ifx\SetFigFont\undefined%
\gdef\SetFigFont#1#2#3#4#5{%
  \reset@font\fontsize{#1}{#2pt}%
  \fontfamily{#3}\fontseries{#4}\fontshape{#5}%
  \selectfont}%
\fi\endgroup%
\begin{picture}(6424,2103)(751,-2069)
\put(751,-1102){\makebox(0,0)[lb]{\smash{\SetFigFont{14}{16.8}{\rmdefault}{\mddefault}{\updefault}$\x$}}}
\put(1602,-138){\makebox(0,0)[lb]{\smash{\SetFigFont{14}{16.8}{\rmdefault}{\mddefault}{\updefault}$\w$}}}
\put(5376,-110){\makebox(0,0)[lb]{\smash{\SetFigFont{14}{16.8}{\rmdefault}{\mddefault}{\updefault}$\w$}}}
\put(4468,-1102){\makebox(0,0)[lb]{\smash{\SetFigFont{14}{16.8}{\rmdefault}{\mddefault}{\updefault}$\x$}}}
\put(1426,-2011){\makebox(0,0)[lb]{\smash{\SetFigFont{14}{16.8}{\rmdefault}{\mddefault}{\updefault}$\z$}}}
\put(5176,-2011){\makebox(0,0)[lb]{\smash{\SetFigFont{14}{16.8}{\rmdefault}{\mddefault}{\updefault}$\z$}}}
\put(1576,-2011){\makebox(0,0)[lb]{\smash{\SetFigFont{12}{14.4}{\rmdefault}{\mddefault}{\updefault}(norm $D$)}}}
\put(5326,-2011){\makebox(0,0)[lb]{\smash{\SetFigFont{12}{14.4}{\rmdefault}{\mddefault}{\updefault}(norm $D$)}}}
\end{picture}

%% file: xfigamphid.pstex_t
\begin{picture}(0,0)%
\includegraphics{xfigamphid.pstex}%
\end{picture}%
\setlength{\unitlength}{3552sp}%
\begingroup\makeatletter\ifx\SetFigFont\undefined%
\gdef\SetFigFont#1#2#3#4#5{%
  \reset@font\fontsize{#1}{#2pt}%
  \fontfamily{#3}\fontseries{#4}\fontshape{#5}%
  \selectfont}%
\fi\endgroup%
\begin{picture}(6495,2192)(751,-2219)
\put(751,-1207){\makebox(0,0)[lb]{\smash{\SetFigFont{12}{14.4}{\rmdefault}{\mddefault}{\updefault}$\x$}}}
\put(1645,-195){\makebox(0,0)[lb]{\smash{\SetFigFont{12}{14.4}{\rmdefault}{\mddefault}{\updefault}$\y$}}}
\put(4403,-1207){\makebox(0,0)[lb]{\smash{\SetFigFont{12}{14.4}{\rmdefault}{\mddefault}{\updefault}$\x$}}}
\put(5327,-195){\makebox(0,0)[lb]{\smash{\SetFigFont{12}{14.4}{\rmdefault}{\mddefault}{\updefault}$\y$}}}
\put(5326,-2161){\makebox(0,0)[lb]{\smash{\SetFigFont{12}{14.4}{\rmdefault}{\mddefault}{\updefault}$\z$}}}
\put(5551,-2161){\makebox(0,0)[lb]{\smash{\SetFigFont{11}{13.2}{\rmdefault}{\mddefault}{\updefault}(norm $D$)}}}
\put(1801,-2161){\makebox(0,0)[lb]{\smash{\SetFigFont{11}{13.2}{\rmdefault}{\mddefault}{\updefault}(norm $D$)}}}
\put(1576,-2161){\makebox(0,0)[lb]{\smash{\SetFigFont{12}{14.4}{\rmdefault}{\mddefault}{\updefault}$\z$}}}
\end{picture}

%% file: planesc2.pstex_t
\begin{picture}(0,0)%
\includegraphics{planesc2.pstex}%
\end{picture}%
\setlength{\unitlength}{1579sp}%
\begingroup\makeatletter\ifx\SetFigFont\undefined%
\gdef\SetFigFont#1#2#3#4#5{%
  \reset@font\fontsize{#1}{#2pt}%
  \fontfamily{#3}\fontseries{#4}\fontshape{#5}%
  \selectfont}%
\fi\endgroup%
\begin{picture}(6012,3795)(1,-3421)
\put(3151, 14){\makebox(0,0)[lb]{\smash{\SetFigFont{12}{14.4}{\rmdefault}{\mddefault}{\updefault}$\mb{y}$}}}
\put(2101,-961){\makebox(0,0)[lb]{\smash{\SetFigFont{12}{14.4}{\rmdefault}{\mddefault}{\updefault}$\frac12\mb{y}$}}}
\put(5701,-3286){\makebox(0,0)[lb]{\smash{\SetFigFont{12}{14.4}{\rmdefault}{\mddefault}{\updefault}$\mb{x}$}}}
\put(  1,-1336){\makebox(0,0)[lb]{\smash{\SetFigFont{12}{14.4}{\rmdefault}{\mddefault}{\updefault}$1sK$}}}
\put(  1,-2911){\makebox(0,0)[lb]{\smash{\SetFigFont{12}{14.4}{\rmdefault}{\mddefault}{\updefault}$1sK$}}}
\put(  1,-2161){\makebox(0,0)[lb]{\smash{\SetFigFont{12}{14.4}{\rmdefault}{\mddefault}{\updefault}$2T$}}}
\end{picture}

%% file: xfigamphis.pstex_t
\begin{picture}(0,0)%
\includegraphics{xfigamphis.pstex}%
\end{picture}%
\setlength{\unitlength}{3947sp}%
\begingroup\makeatletter\ifx\SetFigFont\undefined%
\gdef\SetFigFont#1#2#3#4#5{%
  \reset@font\fontsize{#1}{#2pt}%
  \fontfamily{#3}\fontseries{#4}\fontshape{#5}%
  \selectfont}%
\fi\endgroup%
\begin{picture}(6687,1831)(751,-3785)
\put(1576,-2986){\makebox(0,0)[lb]{\smash{\SetFigFont{14}{16.8}{\rmdefault}{\mddefault}{\updefault}{\bf o}}}}
\put(1576,-3286){\makebox(0,0)[lb]{\smash{\SetFigFont{14}{16.8}{\rmdefault}{\mddefault}{\updefault}{\bf o}}}}
\put(1276,-3286){\makebox(0,0)[lb]{\smash{\SetFigFont{14}{16.8}{\rmdefault}{\mddefault}{\updefault}{\bf d}}}}
\put(1951,-2461){\makebox(0,0)[lb]{\smash{\SetFigFont{14}{16.8}{\rmdefault}{\mddefault}{\updefault}{\bf o}}}}
\put(1951,-2761){\makebox(0,0)[lb]{\smash{\SetFigFont{14}{16.8}{\rmdefault}{\mddefault}{\updefault}{\bf o}}}}
\put(1651,-2761){\makebox(0,0)[lb]{\smash{\SetFigFont{14}{16.8}{\rmdefault}{\mddefault}{\updefault}{\bf b}}}}
\put(1276,-2986){\makebox(0,0)[lb]{\smash{\SetFigFont{14}{16.8}{\rmdefault}{\mddefault}{\updefault}{\bf b}}}}
\put(5176,-2986){\makebox(0,0)[lb]{\smash{\SetFigFont{14}{16.8}{\rmdefault}{\mddefault}{\updefault}{\bf o}}}}
\put(5176,-3286){\makebox(0,0)[lb]{\smash{\SetFigFont{14}{16.8}{\rmdefault}{\mddefault}{\updefault}{\bf o}}}}
\put(4876,-3286){\makebox(0,0)[lb]{\smash{\SetFigFont{14}{16.8}{\rmdefault}{\mddefault}{\updefault}{\bf d}}}}
\put(5551,-2461){\makebox(0,0)[lb]{\smash{\SetFigFont{14}{16.8}{\rmdefault}{\mddefault}{\updefault}{\bf o}}}}
\put(5551,-2761){\makebox(0,0)[lb]{\smash{\SetFigFont{14}{16.8}{\rmdefault}{\mddefault}{\updefault}{\bf o}}}}
\put(5251,-2461){\makebox(0,0)[lb]{\smash{\SetFigFont{14}{16.8}{\rmdefault}{\mddefault}{\updefault}{\bf q}}}}
\put(5251,-2761){\makebox(0,0)[lb]{\smash{\SetFigFont{14}{16.8}{\rmdefault}{\mddefault}{\updefault}{\bf p}}}}
\put(4876,-2986){\makebox(0,0)[lb]{\smash{\SetFigFont{14}{16.8}{\rmdefault}{\mddefault}{\updefault}{\bf b}}}}
\put(3226,-2986){\makebox(0,0)[lb]{\smash{\SetFigFont{14}{16.8}{\rmdefault}{\mddefault}{\updefault}{\bf o}}}}
\put(3226,-3286){\makebox(0,0)[lb]{\smash{\SetFigFont{14}{16.8}{\rmdefault}{\mddefault}{\updefault}{\bf o}}}}
\put(2926,-2986){\makebox(0,0)[lb]{\smash{\SetFigFont{14}{16.8}{\rmdefault}{\mddefault}{\updefault}{\bf b}}}}
\put(2926,-3286){\makebox(0,0)[lb]{\smash{\SetFigFont{14}{16.8}{\rmdefault}{\mddefault}{\updefault}{\bf d}}}}
\put(3601,-2461){\makebox(0,0)[lb]{\smash{\SetFigFont{14}{16.8}{\rmdefault}{\mddefault}{\updefault}{\bf d}}}}
\put(3601,-2761){\makebox(0,0)[lb]{\smash{\SetFigFont{14}{16.8}{\rmdefault}{\mddefault}{\updefault}{\bf b}}}}
\put(3301,-2461){\makebox(0,0)[lb]{\smash{\SetFigFont{14}{16.8}{\rmdefault}{\mddefault}{\updefault}{\bf o}}}}
\put(3301,-2761){\makebox(0,0)[lb]{\smash{\SetFigFont{14}{16.8}{\rmdefault}{\mddefault}{\updefault}{\bf o}}}}
\put(6826,-2986){\makebox(0,0)[lb]{\smash{\SetFigFont{14}{16.8}{\rmdefault}{\mddefault}{\updefault}{\bf o}}}}
\put(6826,-3286){\makebox(0,0)[lb]{\smash{\SetFigFont{14}{16.8}{\rmdefault}{\mddefault}{\updefault}{\bf o}}}}
\put(6526,-2986){\makebox(0,0)[lb]{\smash{\SetFigFont{14}{16.8}{\rmdefault}{\mddefault}{\updefault}{\bf b}}}}
\put(6526,-3286){\makebox(0,0)[lb]{\smash{\SetFigFont{14}{16.8}{\rmdefault}{\mddefault}{\updefault}{\bf d}}}}
\put(7201,-2461){\makebox(0,0)[lb]{\smash{\SetFigFont{14}{16.8}{\rmdefault}{\mddefault}{\updefault}{\bf q}}}}
\put(7201,-2761){\makebox(0,0)[lb]{\smash{\SetFigFont{14}{16.8}{\rmdefault}{\mddefault}{\updefault}{\bf p}}}}
\put(6901,-2761){\makebox(0,0)[lb]{\smash{\SetFigFont{14}{16.8}{\rmdefault}{\mddefault}{\updefault}{\bf o}}}}
\put(6901,-2461){\makebox(0,0)[lb]{\smash{\SetFigFont{14}{16.8}{\rmdefault}{\mddefault}{\updefault}{\bf o}}}}
\put(1651,-2461){\makebox(0,0)[lb]{\smash{\SetFigFont{14}{16.8}{\rmdefault}{\mddefault}{\updefault}{\bf d}}}}
\put(1426,-3736){\makebox(0,0)[lb]{\smash{\SetFigFont{12}{14.4}{\rmdefault}{\mddefault}{\updefault}$Z$}}}
\put(3076,-3736){\makebox(0,0)[lb]{\smash{\SetFigFont{12}{14.4}{\rmdefault}{\mddefault}{\updefault}$Z$}}}
\put(5026,-3736){\makebox(0,0)[lb]{\smash{\SetFigFont{12}{14.4}{\rmdefault}{\mddefault}{\updefault}$Z$}}}
\put(6676,-3736){\makebox(0,0)[lb]{\smash{\SetFigFont{12}{14.4}{\rmdefault}{\mddefault}{\updefault}$Z$}}}
\put(751,-3136){\makebox(0,0)[lb]{\smash{\SetFigFont{12}{14.4}{\rmdefault}{\mddefault}{\updefault}$X$}}}
\put(2401,-3136){\makebox(0,0)[lb]{\smash{\SetFigFont{12}{14.4}{\rmdefault}{\mddefault}{\updefault}$X$}}}
\put(4351,-3136){\makebox(0,0)[lb]{\smash{\SetFigFont{12}{14.4}{\rmdefault}{\mddefault}{\updefault}$X$}}}
\put(6001,-3136){\makebox(0,0)[lb]{\smash{\SetFigFont{12}{14.4}{\rmdefault}{\mddefault}{\updefault}$X$}}}
\put(1201,-2686){\makebox(0,0)[lb]{\smash{\SetFigFont{12}{14.4}{\rmdefault}{\mddefault}{\updefault}$W$}}}
\put(2926,-2686){\makebox(0,0)[lb]{\smash{\SetFigFont{12}{14.4}{\rmdefault}{\mddefault}{\updefault}$W$}}}
\put(4801,-2686){\makebox(0,0)[lb]{\smash{\SetFigFont{12}{14.4}{\rmdefault}{\mddefault}{\updefault}$W$}}}
\put(6526,-2686){\makebox(0,0)[lb]{\smash{\SetFigFont{12}{14.4}{\rmdefault}{\mddefault}{\updefault}$W$}}}
\put(1201,-2086){\makebox(0,0)[lb]{\smash{\SetFigFont{12}{14.4}{\rmdefault}{\mddefault}{\updefault}$+a1$:}}}
\put(2851,-2086){\makebox(0,0)[lb]{\smash{\SetFigFont{12}{14.4}{\rmdefault}{\mddefault}{\updefault}$-a1$:}}}
\put(4801,-2086){\makebox(0,0)[lb]{\smash{\SetFigFont{12}{14.4}{\rmdefault}{\mddefault}{\updefault}$+a2$:}}}
\put(6451,-2086){\makebox(0,0)[lb]{\smash{\SetFigFont{12}{14.4}{\rmdefault}{\mddefault}{\updefault}$-a2$:}}}
\end{picture}

%% file: symmetries.pstex_t
\begin{picture}(0,0)%
\includegraphics{symmetries.pstex}%
\end{picture}%
\setlength{\unitlength}{2960sp}%
\begingroup\makeatletter\ifx\SetFigFont\undefined%
\gdef\SetFigFont#1#2#3#4#5{%
  \reset@font\fontsize{#1}{#2pt}%
  \fontfamily{#3}\fontseries{#4}\fontshape{#5}%
  \selectfont}%
\fi\endgroup%
\begin{picture}(7598,3824)(-890,-2548)
\put(2647,-1711){\makebox(0,0)[lb]{\smash{\SetFigFont{17}{20.4}{\rmdefault}{\mddefault}{\updefault}$\m{T}$}}}
\put(1951,464){\makebox(0,0)[lb]{\smash{\SetFigFont{14}{16.8}{\familydefault}{\mddefault}{\updefault}$B_1$}}}
\put(2701,539){\makebox(0,0)[lb]{\smash{\SetFigFont{14}{16.8}{\familydefault}{\mddefault}{\updefault}$B_2$}}}
\end{picture}

%% file: con-von+a1.pstex_t
\begin{picture}(0,0)%
\includegraphics{con-von+a1.pstex}%
\end{picture}%
\setlength{\unitlength}{2131sp}%
\begingroup\makeatletter\ifx\SetFigFont\undefined%
\gdef\SetFigFont#1#2#3#4#5{%
  \reset@font\fontsize{#1}{#2pt}%
  \fontfamily{#3}\fontseries{#4}\fontshape{#5}%
  \selectfont}%
\fi\endgroup%
\begin{picture}(9538,2689)(-449,-2294)
\put(9076,164){\makebox(0,0)[lb]{\smash{\SetFigFont{9}{10.8}{\rmdefault}{\mddefault}{\updefault}$A\!+\!B$}}}
\put(3386,-2226){\makebox(0,0)[lb]{\smash{\SetFigFont{9}{10.8}{\rmdefault}{\mddefault}{\updefault}$C$}}}
\put(2786,-981){\makebox(0,0)[lb]{\smash{\SetFigFont{9}{10.8}{\rmdefault}{\mddefault}{\updefault}$0$}}}
\put(1971,164){\makebox(0,0)[lb]{\smash{\SetFigFont{9}{10.8}{\rmdefault}{\mddefault}{\updefault}$D$}}}
\put(601,-2236){\makebox(0,0)[lb]{\smash{\SetFigFont{9}{10.8}{\rmdefault}{\mddefault}{\updefault}$A$}}}
\put(1951,-2236){\makebox(0,0)[lb]{\smash{\SetFigFont{9}{10.8}{\rmdefault}{\mddefault}{\updefault}$B$}}}
\put(-449,-586){\makebox(0,0)[lb]{\smash{\SetFigFont{11}{13.2}{\rmdefault}{\mddefault}{\updefault}\text{conorms:}}}}
\put(6001,164){\makebox(0,0)[lb]{\smash{\SetFigFont{9}{10.8}{\rmdefault}{\mddefault}{\updefault}$B\!+\!C$}}}
\put(7501,239){\makebox(0,0)[lb]{\smash{\SetFigFont{9}{10.8}{\rmdefault}{\mddefault}{\updefault}$A\!+\!C$}}}
\put(7726,-2236){\makebox(0,0)[lb]{\smash{\SetFigFont{9}{10.8}{\rmdefault}{\mddefault}{\updefault}$D$}}}
\put(1876,-1111){\makebox(0,0)[lb]{\smash{\SetFigFont{9}{10.8}{\rmdefault}{\mddefault}{\updefault}$0$}}}
\put(1276,-961){\makebox(0,0)[lb]{\smash{\SetFigFont{9}{10.8}{\rmdefault}{\mddefault}{\updefault}$0$}}}
\put(4426,-586){\makebox(0,0)[lb]{\smash{\SetFigFont{11}{13.2}{\rmdefault}{\mddefault}{\updefault}\text{vonorms:}}}}
\put(6001,-1111){\makebox(0,0)[lb]{\smash{\SetFigFont{9}{10.8}{\rmdefault}{\mddefault}{\updefault}$B\!+\!C\!+\!D$}}}
\put(8551,-1036){\makebox(0,0)[lb]{\smash{\SetFigFont{9}{10.8}{\rmdefault}{\mddefault}{\updefault}$A\!+\!B\!+\!D$}}}
\put(7276,-501){\makebox(0,0)[lb]{\smash{\SetFigFont{9}{10.8}{\rmdefault}{\mddefault}{\updefault}$A\!+\!C\!+\!D$}}}
\end{picture}

%% file: kb.pstex_t
\begin{picture}(0,0)%
\includegraphics{kb.pstex}%
\end{picture}%
\setlength{\unitlength}{1381sp}%
\begingroup\makeatletter\ifx\SetFigFont\undefined%
\gdef\SetFigFont#1#2#3#4#5{%
  \reset@font\fontsize{#1}{#2pt}%
  \fontfamily{#3}\fontseries{#4}\fontshape{#5}%
  \selectfont}%
\fi\endgroup%
\begin{picture}(12024,6011)(1189,-6347)
\put(4951,-4501){\makebox(0,0)[lb]{\smash{\SetFigFont{11}{13.2}{\rmdefault}{\mddefault}{\updefault}$\mb{p_0}$}}}
\put(3601,-4005){\makebox(0,0)[lb]{\smash{\SetFigFont{11}{13.2}{\rmdefault}{\mddefault}{\updefault}$\mb{p_1}$}}}
\end{picture}

%% file: con-von.pstex_t
\begin{picture}(0,0)%
\includegraphics{con-von.pstex}%
\end{picture}%
\setlength{\unitlength}{2171sp}%
\begingroup\makeatletter\ifx\SetFigFont\undefined%
\gdef\SetFigFont#1#2#3#4#5{%
  \reset@font\fontsize{#1}{#2pt}%
  \fontfamily{#3}\fontseries{#4}\fontshape{#5}%
  \selectfont}%
\fi\endgroup%
\begin{picture}(9538,2689)(-449,-2294)
\put(3451,-886){\makebox(0,0)[lb]{\smash{\SetFigFont{11}{13.2}{\rmdefault}{\mddefault}{\updefault}\text{has vonorms:}}}}
\put(3386,-2226){\makebox(0,0)[lb]{\smash{\SetFigFont{9}{10.8}{\rmdefault}{\mddefault}{\updefault}$C$}}}
\put(2786,-981){\makebox(0,0)[lb]{\smash{\SetFigFont{9}{10.8}{\rmdefault}{\mddefault}{\updefault}$F$}}}
\put(1971,164){\makebox(0,0)[lb]{\smash{\SetFigFont{9}{10.8}{\rmdefault}{\mddefault}{\updefault}$G$}}}
\put(601,-2236){\makebox(0,0)[lb]{\smash{\SetFigFont{9}{10.8}{\rmdefault}{\mddefault}{\updefault}$A$}}}
\put(1801,-1111){\makebox(0,0)[lb]{\smash{\SetFigFont{9}{10.8}{\rmdefault}{\mddefault}{\updefault}$E$}}}
\put(1951,-2236){\makebox(0,0)[lb]{\smash{\SetFigFont{9}{10.8}{\rmdefault}{\mddefault}{\updefault}$B$}}}
\put(1201,-961){\makebox(0,0)[lb]{\smash{\SetFigFont{9}{10.8}{\rmdefault}{\mddefault}{\updefault}$D$}}}
\put(8626,-1036){\makebox(0,0)[lb]{\smash{\SetFigFont{9}{10.8}{\rmdefault}{\mddefault}{\updefault}$A\!+\!B\!+\!F\!+\!G$}}}
\put(7126,239){\makebox(0,0)[lb]{\smash{\SetFigFont{9}{10.8}{\rmdefault}{\mddefault}{\updefault}$C\!+\!F\!+\!A\!+\!D$}}}
\put(-449,-886){\makebox(0,0)[lb]{\smash{\SetFigFont{11}{13.2}{\rmdefault}{\mddefault}{\updefault}\text{conorms:}}}}
\put(9076,164){\makebox(0,0)[lb]{\smash{\SetFigFont{9}{10.8}{\rmdefault}{\mddefault}{\updefault}$A\!+\!D\!+\!B\!+\!E$}}}
\put(5326,164){\makebox(0,0)[lb]{\smash{\SetFigFont{9}{10.8}{\rmdefault}{\mddefault}{\updefault}$B\!+\!E\!+\!C\!+\!F$}}}
\put(5626,-1111){\makebox(0,0)[lb]{\smash{\SetFigFont{9}{10.8}{\rmdefault}{\mddefault}{\updefault}$D\!+\!B\!+\!C\!+\!G$}}}
\put(7201,-511){\makebox(0,0)[lb]{\smash{\SetFigFont{9}{10.8}{\rmdefault}{\mddefault}{\updefault}$A\!+\!E\!+\!C\!+\!G$}}}
\put(7051,-2236){\makebox(0,0)[lb]{\smash{\SetFigFont{9}{10.8}{\rmdefault}{\mddefault}{\updefault}$D\!+\!E\!+\!F\!+\!G$}}}
\end{picture}

%% file: sensual.pstex_t
\begin{picture}(0,0)%
\includegraphics{sensual.pstex}%
\end{picture}%
\setlength{\unitlength}{3158sp}%
\begingroup\makeatletter\ifx\SetFigFont\undefined%
\gdef\SetFigFont#1#2#3#4#5{%
  \reset@font\fontsize{#1}{#2pt}%
  \fontfamily{#3}\fontseries{#4}\fontshape{#5}%
  \selectfont}%
\fi\endgroup%
\begin{picture}(6919,7690)(-2324,-6707)
\put(2851,-670){\makebox(0,0)[lb]{\smash{\SetFigFont{10}{12.0}{\familydefault}{\mddefault}{\updefault}$\al$}}}
\put(3601,-725){\makebox(0,0)[lb]{\smash{\SetFigFont{10}{12.0}{\familydefault}{\mddefault}{\updefault}$\be$}}}
\put(4360,-670){\makebox(0,0)[lb]{\smash{\SetFigFont{10}{12.0}{\familydefault}{\mddefault}{\updefault}$\g$}}}
\put(3141, 14){\makebox(0,0)[lb]{\smash{\SetFigFont{10}{12.0}{\familydefault}{\mddefault}{\updefault}$a$}}}
\put(3621,665){\makebox(0,0)[lb]{\smash{\SetFigFont{10}{12.0}{\familydefault}{\mddefault}{\updefault}$0$}}}
\put(4066, 14){\makebox(0,0)[lb]{\smash{\SetFigFont{10}{12.0}{\familydefault}{\mddefault}{\updefault}$c$}}}
\put(3541,-61){\makebox(0,0)[lb]{\smash{\SetFigFont{10}{12.0}{\familydefault}{\mddefault}{\updefault}$b$}}}
\put(-1499,-61){\makebox(0,0)[lb]{\smash{\SetFigFont{11}{13.2}{\familydefault}{\mddefault}{\updefault}truncated}}}
\put(-1499,-286){\makebox(0,0)[lb]{\smash{\SetFigFont{11}{13.2}{\familydefault}{\mddefault}{\updefault}Octahedron}}}
\put(  1,-5911){\makebox(0,0)[lb]{\smash{\SetFigFont{11}{13.2}{\familydefault}{\mddefault}{\updefault}rectangular}}}
\put(  1,-6136){\makebox(0,0)[lb]{\smash{\SetFigFont{11}{13.2}{\familydefault}{\mddefault}{\updefault}Cuboid}}}
\put(-899,-2386){\makebox(0,0)[lb]{\smash{\SetFigFont{11}{13.2}{\familydefault}{\mddefault}{\updefault}hexarhombic}}}
\put(-899,-2611){\makebox(0,0)[lb]{\smash{\SetFigFont{11}{13.2}{\familydefault}{\mddefault}{\updefault}Dodecahedron}}}
\put(  1,-3736){\makebox(0,0)[lb]{\smash{\SetFigFont{10}{12.0}{\familydefault}{\mddefault}{\updefault}$b=0$}}}
\put(1801,-3736){\makebox(0,0)[lb]{\smash{\SetFigFont{10}{12.0}{\familydefault}{\mddefault}{\updefault}$\al=0$}}}
\put(601,-4861){\makebox(0,0)[lb]{\smash{\SetFigFont{10}{12.0}{\familydefault}{\mddefault}{\updefault}$\al=0$}}}
\put(1576,-4861){\makebox(0,0)[lb]{\smash{\SetFigFont{10}{12.0}{\familydefault}{\mddefault}{\updefault}$b=0$}}}
\put(1801,-3061){\makebox(0,0)[lb]{\smash{\SetFigFont{10}{12.0}{\familydefault}{\mddefault}{\updefault}$b$}}}
\put(1126,839){\makebox(0,0)[lb]{\smash{\SetFigFont{11}{13.2}{\familydefault}{\mddefault}{\updefault}$\be$}}}
\put(2051,239){\makebox(0,0)[lb]{\smash{\SetFigFont{11}{13.2}{\familydefault}{\mddefault}{\updefault}$\g$}}}
\put(151,-811){\makebox(0,0)[lb]{\smash{\SetFigFont{11}{13.2}{\familydefault}{\mddefault}{\updefault}$\g$}}}
\put(1126,-1111){\makebox(0,0)[lb]{\smash{\SetFigFont{11}{13.2}{\familydefault}{\mddefault}{\updefault}$\be$}}}
\put(511, 89){\makebox(0,0)[lb]{\smash{\SetFigFont{11}{13.2}{\familydefault}{\mddefault}{\updefault}$\al$}}}
\put(526,389){\makebox(0,0)[lb]{\smash{\SetFigFont{11}{13.2}{\familydefault}{\mddefault}{\updefault}$a$}}}
\put(1126,464){\makebox(0,0)[lb]{\smash{\SetFigFont{11}{13.2}{\familydefault}{\mddefault}{\updefault}$\be$}}}
\put(1751, 89){\makebox(0,0)[lb]{\smash{\SetFigFont{11}{13.2}{\familydefault}{\mddefault}{\updefault}$\g$}}}
\put(2271,-286){\makebox(0,0)[lb]{\smash{\SetFigFont{11}{13.2}{\familydefault}{\mddefault}{\updefault}$b$}}}
\put(1726,-661){\makebox(0,0)[lb]{\smash{\SetFigFont{11}{13.2}{\familydefault}{\mddefault}{\updefault}$\al$}}}
\put(1576,-1036){\makebox(0,0)[lb]{\smash{\SetFigFont{11}{13.2}{\familydefault}{\mddefault}{\updefault}$b$}}}
\put(526,-586){\makebox(0,0)[lb]{\smash{\SetFigFont{11}{13.2}{\familydefault}{\mddefault}{\updefault}$\g$}}}
\put(  1,-286){\makebox(0,0)[lb]{\smash{\SetFigFont{11}{13.2}{\familydefault}{\mddefault}{\updefault}$b$}}}
\put(151,239){\makebox(0,0)[lb]{\smash{\SetFigFont{11}{13.2}{\familydefault}{\mddefault}{\updefault}$\al$}}}
\put(301,-61){\makebox(0,0)[lb]{\smash{\SetFigFont{11}{13.2}{\familydefault}{\mddefault}{\updefault}$a$}}}
\put(1671,689){\makebox(0,0)[lb]{\smash{\SetFigFont{11}{13.2}{\familydefault}{\mddefault}{\updefault}$a$}}}
\put(551,709){\makebox(0,0)[lb]{\smash{\SetFigFont{11}{13.2}{\familydefault}{\mddefault}{\updefault}$c$}}}
\put(1126,-811){\makebox(0,0)[lb]{\smash{\SetFigFont{11}{13.2}{\familydefault}{\mddefault}{\updefault}$\be$}}}
\put(2051,-811){\makebox(0,0)[lb]{\smash{\SetFigFont{11}{13.2}{\familydefault}{\mddefault}{\updefault}$\al$}}}
\put(526,-1186){\makebox(0,0)[lb]{\smash{\SetFigFont{11}{13.2}{\familydefault}{\mddefault}{\updefault}$a$}}}
\put(731,-1036){\makebox(0,0)[lb]{\smash{\SetFigFont{11}{13.2}{\familydefault}{\mddefault}{\updefault}$b$}}}
\put(1726,-1186){\makebox(0,0)[lb]{\smash{\SetFigFont{11}{13.2}{\familydefault}{\mddefault}{\updefault}$c$}}}
\put(1981,-81){\makebox(0,0)[lb]{\smash{\SetFigFont{11}{13.2}{\familydefault}{\mddefault}{\updefault}$c$}}}
\put(1621,529){\makebox(0,0)[lb]{\smash{\SetFigFont{11}{13.2}{\familydefault}{\mddefault}{\updefault}$c$}}}
\put(1276,-1711){\makebox(0,0)[lb]{\smash{\SetFigFont{11}{13.2}{\familydefault}{\mddefault}{\updefault}$\be=0$}}}
\put(2776,-3736){\makebox(0,0)[lb]{\smash{\SetFigFont{11}{13.2}{\familydefault}{\mddefault}{\updefault}hexagonal}}}
\put(-2099,-3886){\makebox(0,0)[lb]{\smash{\SetFigFont{11}{13.2}{\familydefault}{\mddefault}{\updefault}rhombic}}}
\put(-2099,-4111){\makebox(0,0)[lb]{\smash{\SetFigFont{11}{13.2}{\familydefault}{\mddefault}{\updefault}Dodecahedron}}}
\put(601,-2536){\makebox(0,0)[lb]{\smash{\SetFigFont{10}{12.0}{\familydefault}{\mddefault}{\updefault}$\al$}}}
\put(1651,-2536){\makebox(0,0)[lb]{\smash{\SetFigFont{10}{12.0}{\familydefault}{\mddefault}{\updefault}$\g$}}}
\put(2251,-3012){\makebox(0,0)[lb]{\smash{\SetFigFont{9}{10.8}{\familydefault}{\mddefault}{\updefault}$\al$}}}
\put(2911,-3060){\makebox(0,0)[lb]{\smash{\SetFigFont{9}{10.8}{\familydefault}{\mddefault}{\updefault}$0$}}}
\put(3579,-3012){\makebox(0,0)[lb]{\smash{\SetFigFont{9}{10.8}{\familydefault}{\mddefault}{\updefault}$\g$}}}
\put(2506,-2411){\makebox(0,0)[lb]{\smash{\SetFigFont{9}{10.8}{\familydefault}{\mddefault}{\updefault}$a$}}}
\put(2929,-1838){\makebox(0,0)[lb]{\smash{\SetFigFont{9}{10.8}{\familydefault}{\mddefault}{\updefault}$0$}}}
\put(3320,-2411){\makebox(0,0)[lb]{\smash{\SetFigFont{9}{10.8}{\familydefault}{\mddefault}{\updefault}$c$}}}
\put(2858,-2476){\makebox(0,0)[lb]{\smash{\SetFigFont{9}{10.8}{\familydefault}{\mddefault}{\updefault}$b$}}}
\put(3376,-5048){\makebox(0,0)[lb]{\smash{\SetFigFont{8}{9.6}{\familydefault}{\mddefault}{\updefault}$0$}}}
\put(3976,-5092){\makebox(0,0)[lb]{\smash{\SetFigFont{8}{9.6}{\familydefault}{\mddefault}{\updefault}$0$}}}
\put(4583,-5048){\makebox(0,0)[lb]{\smash{\SetFigFont{8}{9.6}{\familydefault}{\mddefault}{\updefault}$\g$}}}
\put(3608,-4502){\makebox(0,0)[lb]{\smash{\SetFigFont{8}{9.6}{\familydefault}{\mddefault}{\updefault}$a$}}}
\put(3992,-3981){\makebox(0,0)[lb]{\smash{\SetFigFont{8}{9.6}{\familydefault}{\mddefault}{\updefault}$0$}}}
\put(4348,-4502){\makebox(0,0)[lb]{\smash{\SetFigFont{8}{9.6}{\familydefault}{\mddefault}{\updefault}$c$}}}
\put(3928,-4561){\makebox(0,0)[lb]{\smash{\SetFigFont{8}{9.6}{\familydefault}{\mddefault}{\updefault}$b$}}}
\put(-2324,-5498){\makebox(0,0)[lb]{\smash{\SetFigFont{8}{9.6}{\familydefault}{\mddefault}{\updefault}$\al$}}}
\put(-1724,-5542){\makebox(0,0)[lb]{\smash{\SetFigFont{8}{9.6}{\familydefault}{\mddefault}{\updefault}$0$}}}
\put(-1117,-5498){\makebox(0,0)[lb]{\smash{\SetFigFont{8}{9.6}{\familydefault}{\mddefault}{\updefault}$\g$}}}
\put(-2092,-4952){\makebox(0,0)[lb]{\smash{\SetFigFont{8}{9.6}{\familydefault}{\mddefault}{\updefault}$a$}}}
\put(-1708,-4431){\makebox(0,0)[lb]{\smash{\SetFigFont{8}{9.6}{\familydefault}{\mddefault}{\updefault}$0$}}}
\put(-1352,-4952){\makebox(0,0)[lb]{\smash{\SetFigFont{8}{9.6}{\familydefault}{\mddefault}{\updefault}$c$}}}
\put(-1772,-5011){\makebox(0,0)[lb]{\smash{\SetFigFont{8}{9.6}{\familydefault}{\mddefault}{\updefault}$0$}}}
\put(2026,-6623){\makebox(0,0)[lb]{\smash{\SetFigFont{8}{9.6}{\familydefault}{\mddefault}{\updefault}$0$}}}
\put(2626,-6667){\makebox(0,0)[lb]{\smash{\SetFigFont{8}{9.6}{\familydefault}{\mddefault}{\updefault}$0$}}}
\put(3233,-6623){\makebox(0,0)[lb]{\smash{\SetFigFont{8}{9.6}{\familydefault}{\mddefault}{\updefault}$\g$}}}
\put(2258,-6077){\makebox(0,0)[lb]{\smash{\SetFigFont{8}{9.6}{\familydefault}{\mddefault}{\updefault}$a$}}}
\put(2642,-5556){\makebox(0,0)[lb]{\smash{\SetFigFont{8}{9.6}{\familydefault}{\mddefault}{\updefault}$0$}}}
\put(2998,-6077){\makebox(0,0)[lb]{\smash{\SetFigFont{8}{9.6}{\familydefault}{\mddefault}{\updefault}$c$}}}
\put(2578,-6136){\makebox(0,0)[lb]{\smash{\SetFigFont{8}{9.6}{\familydefault}{\mddefault}{\updefault}$0$}}}
\put(-974,-4411){\makebox(0,0)[lb]{\smash{\SetFigFont{10}{12.0}{\familydefault}{\mddefault}{\updefault}$\al$}}}
\put(151,-4411){\makebox(0,0)[lb]{\smash{\SetFigFont{10}{12.0}{\familydefault}{\mddefault}{\updefault}$\g$}}}
\put(3001,-4486){\makebox(0,0)[lb]{\smash{\SetFigFont{10}{12.0}{\familydefault}{\mddefault}{\updefault}$\g$}}}
\put(3151,-4936){\makebox(0,0)[lb]{\smash{\SetFigFont{10}{12.0}{\familydefault}{\mddefault}{\updefault}$b$}}}
\put(1801,-5686){\makebox(0,0)[lb]{\smash{\SetFigFont{10}{12.0}{\familydefault}{\mddefault}{\updefault}$\g$}}}
\put(-524,-3961){\makebox(0,0)[lb]{\smash{\SetFigFont{10}{12.0}{\familydefault}{\mddefault}{\updefault}$c$}}}
\put(2326,-4036){\makebox(0,0)[lb]{\smash{\SetFigFont{10}{12.0}{\familydefault}{\mddefault}{\updefault}$c$}}}
\put(1126,-5311){\makebox(0,0)[lb]{\smash{\SetFigFont{10}{12.0}{\familydefault}{\mddefault}{\updefault}$c$}}}
\put(2701,-4036){\makebox(0,0)[lb]{\smash{\SetFigFont{10}{12.0}{\familydefault}{\mddefault}{\updefault}$a$}}}
\put(-149,-3961){\makebox(0,0)[lb]{\smash{\SetFigFont{10}{12.0}{\familydefault}{\mddefault}{\updefault}$a$}}}
\put(1351,-2161){\makebox(0,0)[lb]{\smash{\SetFigFont{10}{12.0}{\familydefault}{\mddefault}{\updefault}$a$}}}
\put(976,-2161){\makebox(0,0)[lb]{\smash{\SetFigFont{10}{12.0}{\familydefault}{\mddefault}{\updefault}$c$}}}
\put(1426,-5311){\makebox(0,0)[lb]{\smash{\SetFigFont{10}{12.0}{\familydefault}{\mddefault}{\updefault}$a$}}}
\put(3151,-3961){\makebox(0,0)[lb]{\smash{\SetFigFont{11}{13.2}{\familydefault}{\mddefault}{\updefault}Prism}}}
\end{picture}

%% file: conbravo.pstex_t
\begin{picture}(0,0)%
\includegraphics{conbravo.pstex}%
\end{picture}%
\setlength{\unitlength}{3355sp}%
\begingroup\makeatletter\ifx\SetFigFont\undefined%
\gdef\SetFigFont#1#2#3#4#5{%
  \reset@font\fontsize{#1}{#2pt}%
  \fontfamily{#3}\fontseries{#4}\fontshape{#5}%
  \selectfont}%
\fi\endgroup%
\begin{picture}(6050,3501)(2101,-2894)
\put(2101,-1036){\makebox(0,0)[lb]{\smash{\SetFigFont{10}{12.0}{\familydefault}{\mddefault}{\updefault}$-\ep$}}}
\put(2951,-1048){\makebox(0,0)[lb]{\smash{\SetFigFont{10}{12.0}{\familydefault}{\mddefault}{\updefault}$-\ep$}}}
\put(3686,-1036){\makebox(0,0)[lb]{\smash{\SetFigFont{10}{12.0}{\familydefault}{\mddefault}{\updefault}$-\ep$}}}
\put(2941,-391){\makebox(0,0)[lb]{\smash{\SetFigFont{10}{12.0}{\familydefault}{\mddefault}{\updefault}$a$}}}
\put(3476,-261){\makebox(0,0)[lb]{\smash{\SetFigFont{10}{12.0}{\familydefault}{\mddefault}{\updefault}$a$}}}
\put(2995,419){\makebox(0,0)[lb]{\smash{\SetFigFont{10}{12.0}{\familydefault}{\mddefault}{\updefault}$0$}}}
\put(2541,-281){\makebox(0,0)[lb]{\smash{\SetFigFont{10}{12.0}{\familydefault}{\mddefault}{\updefault}$a$}}}
\put(5026,419){\makebox(0,0)[lb]{\smash{\SetFigFont{10}{12.0}{\familydefault}{\mddefault}{\updefault}$-\ep$}}}
\put(4426,-211){\makebox(0,0)[lb]{\smash{\SetFigFont{10}{12.0}{\familydefault}{\mddefault}{\updefault}$a\!-\!\ep$}}}
\put(4946,-321){\makebox(0,0)[lb]{\smash{\SetFigFont{10}{12.0}{\familydefault}{\mddefault}{\updefault}$a\!-3\ep$}}}
\put(5575,-251){\makebox(0,0)[lb]{\smash{\SetFigFont{10}{12.0}{\familydefault}{\mddefault}{\updefault}$a\!-\!3\ep$}}}
\put(4209,-1036){\makebox(0,0)[lb]{\smash{\SetFigFont{10}{12.0}{\familydefault}{\mddefault}{\updefault}$2\ep$}}}
\put(5111,-1036){\makebox(0,0)[lb]{\smash{\SetFigFont{10}{12.0}{\familydefault}{\mddefault}{\updefault}$0$}}}
\put(5944,-1036){\makebox(0,0)[lb]{\smash{\SetFigFont{10}{12.0}{\familydefault}{\mddefault}{\updefault}$0$}}}
\put(7309,451){\makebox(0,0)[lb]{\smash{\SetFigFont{10}{12.0}{\familydefault}{\mddefault}{\updefault}$\ep$}}}
\put(6521,-276){\makebox(0,0)[lb]{\smash{\SetFigFont{10}{12.0}{\familydefault}{\mddefault}{\updefault}$a\!-\!3\ep$}}}
\put(7126,-321){\makebox(0,0)[lb]{\smash{\SetFigFont{10}{12.0}{\familydefault}{\mddefault}{\updefault}$a\!-3\ep$}}}
\put(7788,-276){\makebox(0,0)[lb]{\smash{\SetFigFont{10}{12.0}{\familydefault}{\mddefault}{\updefault}$a\!-\!3\ep$}}}
\put(8151,-1057){\makebox(0,0)[lb]{\smash{\SetFigFont{10}{12.0}{\familydefault}{\mddefault}{\updefault}$0$}}}
\put(7323,-1057){\makebox(0,0)[lb]{\smash{\SetFigFont{10}{12.0}{\familydefault}{\mddefault}{\updefault}$0$}}}
\put(6488,-1057){\makebox(0,0)[lb]{\smash{\SetFigFont{10}{12.0}{\familydefault}{\mddefault}{\updefault}$0$}}}
\put(6126,-1365){\makebox(0,0)[lb]{\smash{\SetFigFont{10}{12.0}{\familydefault}{\mddefault}{\updefault}$0$}}}
\put(6571,-2021){\makebox(0,0)[lb]{\smash{\SetFigFont{10}{12.0}{\familydefault}{\mddefault}{\updefault}$a\!-\!4\ep$}}}
\put(6937,-2817){\makebox(0,0)[lb]{\smash{\SetFigFont{10}{12.0}{\familydefault}{\mddefault}{\updefault}$-\ep$}}}
\put(6136,-2836){\makebox(0,0)[lb]{\smash{\SetFigFont{10}{12.0}{\familydefault}{\mddefault}{\updefault}$\ep$}}}
\put(5951,-2086){\makebox(0,0)[lb]{\smash{\SetFigFont{10}{12.0}{\familydefault}{\mddefault}{\updefault}$a\!-2\ep$}}}
\put(5326,-2836){\makebox(0,0)[lb]{\smash{\SetFigFont{10}{12.0}{\familydefault}{\mddefault}{\updefault}$\ep$}}}
\put(4726,-2836){\makebox(0,0)[lb]{\smash{\SetFigFont{10}{12.0}{\familydefault}{\mddefault}{\updefault}$-2\ep$}}}
\put(3881,-2836){\makebox(0,0)[lb]{\smash{\SetFigFont{10}{12.0}{\familydefault}{\mddefault}{\updefault}$-2\ep$}}}
\put(3133,-2828){\makebox(0,0)[lb]{\smash{\SetFigFont{10}{12.0}{\familydefault}{\mddefault}{\updefault}$0$}}}
\put(3301,-2011){\makebox(0,0)[lb]{\smash{\SetFigFont{10}{12.0}{\familydefault}{\mddefault}{\updefault}$a\!+\!\ep$}}}
\put(3831,-2121){\makebox(0,0)[lb]{\smash{\SetFigFont{10}{12.0}{\familydefault}{\mddefault}{\updefault}$a\!-\ep$}}}
\put(4449,-2021){\makebox(0,0)[lb]{\smash{\SetFigFont{10}{12.0}{\familydefault}{\mddefault}{\updefault}$a\!-\!\ep$}}}
\put(3995,-1385){\makebox(0,0)[lb]{\smash{\SetFigFont{10}{12.0}{\familydefault}{\mddefault}{\updefault}$\ep$}}}
\put(5341,-2011){\makebox(0,0)[lb]{\smash{\SetFigFont{10}{12.0}{\familydefault}{\mddefault}{\updefault}$a\!-\!2\ep$}}}
\end{picture}

%% file: figw.pstex_t
\begin{picture}(0,0)%
\includegraphics{figw.pstex}%
\end{picture}%
\setlength{\unitlength}{3355sp}%
\begingroup\makeatletter\ifx\SetFigFont\undefined%
\gdef\SetFigFont#1#2#3#4#5{%
  \reset@font\fontsize{#1}{#2pt}%
  \fontfamily{#3}\fontseries{#4}\fontshape{#5}%
  \selectfont}%
\fi\endgroup%
\begin{picture}(7050,1635)(901,-1711)
\put(7576,-961){\makebox(0,0)[lb]{\smash{\SetFigFont{10}{12.0}{\rmdefault}{\mddefault}{\updefault}I}}}
\put(7126,-211){\makebox(0,0)[lb]{\smash{\SetFigFont{10}{12.0}{\rmdefault}{\mddefault}{\updefault}A}}}
\put(7951,-1711){\makebox(0,0)[lb]{\smash{\SetFigFont{10}{12.0}{\rmdefault}{\mddefault}{\updefault}N}}}
\put(7426,-436){\makebox(0,0)[lb]{\smash{\SetFigFont{10}{12.0}{\familydefault}{\mddefault}{\updefault}$c=0$}}}
\put(7776,-1111){\makebox(0,0)[lb]{\smash{\SetFigFont{10}{12.0}{\familydefault}{\mddefault}{\updefault}$\g=0$}}}
\put(1811,-1561){\makebox(0,0)[lb]{\smash{\SetFigFont{10}{12.0}{\rmdefault}{\mddefault}{\updefault}O}}}
\put(3601,-361){\makebox(0,0)[lb]{\smash{\SetFigFont{10}{12.0}{\rmdefault}{\mddefault}{\updefault}D}}}
\put(2421,-361){\makebox(0,0)[lb]{\smash{\SetFigFont{10}{12.0}{\rmdefault}{\mddefault}{\updefault}B}}}
\put(3031,-1561){\makebox(0,0)[lb]{\smash{\SetFigFont{10}{12.0}{\rmdefault}{\mddefault}{\updefault}J}}}
\put(1201,-361){\makebox(0,0)[lb]{\smash{\SetFigFont{10}{12.0}{\rmdefault}{\mddefault}{\updefault}K}}}
\put(2671,-586){\makebox(0,0)[lb]{\smash{\SetFigFont{10}{12.0}{\familydefault}{\mddefault}{\updefault}$c=0$}}}
\put(2821,-886){\makebox(0,0)[lb]{\smash{\SetFigFont{10}{12.0}{\familydefault}{\mddefault}{\updefault}$\g=0$}}}
\put(3526,-736){\makebox(0,0)[lb]{\smash{\SetFigFont{10}{12.0}{\familydefault}{\mddefault}{\updefault}$\be=0$}}}
\put(1876,-886){\makebox(0,0)[lb]{\smash{\SetFigFont{10}{12.0}{\familydefault}{\mddefault}{\updefault}$b=0$}}}
\put(901,-736){\makebox(0,0)[lb]{\smash{\SetFigFont{10}{12.0}{\familydefault}{\mddefault}{\updefault}$\g=0$}}}
\put(1801,-586){\makebox(0,0)[lb]{\smash{\SetFigFont{10}{12.0}{\familydefault}{\mddefault}{\updefault}$a=0$}}}
\put(1951,-736){\makebox(0,0)[lb]{\smash{\SetFigFont{10}{12.0}{\rmdefault}{\mddefault}{\updefault}or}}}
\put(2851,-736){\makebox(0,0)[lb]{\smash{\SetFigFont{10}{12.0}{\rmdefault}{\mddefault}{\updefault}or}}}
\put(4876,-361){\makebox(0,0)[lb]{\smash{\SetFigFont{10}{12.0}{\rmdefault}{\mddefault}{\updefault}C}}}
\put(5486,-1561){\makebox(0,0)[lb]{\smash{\SetFigFont{10}{12.0}{\rmdefault}{\mddefault}{\updefault}Q}}}
\put(6086,-361){\makebox(0,0)[lb]{\smash{\SetFigFont{10}{12.0}{\rmdefault}{\mddefault}{\updefault}L}}}
\put(4651,-811){\makebox(0,0)[lb]{\smash{\SetFigFont{10}{12.0}{\familydefault}{\mddefault}{\updefault}$c=0$}}}
\put(6001,-811){\makebox(0,0)[lb]{\smash{\SetFigFont{10}{12.0}{\familydefault}{\mddefault}{\updefault}$\g=0$}}}
\end{picture}

%% file: telegraphic.pstex_t
\begin{picture}(0,0)%
\includegraphics{telegraphic.pstex}%
\end{picture}%
\setlength{\unitlength}{2960sp}%
\begingroup\makeatletter\ifx\SetFigFont\undefined%
\gdef\SetFigFont#1#2#3#4#5{%
  \reset@font\fontsize{#1}{#2pt}%
  \fontfamily{#3}\fontseries{#4}\fontshape{#5}%
  \selectfont}%
\fi\endgroup%
\begin{picture}(7587,6837)(-74,-6523)
\put(2761,-1316){\makebox(0,0)[lb]{\smash{\SetFigFont{17}{20.4}{\rmdefault}{\mddefault}{\updefault}?}}}
\put(901,-811){\makebox(0,0)[lb]{\smash{\SetFigFont{11}{13.2}{\rmdefault}{\mddefault}{\updefault}no}}}
\put(1931,-1861){\makebox(0,0)[lb]{\smash{\SetFigFont{11}{13.2}{\rmdefault}{\mddefault}{\updefault}yes}}}
\put(2781,-4541){\makebox(0,0)[lb]{\smash{\SetFigFont{17}{20.4}{\rmdefault}{\mddefault}{\updefault}?}}}
\put(1801,-6136){\makebox(0,0)[lb]{\smash{\SetFigFont{11}{13.2}{\rmdefault}{\mddefault}{\updefault}of period 2.}}}
\put(1651,-5876){\makebox(0,0)[lb]{\smash{\SetFigFont{11}{13.2}{\rmdefault}{\mddefault}{\updefault}orthogonal and }}}
\put(1726,-5611){\makebox(0,0)[lb]{\smash{\SetFigFont{11}{13.2}{\rmdefault}{\mddefault}{\updefault}then they are}}}
\put(4101,-6406){\makebox(0,0)[lb]{\smash{\SetFigFont{11}{13.2}{\rmdefault}{\mddefault}{\updefault}are orthogonal.}}}
\put(4101,-6211){\makebox(0,0)[lb]{\smash{\SetFigFont{11}{13.2}{\rmdefault}{\mddefault}{\updefault}and their glide mirrors}}}
\put(4101,-5986){\makebox(0,0)[lb]{\smash{\SetFigFont{11}{13.2}{\rmdefault}{\mddefault}{\updefault}of one is horizontal, }}}
\put(5401,-5086){\makebox(0,0)[lb]{\smash{\SetFigFont{11}{13.2}{\rmdefault}{\mddefault}{\updefault}no}}}
\put(2551,-5086){\makebox(0,0)[lb]{\smash{\SetFigFont{11}{13.2}{\rmdefault}{\mddefault}{\updefault}no}}}
\put(2026,-4541){\makebox(0,0)[lb]{\smash{\SetFigFont{11}{13.2}{\rmdefault}{\mddefault}{\updefault}parallel}}}
\put(2176,-4291){\makebox(0,0)[lb]{\smash{\SetFigFont{11}{13.2}{\rmdefault}{\mddefault}{\updefault}all}}}
\put(1951,-4036){\makebox(0,0)[lb]{\smash{\SetFigFont{11}{13.2}{\rmdefault}{\mddefault}{\updefault}are they }}}
\put(1126,-4036){\makebox(0,0)[lb]{\smash{\SetFigFont{11}{13.2}{\rmdefault}{\mddefault}{\updefault}yes}}}
\put(4651,-4036){\makebox(0,0)[lb]{\smash{\SetFigFont{11}{13.2}{\rmdefault}{\mddefault}{\updefault}are all their }}}
\put(4421,-2971){\makebox(0,0)[lb]{\smash{\SetFigFont{11}{13.2}{\rmdefault}{\mddefault}{\updefault}a screw motion.}}}
\put(4281,-2726){\makebox(0,0)[lb]{\smash{\SetFigFont{11}{13.2}{\rmdefault}{\mddefault}{\updefault}is a translation or}}}
\put(1951,-1021){\makebox(0,0)[lb]{\smash{\SetFigFont{11}{13.2}{\rmdefault}{\mddefault}{\updefault}  screw-}}}
\put(1876,-736){\makebox(0,0)[lb]{\smash{\SetFigFont{11}{13.2}{\rmdefault}{\mddefault}{\updefault}Are there }}}
\put(1951,-1316){\makebox(0,0)[lb]{\smash{\SetFigFont{11}{13.2}{\rmdefault}{\mddefault}{\updefault}motions}}}
\put(3676,-811){\makebox(0,0)[lb]{\smash{\SetFigFont{11}{13.2}{\rmdefault}{\mddefault}{\updefault}no}}}
\put(5421,-1861){\makebox(0,0)[lb]{\smash{\SetFigFont{11}{13.2}{\rmdefault}{\mddefault}{\updefault}yes}}}
\put(4701,-1291){\makebox(0,0)[lb]{\smash{\SetFigFont{11}{13.2}{\rmdefault}{\mddefault}{\updefault}reflections}}}
\put(4901,-1011){\makebox(0,0)[lb]{\smash{\SetFigFont{11}{13.2}{\rmdefault}{\mddefault}{\updefault}glide-}}}
\put(4726,-736){\makebox(0,0)[lb]{\smash{\SetFigFont{11}{13.2}{\rmdefault}{\mddefault}{\updefault}Are there}}}
\put(4726,164){\makebox(0,0)[lb]{\smash{\SetFigFont{11}{13.2}{\rmdefault}{\mddefault}{\updefault}Start Here}}}
\put(6451,-4036){\makebox(0,0)[lb]{\smash{\SetFigFont{11}{13.2}{\rmdefault}{\mddefault}{\updefault}yes}}}
\put(5736,-1276){\makebox(0,0)[lb]{\smash{\SetFigFont{17}{20.4}{\rmdefault}{\mddefault}{\updefault}?}}}
\put(4021,-5536){\makebox(0,0)[lb]{\smash{\SetFigFont{11}{13.2}{\rmdefault}{\mddefault}{\updefault}then they are in vertical}}}
\put(1726,-2746){\makebox(0,0)[lb]{\smash{\SetFigFont{11}{13.2}{\rmdefault}{\mddefault}{\updefault}decomposes!}}}
\put(1621,-2461){\makebox(0,0)[lb]{\smash{\SetFigFont{11}{13.2}{\rmdefault}{\mddefault}{\updefault}then the lattice}}}
\put(4261,-2461){\makebox(0,0)[lb]{\smash{\SetFigFont{11}{13.2}{\rmdefault}{\mddefault}{\updefault}then their product}}}
\put(4051,-5761){\makebox(0,0)[lb]{\smash{\SetFigFont{11}{13.2}{\rmdefault}{\mddefault}{\updefault}planes, the glide vector}}}
\put(4581,-4261){\makebox(0,0)[lb]{\smash{\SetFigFont{11}{13.2}{\rmdefault}{\mddefault}{\updefault}glide mirrors}}}
\put(4851,-4561){\makebox(0,0)[lb]{\smash{\SetFigFont{11}{13.2}{\rmdefault}{\mddefault}{\updefault}parallel}}}
\put(5601,-4561){\makebox(0,0)[lb]{\smash{\SetFigFont{17}{20.4}{\rmdefault}{\mddefault}{\updefault}?}}}
\put(-74,-4036){\makebox(0,0)[lb]{\smash{\SetFigFont{14}{16.8}{\rmdefault}{\mddefault}{\updefault}$c4$}}}
\put(-74,-4936){\makebox(0,0)[lb]{\smash{\SetFigFont{14}{16.8}{\rmdefault}{\mddefault}{\updefault}$c6$}}}
\put(7426,-6361){\makebox(0,0)[lb]{\smash{\SetFigFont{14}{16.8}{\rmdefault}{\mddefault}{\updefault}$-a2$}}}
\put(7426,-5611){\makebox(0,0)[lb]{\smash{\SetFigFont{14}{16.8}{\rmdefault}{\mddefault}{\updefault}$+a2$}}}
\put(7426,-4636){\makebox(0,0)[lb]{\smash{\SetFigFont{14}{16.8}{\rmdefault}{\mddefault}{\updefault}$-a1$}}}
\put(7426,-3886){\makebox(0,0)[lb]{\smash{\SetFigFont{14}{16.8}{\rmdefault}{\mddefault}{\updefault}$+a1$}}}
\put(-74,-1036){\makebox(0,0)[lb]{\smash{\SetFigFont{14}{16.8}{\rmdefault}{\mddefault}{\updefault}$c1$}}}
\put(-74,-2236){\makebox(0,0)[lb]{\smash{\SetFigFont{14}{16.8}{\rmdefault}{\mddefault}{\updefault}$c2$}}}
\put(-74,-3136){\makebox(0,0)[lb]{\smash{\SetFigFont{14}{16.8}{\rmdefault}{\mddefault}{\updefault}$c3$}}}
\put(-74,-5836){\makebox(0,0)[lb]{\smash{\SetFigFont{14}{16.8}{\rmdefault}{\mddefault}{\updefault}$c22$}}}
\end{picture}

%% file: fig1.pstex_t
\begin{picture}(0,0)%
\includegraphics{fig1.pstex}%
\end{picture}%
\setlength{\unitlength}{3355sp}%
\begingroup\makeatletter\ifx\SetFigFont\undefined%
\gdef\SetFigFont#1#2#3#4#5{%
  \reset@font\fontsize{#1}{#2pt}%
  \fontfamily{#3}\fontseries{#4}\fontshape{#5}%
  \selectfont}%
\fi\endgroup%
\begin{picture}(1797,1792)(1876,-1877)
\put(1876,-1814){\makebox(0,0)[lb]{\smash{\SetFigFont{11}{13.2}{\rmdefault}{\mddefault}{\updefault}$p_{01}$}}}
\put(3673,-1814){\makebox(0,0)[lb]{\smash{\SetFigFont{11}{13.2}{\rmdefault}{\mddefault}{\updefault}$p_{03}$}}}
\put(2774,-1814){\makebox(0,0)[lb]{\smash{\SetFigFont{11}{13.2}{\rmdefault}{\mddefault}{\updefault}$p_{13}$}}}
\put(2774,-253){\makebox(0,0)[lb]{\smash{\SetFigFont{11}{13.2}{\rmdefault}{\mddefault}{\updefault}$p_{02}$}}}
\put(2161,-963){\makebox(0,0)[lb]{\smash{\SetFigFont{11}{13.2}{\rmdefault}{\mddefault}{\updefault}$p_{12}$}}}
\put(3341,-963){\makebox(0,0)[lb]{\smash{\SetFigFont{11}{13.2}{\rmdefault}{\mddefault}{\updefault}$p_{23}$}}}
\put(2754,-1154){\makebox(0,0)[lb]{\smash{\SetFigFont{11}{13.2}{\rmdefault}{\mddefault}{\updefault}$0$}}}
\end{picture}

%% file: fig0.pstex_t
\begin{picture}(0,0)%
\includegraphics{fig0.pstex}%
\end{picture}%
\setlength{\unitlength}{3355sp}%
\begingroup\makeatletter\ifx\SetFigFont\undefined%
\gdef\SetFigFont#1#2#3#4#5{%
  \reset@font\fontsize{#1}{#2pt}%
  \fontfamily{#3}\fontseries{#4}\fontshape{#5}%
  \selectfont}%
\fi\endgroup%
\begin{picture}(5956,3427)(2195,-2868)
\put(2195,-1048){\makebox(0,0)[lb]{\smash{\SetFigFont{8}{9.6}{\rmdefault}{\mddefault}{\updefault}$4$}}}
\put(3861,-1048){\makebox(0,0)[lb]{\smash{\SetFigFont{8}{9.6}{\rmdefault}{\mddefault}{\updefault}$6$}}}
\put(2985,403){\makebox(0,0)[lb]{\smash{\SetFigFont{8}{9.6}{\rmdefault}{\mddefault}{\updefault}$5$}}}
\put(2951,-395){\makebox(0,0)[lb]{\smash{\SetFigFont{8}{9.6}{\rmdefault}{\mddefault}{\updefault}$0$}}}
\put(2966,-1048){\makebox(0,0)[lb]{\smash{\SetFigFont{8}{9.6}{\rmdefault}{\mddefault}{\updefault}$-1$}}}
\put(2463,-251){\makebox(0,0)[lb]{\smash{\SetFigFont{8}{9.6}{\rmdefault}{\mddefault}{\updefault}$-1$}}}
\put(3468,-251){\makebox(0,0)[lb]{\smash{\SetFigFont{8}{9.6}{\rmdefault}{\mddefault}{\updefault}$-1$}}}
\put(2391,324){\makebox(0,0)[lb]{\smash{\SetFigFont{9}{10.8}{\rmdefault}{\mddefault}{\updefault}(a)}}}
\put(6126,-1365){\makebox(0,0)[lb]{\smash{\SetFigFont{8}{9.6}{\rmdefault}{\mddefault}{\updefault}$3$}}}
\put(6937,-2817){\makebox(0,0)[lb]{\smash{\SetFigFont{8}{9.6}{\rmdefault}{\mddefault}{\updefault}$2$}}}
\put(5259,-2817){\makebox(0,0)[lb]{\smash{\SetFigFont{8}{9.6}{\rmdefault}{\mddefault}{\updefault}$2$}}}
\put(6054,-2817){\makebox(0,0)[lb]{\smash{\SetFigFont{8}{9.6}{\rmdefault}{\mddefault}{\updefault}$1$}}}
\put(6049,-2118){\makebox(0,0)[lb]{\smash{\SetFigFont{8}{9.6}{\rmdefault}{\mddefault}{\updefault}$0$}}}
\put(5546,-2046){\makebox(0,0)[lb]{\smash{\SetFigFont{8}{9.6}{\rmdefault}{\mddefault}{\updefault}$-1$}}}
\put(6571,-2046){\makebox(0,0)[lb]{\smash{\SetFigFont{8}{9.6}{\rmdefault}{\mddefault}{\updefault}$1$}}}
\put(5944,-1048){\makebox(0,0)[lb]{\smash{\SetFigFont{8}{9.6}{\rmdefault}{\mddefault}{\updefault}$3$}}}
\put(5575,-251){\makebox(0,0)[lb]{\smash{\SetFigFont{8}{9.6}{\rmdefault}{\mddefault}{\updefault}$0$}}}
\put(4209,-1048){\makebox(0,0)[lb]{\smash{\SetFigFont{8}{9.6}{\rmdefault}{\mddefault}{\updefault}$1$}}}
\put(5031,-1048){\makebox(0,0)[lb]{\smash{\SetFigFont{8}{9.6}{\rmdefault}{\mddefault}{\updefault}$2$}}}
\put(5109,395){\makebox(0,0)[lb]{\smash{\SetFigFont{8}{9.6}{\rmdefault}{\mddefault}{\updefault}$4$}}}
\put(4956,-371){\makebox(0,0)[lb]{\smash{\SetFigFont{8}{9.6}{\rmdefault}{\mddefault}{\updefault}$-1$}}}
\put(4616,-251){\makebox(0,0)[lb]{\smash{\SetFigFont{8}{9.6}{\rmdefault}{\mddefault}{\updefault}$0$}}}
\put(4381,324){\makebox(0,0)[lb]{\smash{\SetFigFont{9}{10.8}{\rmdefault}{\mddefault}{\updefault}(c)}}}
\put(5176,-2386){\makebox(0,0)[lb]{\smash{\SetFigFont{9}{10.8}{\rmdefault}{\mddefault}{\updefault}(d)}}}
\put(4823,-2828){\makebox(0,0)[lb]{\smash{\SetFigFont{8}{9.6}{\rmdefault}{\mddefault}{\updefault}$5$}}}
\put(3133,-2828){\makebox(0,0)[lb]{\smash{\SetFigFont{8}{9.6}{\rmdefault}{\mddefault}{\updefault}$3$}}}
\put(3995,-2828){\makebox(0,0)[lb]{\smash{\SetFigFont{8}{9.6}{\rmdefault}{\mddefault}{\updefault}$0$}}}
\put(4449,-2031){\makebox(0,0)[lb]{\smash{\SetFigFont{8}{9.6}{\rmdefault}{\mddefault}{\updefault}$-2$}}}
\put(3451,-2031){\makebox(0,0)[lb]{\smash{\SetFigFont{8}{9.6}{\rmdefault}{\mddefault}{\updefault}$-2$}}}
\put(3923,-2175){\makebox(0,0)[lb]{\smash{\SetFigFont{8}{9.6}{\rmdefault}{\mddefault}{\updefault}$1$}}}
\put(3995,-1385){\makebox(0,0)[lb]{\smash{\SetFigFont{8}{9.6}{\rmdefault}{\mddefault}{\updefault}$6$}}}
\put(2926,-2441){\makebox(0,0)[lb]{\smash{\SetFigFont{9}{10.8}{\rmdefault}{\mddefault}{\updefault}(b)}}}
\put(7323,-1057){\makebox(0,0)[lb]{\smash{\SetFigFont{8}{9.6}{\rmdefault}{\mddefault}{\updefault}$0$}}}
\put(6488,-1057){\makebox(0,0)[lb]{\smash{\SetFigFont{8}{9.6}{\rmdefault}{\mddefault}{\updefault}$1$}}}
\put(8151,-1057){\makebox(0,0)[lb]{\smash{\SetFigFont{8}{9.6}{\rmdefault}{\mddefault}{\updefault}$3$}}}
\put(7309,451){\makebox(0,0)[lb]{\smash{\SetFigFont{8}{9.6}{\rmdefault}{\mddefault}{\updefault}$2$}}}
\put(7228,-349){\makebox(0,0)[lb]{\smash{\SetFigFont{8}{9.6}{\rmdefault}{\mddefault}{\updefault}$1$}}}
\put(7788,-276){\makebox(0,0)[lb]{\smash{\SetFigFont{8}{9.6}{\rmdefault}{\mddefault}{\updefault}$0$}}}
\put(6807,-276){\makebox(0,0)[lb]{\smash{\SetFigFont{8}{9.6}{\rmdefault}{\mddefault}{\updefault}$0$}}}
\put(6653,305){\makebox(0,0)[lb]{\smash{\SetFigFont{9}{10.8}{\rmdefault}{\mddefault}{\updefault}(e)}}}
\end{picture}

%% file: dp.bbl
\begin{thebibliography}{URLW}

\bibitem[Ch]{Ch} {\sc Charlap, L.S.}, Bieberbach groups and flat manifolds.
{\it Universitext. Springer-Verlag, New York}, 1986.

\bibitem[ChV]{ChV} {\sc Charlap, L. S.; Vasquez, A. T.}, Compact flat riemannian manifolds. III. The group of affinities.
{\it Amer. J. Math.} {\bf 95} (1973), 471--494.

\bibitem[CiS]{CS} {\sc Cid, C.; Schulz, T.}, Computation of Five and Six dimensional Bieberbach groups.
{\it Experiment Math.} {\bf 10} (2001), 109--115.

\bibitem[Co1]{orbnot} {\sc Conway, J.H.},
The orbifold notation for surface groups.
In {\it Groups, combinatorics and geometry} (Durham, 1990), 438--447,
London Math. Soc. Lecture Note Ser., 165,
Cambridge Univ. Press, Cambridge, 1992.

\bibitem[Co2]{sensualform} {\sc Conway, J.H.} (with the assistance of {\sc F.Y.C. Fung}),
The sensual (quadratic) form,
{\it Carus Mathematical Monographs, {\bf 26}
Mathematical Association of America, Washington, DC} (1997).

\bibitem[CDHT]{CDHT} {\sc Conway, J.H.; Delgado Friedrichs, O.; Huson, D.H.; Thurston, W.P.},
On three-dimensional space groups.
{\it Beitr\"age Algebra Geom.} {\bf 42} (2001), no. 2, 475--507.

\bibitem[CoS]{CSVI} {\sc Conway, J.H.; Sloane, N.J.A.},
Low-dimensional lattices. VI. Vorono\u{\i} reduction of three-dimensional lattices,
{\it Proc. Roy. Soc. London Ser. A} {\bf 436} (1992), no. 1896, 55--68.

\bibitem[DR]{DR} {\sc Doyle, P.G.; Rossetti, J.P.}, Tetra and Didi, the cosmic spectral twins.
{\it Preprint}.

\bibitem[E]{E}{\sc Efstathiou, G.},
Is the Low CMB Quadrupole a Signature of Spatial Curvature?
Submitted to {\it MNRAS}, arXiv:astro-ph/0303127.

\bibitem[HS]{HS} {\sc Hajian, A.; Souradeep, T.},
Statistical Isotropy of CMB and Cosmic Topology.
Submitted to {\it Phys. Rev. Lett.}, arXiv: astro-ph/0301590.

\bibitem[HW]{HW} {\sc Hantzsche W.; Wendt H.},
Dreidimensionale euklidische Raumformen.
{\it Math. Annalen.}  {\bf 110} (1934--35), 593--611.

\bibitem[Hi]{Hi} {\sc Hillman, J.A.}, Flat $4$-manifold groups.
{\it New Zealand J. Math.} {\bf 24} (1995), no. 1, 29--40.

\bibitem[LSW]{sciamer} {\sc Luminet, J.P.; Starkman, G.D.; Weeks, J.R.},
Is Space Finite? {\it Scientific American}  (1999), no. April, 90--97.

\bibitem[MS]{MS} {\sc Malfait, W.; Szczepa\'nski, A.},
The structure of the (outer) automorphism group of a Bieberbach group.
{\it Compositio Math.} {\bf 136} (2003), no. 1, 89--101.

\bibitem[No]{No} {\sc Nowacki, W.}, Die euklidischen, dreidimensionalen, geschlossenen und offenen Raumformen.
{\it Comm. Math. Helv.}  {\bf 7} (1934), 81--93.

\bibitem[NYT]{NYT} {\sc Dennis Overbye}, Universe as Doughnut: New Data, New Debate.
{\it New York Times} March 11, 2003.

\bibitem[RC]{hearing} {\sc Rossetti, J.P.; Conway, J.H.}, Hearing the Platycosms. {\it Preprint}.

\bibitem[Spe]{Spe} {\sc Spergel, D.N.; et al.},
First Year Wilkinson Microwave Anisotropy Probe (WMAP) Observations: Determination of Cosmological Parameters.
To appear in  {\it The Astrophysical Journal}, astro-ph/0302209.

\bibitem[TOH]{TOH} {\sc Tegmark, M.; de Oliveira-Costa, A.; Hamilton, A.},
A high resolution foreground cleaned CMB map from WMAP.
arXiv:astro-ph/0302496.

\bibitem[Th]{Th} {\sc Thurston, W.P.}, Three-dimensional geometry and topology. Vol. 1.
{\it Edited by Silvio Levy. Princeton Mathematical Series, 35. Princeton University Press, Princeton, NJ}, 1997.

\bibitem[URLW]{URLW}{\sc Uzan, J.-P.; Riazuelo, A.; Lehoucq, R; Weeks, J.},
Cosmic microwave background constraints on multi-connected spherical spaces.
Submitted to {\it Phys.Rev.Lett.}, arXiv:astro-ph/0303580.

\bibitem[We1]{We}{\sc Weeks, J.R.} The shape of space.
{\it Second edition. Monographs and Textbooks in Pure and Applied Mathematics, 249. Marcel Dekker, Inc., New York}, 2002.

\bibitem[We2]{We2}{\sc Weeks, J.R.} Real-Time Rendering in Curved Spaces.
{\it IEEE Computer Graphics and Applications} {\bf 22} (2002), no 6, 90--99.

\bibitem[Wo]{Wo} {\sc Wolf, J.A.}, Spaces of constant curvature,
{\it Fifth edition. Publish or Perish, Inc., Houston, TX}, 1984.

\bibitem[Zi]{Zim} {\sc Zimmermann, B.}, On the Hantzsche-Wendt manifold.
{\it Monatsh. Math.} {\bf 110} (1990), no. 3-4, 321--327.


\end{thebibliography}
